\documentclass[12pt]{amsart}
\usepackage{amssymb,amsfonts,amsmath,amsopn,amstext,amscd,color,latexsym,mathrsfs,skak,stackrel,url,verbatim}
\usepackage{mathabx}
\usepackage[utf8]{inputenc}
\usepackage[all, cmtip]{xy}
\usepackage{rotating}
\usepackage[active]{srcltx}
\usepackage[pdftex, colorlinks=true, linktoc=all]{hyperref}

\usepackage{cleveref}

\numberwithin{equation}{subsection}

\newtheorem{proposition}[equation]{Proposition}

\theoremstyle{definition}

\theoremstyle{remark}

\theoremstyle{observation}

\newcommand{\hide}[1]{}

\newcommand{\isom}{\cong}

\usepackage{amsmath,epsf}
\usepackage{amsfonts}
\usepackage{enumerate}
\usepackage{graphicx,version}

\begin{document}

\title[Addendum/Erratum to ``Weyl Groups and Birational Transformations'']{Addendum/Erratum to \linebreak ``Weyl Groups and Birational Transformations among Minimal Models''}

\author{Kenji Matsuki}
\address{Department of Mathematics, Purdue University,
150 N. University Street, West Lafayette, IN 47907, U.S.A.}
\email{matsuki@purdue.edu}

\begin{abstract}
We present an addendum/erratum to the paper ``Weyl Groups and Birational Transformations among Minimal Models'' \cite{WGBTMM95} written by the author and published in 1995, adding the analysis of the ``88-th'' deformation type of a smooth Fano 3-fold with $B_2 = 4$ denoted as $n^o\ 13$, which was missing from the original classification table by Mori-Mukai \cite{MMManu81}\cite{MMAdv83} and added later \cite{MMManu03} to the list of smooth Fano 3-folds with $B_2 \geq 2$, while correcting the mistake pointed out by Dr. Eric Jovinelly (also noticed earlier by Kento Fujita).  We also correct other typos and miscalculations, clarifying some points of ambiguity.
\end{abstract}
\maketitle

\tableofcontents

\section{Outline of the addendum/erratum}

The paper ``Weyl Groups and Birational Transformations among Minimal Models'' \cite{WGBTMM95} claims to provide a complete and explicit list of all the extremal rays for each deformation type of smooth Fano 3-folds with $B_2 \geq 2$, classified by Mori-Mukai \cite{MMManu81, MMAdv83}.  Recently, however, Dr. Eric Jovinelly kindly informed the author that one extremal ray was missing from the list for the Fano 3-fold denoted as $n^o\ 1$  with $B_2 = 5$ in \cite{WGBTMM95}  (also denoted as $n^o\ 1$ in Table 5 for Fano 3-folds with $B_2 \geq 5$ in \cite{MMManu81}).  Jovinelly later informed me that the missing extremal ray was also noticed earlier by K. Fujita \cite{Fujita15}.  In \cite{WGBTMM95} we use Claim III-3-2 as a criterion to judge if the given set of extremal rays should exhaust all of them.  Therefore, finding a missing extremal ray might have raised a question about the validity of the criterion itself.  Fortunately, Jovinelly confirms the validity of the criterion Claim III-3-2, even though he indicates that some of the statements in the proof of Claim III-3-2 need more explanation/clarification.  The missing extremal ray was simply the result of an oversight of the author.

\vskip.1in

The ouline of this addendum/erratum goes as follows.

\vskip.1in

In section 2, we revisit Claim III-3-2 with some detailed explanation/clarification of the statements used in the original proof.

In section 3, we present an analysis of the Fano 3-fold denoted as $n^o\ 1$  with $B_2 = 5$ in \cite{WGBTMM95}  (also denoted as $n^o\ 1$ in Table 5 for Fano 3-folds with $B_2 \geq 5$ in \cite{MMManu81}) and the description of the missing extremal ray.  We mention the reason of the oversight of the author (why the autor missed this extremal ray), while applying the criterion described in Claim III-3-2.

In section 4, we present an analysis of the Fano 3-fold denoted as $n^o\ 13$ in Table 4 for Fano 3-folds with $B_2 = 4$ in \cite{MMManu03}, which was missing from the original classification table in \cite{MMManu81}.

In section 5, we clarify some misleading statements about the classification of the types of the 4-fold flops connecting all the minimal models which appear in our consideration of the Weyl groups associated with the smooth Fano 3-folds.

In section 6, we list the mistakes found in the tables of the intersection pairings, and present the corrections.

\vskip.1in

{\bf Acknowledgement}: Biggest thanks go to Dr. Eric Jovinelly, who not only brought the attention of the author to the missing extremal ray in \cite{WGBTMM95} but also showed the author his detailed notes on the analysis of the Fano 3-fold denoted as $n^o\ 1$ with $B_2 = 5$ (even though the analysis presented here is different from his notes).  The author would like to thank Profs. Masaki Kashiwara, Shigefumi Mori, and Shigeru Mukai for their consistent guidance and warm support.  This note was written while the author was staying at KUIAS \footnote{KUIAS: Kyoto University Institute for Advanced Study} in the summers of 2022 and 2023.

\section{Claim III-3-2 revisited}

Let $T$ be a smooth Fano 3-fold.  That is to say, $T$ is a nonsingular projective variety defined over $\mathbb{C}$ of $\dim T = 3$, satisfying the condition that $- K_T$ is ample.

\subsection{Claim III-3-2} The following proposition, which is Claim III-3-2 of  \cite{WGBTMM95}, provides not only a criterion to judge if the given set of extremal rays should exhaust all of them, but also an inductive procedure to find all the extremal rays starting from a given set of extremal rays.

\begin{proposition} Let $L = \{l_1, l_2, \ldots, l_m\}$ be a given set of extremal rays on a smooth Fano 3-fold $T$.  Take and fix an extremal ray $l_i \in L$.  Then the contraction of the extremal ray $l_i$, denoted by $\mathrm{cont}_{l_i}: T \rightarrow U_i$, gives another $\mathbb{Q}$-factorial variety $U_i$.  (We note that $\mathrm{cont}_{l_i}$ cannot be of flipping type by \cite{Mori82}, and hence that it is either of fiber type or of divisorial type.  We see then that the normal variety $U_i$ is $\mathbb{Q}$-factorial by Lemma 5-1-5 or Proposition 5-1-6 in \cite{KMM}, respectively.)

Note that $\overline{\mathrm{NE}}(T)$ is a polyhedral cone by the Cone Theorem (cf. Theorem 4-2-1 in \cite{KMM}), since $T$ is a smooth Fano 3-fold.  Let $\phi_i = (\mathrm{cont}_{l_i})_*:N_1(T) \rightarrow N_1(U_i)$ be the natural linear map.  Then $\overline{\mathrm{NE}}(U_i) = \phi_i(\overline{\mathrm{NE}}(T))$ is also a polyhedral cone, being the image of another polyhedral cone $\overline{\mathrm{NE}}(T)$ under the linear map $\phi_i$.  Let $M_i = \{l_{1i}, l_{2i}, \ldots, l_{ni}\}$ denote the set of the edges of the polyhedral cone $\overline{\mathrm{NE}}(U_i)$.  We call $M_i = \{l_{1i}, l_{2i}, \ldots, l_{ni}\}$ the set of the the extremal rays of the polyhedral cone $\overline{\mathrm{NE}}(U_i)$ by abuse of language, even though some edge $l_{ki}$ may not have negative intersection \text{with $K_{U_i}$.}

Observe that each $l_{ki} \in M_i$ is the image of some extremal ray $l_{\alpha}$ of $\overline{\mathrm{NE}}(T)$.  Remark that the linear map $\phi_i$ crushes the extremal ray $l_i$ to a point, and that the fiber of $\phi_i$ is 1-dimensional, since $\dim N_1(U_i) = \rho(U_i) = \rho(T) - 1 = \dim N_1(T) - 1$.  Therefore, we conclude that $l_{\alpha} (\neq l_i)$ is the unique extremal ray (of $\overline{\mathrm{NE}}(T)$) which maps to $l_{ki}\in M$, and that $l_i$ and $l_{\alpha}$ spans a 2-dimensional face $F_{i\alpha}$ of $\overline{\mathrm{NE}}(T)$.

\begin{enumerate}

\item ({\bf Criterion for exhaustion}) The set $L$ exhausts all the possible extremal rays on $T$ if and only if the following is satisfied:

For each $l_i \in L$ and $l_{ki} \in M_i$ there exists another extremal ray $l_j \in L$ such that $\phi_i(l_j) = l_{ki}$.

(In this case it follows immediately that $\phi_j(l_i) \in M_j$ is an extremal ray of $\overline{\mathrm{NE}}(U_j)$ and that we have $\mathrm{cont}_{l_{ki}}(U_i) = \mathrm{cont}_{\phi_j(l_i)}(U_j)$, which is obtained as the contraction $\mathrm{cont}_{F_{ij}}(T) = \mathrm{cont}_{F_{ji}}(T)$ of the face $F_{ij} = F_{ji}$ of $\overline{\mathrm{NE}}(T)$ \text{(cf. Theorem 3-2-1 in \cite{KMM}).}

\item ({\bf Inductive procedure to find all the extremal rays}) Suppose that the given set $L$ does not exhaust all the possible extremal rays on $T$.  Then according to the criterion above, for some $l_i \in L$ and $l_{ki} \in M_i$, there is no $l_j \in L$ which maps onto $l_{ki}$ via $\phi_i$.   We note that, since $\dim N_1(U_i) = \rho(U_i) = \rho(T) - 1 = \dim N_1(T) - 1$, we may assume we know the set of all the extremal rays $M_i = \{l_{1i}, l_{2i}, \ldots, l_{ni}\}$ of $U_i$ inductively.  (We note that $U_i$ may not be a smooth Fano 3-fold and that $\dim U_i$ may even be smaller than $3$.)  Now, as observed before, $l_{ki} \in M_i$ is the image of some extremal ray $l_{\alpha}$ of $\overline{\mathrm{NE}}(T)$.  Add this new extremal ray $l_{\alpha} = l_{m+1}$ to $L$, and repeat the same procedure.  Since there are only finitely many extremal rays of $\overline{\mathrm{NE}}(T)$, this process must come to an end after finitely many repetitions, providing the complete set of all the extremal rays on $T$, satisfying the criterion given in (i).

\end{enumerate}

\end{proposition}

\subsection{Proof of Claim III-3-2} We consider the content of the proposition above in the dual space $N^1(T)$.

The nef cone $\overline{\mathrm{Amp}}(T)$ is the polyhedral cone dual to $\overline{\mathrm{NE}}(T)$.  Its boundary consists of the codimension-one faces $\overline{\mathrm{Amp}}(T) \cap l^{\perp}$ where $l$ varies \text{among all the extremal rays on $T$.}  Therefore, the given set $L = \{l_1, l_2, \cdots, l_m\}$ of extremal rays exhausts all the possible extremal rays if and only if the codimension-one faces $\overline{\mathrm{Amp}}(T) \cap l_i^{\perp}$, which is isomorphic to $\overline{\mathrm{Amp}}(U_i)$, patch up together along the codimension-two faces so that they form the boundary of $\overline{\mathrm{Amp}}(T)$ with no broken part.  That is to say, for each condimesion-one face $\overline{\mathrm{Amp}}(T) \cap l_i^{\perp} \isom \overline{\mathrm{Amp}}(U_i)$ and its codimension-two face $\overline{\mathrm{Amp}}(U_i) \cap l_{ki}^{\perp}$, there exists another codimension-one face $\overline{\mathrm{Amp}}(T) \cap l_j^{\perp} \isom \overline{\mathrm{Amp}}(U_j)$ such that they share the same codimension-two face $\overline{\mathrm{Amp}}(U_i) \cap l_{ki}^{\perp} = \overline{\mathrm{Amp}}(U_j) \cap \phi_j(l_i)^{\perp}$, with $\phi_j(l_i)$ being an extremal ray on $U_j$.  Switching the roles of $i$ and $j$, we have $\overline{\mathrm{Amp}}(U_j) \cap \phi_j(l_i)^{\perp} = \overline{\mathrm{Amp}}(U_i) \cap \phi_i(l_j)^{\perp}$.  Such an extremal ray $l_j$ is characterized as the one with $\phi_i(l_j) = l_{ki}$.  This is exactly the criterion (i).  Assertion (ii) easily follows from (i).

\section{Analysis of the Fano 3-fold denoted as $n^o 1$ with $B_2 \geq 5$ in Table 5}

In this section, we present an analysis of the Fano 3-fold $T$ denoted as $n^o\ 1$  with $B_2 = 5$ in \cite{WGBTMM95}  (also denoted as $n^o\ 1$ in Table 5 for Fano 3-folds with $B_2 \geq 5$ in \cite{MMManu81}).  

\vskip.1in

$\boxed{B_2 = 5}$ : $n^o\ 1$

\vskip.1in

\begin{enumerate}

\item[(0)] (Description of $T$) Let $Q \subset \mathbb{P}^4$ be a smooth quadric in $\mathbb{P}^4$.  We take a smooth conic $C$ inside of $Q$.  Let $Y$ be the blow up of $Q$ with center $C$.  Now $T$ is the blow up of $Y$ with center being three exceptional lines (call them $\epsilon_1, \epsilon_2, \epsilon_3$) of the blowing up $Y \rightarrow Q$.  

Equivalently, blow up three points $P_1, P_2, P_3$ (which are the images of $\epsilon_1, \epsilon_2, \epsilon_3$ under the map $Y \rightarrow Q$) on the conic $C$ to obtain $Z \rightarrow Q$.  Then blow up the strict transform $C_Z$ of the conic $C$ to obtain $T \rightarrow Z$.

\vskip.1in

\item[(1)] (Description of the extremal rays on $T$)

\begin{itemize}

\item[$l_1$] : the ruling of the exceptional divisor over $\epsilon_1$ $\cdots$ ($E_1$)

\item[ $l_2$] : the ruling of the exceptional divisor over $\epsilon_2$ $\cdots$ ($E_1$)

\item[$l_3$]  : the ruling of the exceptional divisor over $\epsilon_3$ $\cdots$ ($E_1$)

\item[$l_4$] : the strict transform of the generator of the singular quadric cone on $Q$ with vertex $P_1$ $\cdots$ ($E_1$)

\item[$l_5$] : the strict transform of the generator of the singular quadric cone on $Q$ with vertex $P_2$ $\cdots$ ($E_1$)

\item[$l_6$] : the strict transform of the generator of the singular quadric cone on $Q$ with vertex $P_3$ $\cdots$ ($E_1$)

\item[$l_7$]  : the ruling of the exceptional divisor over the conic $\cdots$ ($E_1$)

\item[$l_8$] : {\bf This is the extremal ray which was missing in \cite{WGBTMM95}, and its description is given as follows.}

The smooth conic $C$ is cut by two hyperplane sections of $\mathbb{P}^4$ out of the smooth quadric $Q \subset \mathbb{P}^4$.  Therefore, we see that the normal bundle of $C$ inside of the quadric is isomorphic to $\mathit{N}_{C/Q} \cong (\mathcal{O}_{\mathbb{P}^4}(1)|_C)^{\oplus 2}$.  Therefore, the exceptional divisor $E_C$ of $Y \rightarrow Q$ is isomorphic to $C \times \mathbb{P}^1 = \mathbb{P}^1 \times \mathbb{P}^1$, i.e., $E_C \cong \mathbb{P}^1 \times \mathbb{P}^1$.  The strict transform, which we denote by $D$, of $E_C$ after blowing up three exceptional lines, remains isomorphic to $E_C$, i.e., $D \cong E_C \cong \mathbb{P}^1 \times \mathbb{P}^1$.  The fiber $\{p\} \times \mathbb{P}^1$ of the first projection generates the extremal ray $l_7$.  The extremal ray $l_8$ is generated by the fiber $\mathbb{P}^1 \times \{q\}$ of the second projection $\cdots$ ($E_1$)

\vskip.1in

In fact, let $\mathrm{cont}_{l_7}:T \rightarrow Z$ be the contraction of the extremal ray $l_7$, where the morphism $T \rightarrow Z$ appears in the second description of the construction of the Fano 3-fold $T$.  Observe that $C_Z$ generates an extremal ray of flopping type on $Z$.  Let $Z \rightarrow W$ be the contraction of this extremal ray of flopping type.  (See below.)  By composing $T \rightarrow Z$ with $Z \rightarrow W$, we obtain the morphism \text{$T \rightarrow W$.}  By construction, we see that the exceptional locus of the morphism  is $D$, $\dim N^1(T/W) = 2$, and that $- K_T$ is relatively ample for the morphism $T \rightarrow W$.  Now $\overline{\mathrm{NE}}(T/W)$ has two extremal rays, one is generated by the fiber $\{p\} \times \mathbb{P}^1$ of the first projection, the other by the fiber $\mathbb{P}^1 \times \{q\}$ of the second projection.  Since $\overline{\mathrm{NE}}(T/W)$ corresponds to the face of $\overline{\mathrm{NE}}(T)$ (cf. Lemma 3-2-4 in \cite{KMM}), the fiber $\mathbb{P}^1 \times \{q\}$ of the second projection gives rise to an extremal ray $l_8$ of $\overline{\mathrm{NE}}(T)$.   (Note that the fiber $f = \mathbb{P}^1 \times \{q\}$ of the second projection is cut out by the strict transform of the hyperplane section of the quadric passing through $C$ and by $D$, which is the strict transform of $E_C$.  From this and via the identification $f \cong \mathbb{P}^1$, it follows immediately that $\mathit{N}_{f/T} = \mathcal{O}_{\mathbb{P}^1} \oplus \mathcal{O}_{\mathbb{P}^1}(-1)$ and that $f \cdot K_T = -1$.)

\vskip.1in

Note: We can see why $C_Z$ generates an extremal ray of flopping type on $Z$ as follows.  Let $H_1, H_2$  be the two hyperplane sections of $\mathbb{P}^4$ which cut the conic $C$ out of the smooth quadric $Q \subset \mathbb{P}^4$, i.e., $H_1 \cap H_2 \cap Q = C$.  Set $G_1 = H_1 \cap Q$ and $G_2 = H_2 \cap Q$.  Then $C_Z$ is cut out by the strict transforms of $G_1$ and $G_2$ under the blow up $Z \rightarrow Q$.  Let $G$ be the general member of the pencil spanned by the strict transforms of $G_1$ and $G_2$.  Then we see that $C_Z \cdot G = \mathrm{deg}\left(\mathcal{O}_{\mathbb{P}^4}(1)|_C\right) - 3 = - 1$ and that $C_Z$ is the only curve which has negative intersection with $G$, as $C_Z$ is the fixed locus of the pencil.  Since $\overline{\mathrm{NE}}(Z)$ is a polyhedral cone each of whose edges is generated by a curve, this implies that $C_Z$ generates one edge.  It is straightforward to see that $C_Z \cdot K_Z = 0$.  Now it follows that the contraction of this extremal ray, negative with respect to $K_Z + \epsilon G$ for some $0 < \epsilon << 1$, is a contraction of flopping type, with respect to $K_Z$, whose exceptional locus only consists of $C_Z$.  (Note that the pair $(Z,\epsilon \Delta)$ for $0 < \epsilon << 1$ has Kawamata log terminal singularities (cf. Remark 4-4-3 in \cite{Matsuki02}).)

\includegraphics[width=5in]{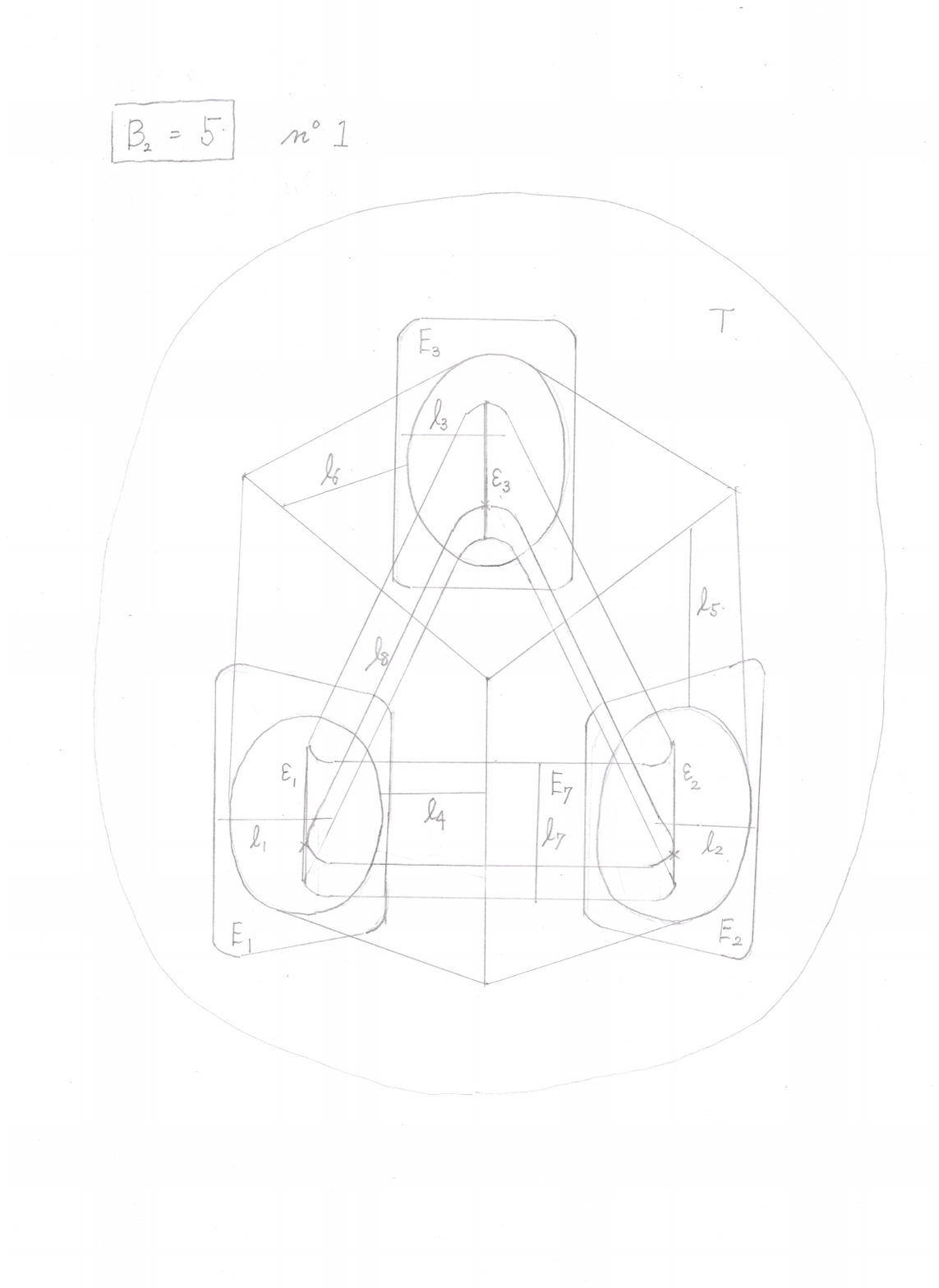}

\text{\bf Oversight: Why the author missed this extremal ray ?}

In \cite{WGBTMM95}, only the extremal rays $l_1, l_2, \cdots, l_7$ were listed.  Let's try to apply the {\bf Criterion for exhaustion} of Claim III-3-2, which should tell if a given set $L$ exhausts all the possible extremal rays.

Let $L$ be the given set of extremal rays.  Then the criterion ($\star$) to check is:

For $l_i \in L$ and the edges $l_{ki} \in M_i$ (Recall that $M_i = \{l_{1i}, l_{2i}, \ldots, l_{ni}\}$ is the set of all the edges of $\overline{\mathrm{NE}}(U_i)$ (i.e., the extremal rays without the negativity condition on the intersection number with $K_{U_i}$) and that $\phi_i: T \rightarrow U_i$ is the contraction of the extremal ray $l_i$.), there exists another extremal ray $l_j \in L$ such that $\phi_i(l_j) = l_{ki}$.

\vskip.1in

{\bf Case: $L = \{l_1, l_2, \ldots, l_7\}$.  The criterion ($\star$) fails.}

\begin{enumerate}

\item[] ($\star$) for $i = 1, 2,3$ checks.  Note that $U_i$ for $i = 1, 2, 3$ is a Fano 3-fold denoted as $B_2 = 4$ with $n^o\ 4$.

\item[] ($\star$) for $i = 4, 5, 6$ fails.  Note that $U_i$ for $i = 4, 5, 6$ is a Fano 3-fold denoted as $B_2 = 4$ with $n^o\ 12$.

The description of the Fano 3-fold denoted as $B_2 = 4$ with $n^o\ 12$ is given as follows: First blow up $\mathbb{P}^3$ with center being a line to obtain $Y \rightarrow \mathbb{P}^3$.  Then blow up $Y$ with center being two exceptional lines $L_1, L_2$ (of the first blow up) to obtain $U_i \rightarrow Y$.  The extremal rays on $U_i$ are:

\begin{itemize}

\item $l_{1i}$: the ruling of the exceptional divisor over $L_1$

\item $l_{2i}$: the ruling of the exceptional divisor over $L_2$

\item $l_{3i}$: the strict transform of the ruling of the exceptional divisor $E_3$ over the line (which is the center of the first blow up) on $\mathbb{P}^3$

\item $l_{4i}$: the the other ruling (than $l_3$) of $E_3 \cong \mathbb{P}^1 \times \mathbb{P}^1$

\end{itemize}

Now for $k = 4$, the criterion ($\star$) fails.  That is to say, there is no $l_j\ (j = 1, 2, \ldots, 7)$ such that $\phi_i(l_j) = l_{4i}$.  In fact, we have $\phi_i(l_8) = l_{4i}$.

\item[] ($\star$) for $i = 7$ fails.  As in the description of the extremal ray $l_8$, $U_i = Z$ has an extremal ray of flopping type.  For $l_{ki}$ being this extremal ray of flopping type, the criterion ($\star$) fails.  That is to say, there is no $l_j\ (j = 1, 2, \ldots, 7)$ such that $\phi_i(l_j) = l_{ki}$.  In fact, we have $\phi_i(l_8) = l_{ki}$.  

\end{enumerate}

{\bf Case: $L = \{l_1, l_2, \ldots, l_7, l_8\}$. The criterion ($\star$) checks.}

\begin{enumerate}

\item[] ($\star$) for $i = 1, 2,3$ checks.  In fact, ($\star$) for $i = 1, 2,3$ checks even with $l_8$ missing from $L$ as seen above.  Therefore, it of course checks when $l_8$ is present in $L$.

\item[] ($\star$) for $i = 4, 5, 6$ checks.  Since we have $\phi_i(l_8) = l_{4i}$, the criterion ($\star$) checks when $l_8$ is present in $L$.

\item[] ($\star$) for $i = 7, 8$ checks.  Here in order to check the criterion ($\star$) for $i = 7, 8$, we face the task of determining the edges (i.e., the extremal rays without the negativity condition on the intersection number with $K_{U_i}$ of $\overline{\mathrm{NE}}(U_i)$.  Since $U_i\ (i = 7,8)$ is NOT a Fano 3-fold (even though $U_i$ is still smooth), one cannot simply use the table of all the possible extremal rays for Fano 3-folds with $B_2 = 4$.  However, using the information of the extemal rays for $U_i\  (i = 1, 2, \ldots, 6)$, which is a Fano 3-fold, one can still carry out this task.  Then check the criterion ($\star$) for the edges of the cone $\overline{\mathrm{NE}}(U_i)$.  Details are left to the reader as an exercise.

\end{enumerate}

\vskip.1in

\end{itemize}

\newpage

\item[(2)] Table for the intersection pairing

\begin{itemize}
\item[] $E_1$: the exceptional divisor for the blow up of $\epsilon_1$
\item[] $E_2$: the exceptional divisor for the blow up of $\epsilon_2$
\item[] $E_3$: the exceptional divisor for the blow up of $\epsilon_3$
\item[] $E_7$: the strict transform of the exceptional divisor for the blow up $Y \rightarrow Q$ of the conic $C$
\item[] $\pi^*\mathcal{O}_Q(1) = \pi^*\{\mathcal{O}_{\mathbb{P}^4}(1)|_Q\}$ where $\pi: T \rightarrow Q$
\end{itemize}

\begin{center}
\begin{tabular}{|c|c|c|c|c|c|c|}
\hline
& $E_1$ & $E_2$ & $E_3$ & $E_7$ & $\pi^*\mathcal{O}_Q(1)$ & $- K_T$ \\
\hline
$l_1$ & $- 1$ & 0 & 0 & 1 & 0 & 1 \\
\hline
$l_2$ & 0 & $- 1$ & 0 & 1 & 0 & 1 \\
\hline
$l_3$ & 0 & 0 & $- 1$ & 1 & 0 & 1 \\
\hline
$l_4$ & 1 & 0 & 0 & 0 & 1 & 1 \\
\hline
$l_5$ & 0 & 1 & 0 & 0 & 1 & 1 \\
\hline
$l_6$ & 0 & 0 & 1 & 0 & 1 & 1 \\
\hline
$l_7$ & 0 & 0 & 0 & $- 1$ & 0 & 1 \\
\hline
$l_8$ & 1 & 1 & 1 & $- 1$ & 2 & 1 \\
\hline
$l_{11}$ & 1 & 0 & 0 & $- 1$ & 0 & $- 1$ \\
\hline
$l_{21}$ & 0 & $- 1$ & 0 & 1 & 0 & 1 \\
\hline
$l_{31}$ & 0 & 0 & $- 1$ & 1 & 0 & 1 \\
\hline
$l_{41}$ & 0 & 0 & 0 & 1 & 1 & 2 \\
\hline
$l_{51}$ & 0 & 1 & 0 & 0 & 1 & 1  \\
\hline
$l_{61}$ & 0 & 0 & 1 & 0 & 1 & 1 \\
\hline
$l_{71}$ & 0 & 0 & 0 & $- 1$ & 0 & 1 \\
\hline
$l_{81}$ & 0 & 1 & 1 & 0 & 2 & 2 \\
\hline
$l_{12}$ & $- 1$  & 0  & 0 & 1 & 0 & 1  \\
\hline
$l_{22}$ & 0 & 1 & 0 & $- 1$ & 0 & $- 1$ \\
\hline
$l_{32}$ & 0  & 0 & $- 1$ & 1 & 0 & 1 \\
\hline
$l_{42}$ & 1 & 0 & 0  & 0 & 1 & 1 \\
\hline
$l_{52}$ & 0 & 0 & 0 & 1 & 1 & 2 \\
\hline
$l_{62}$ & 0 & 0 & 1 & 0 & 1 & 1 \\
\hline
$l_{72}$ & 0 & 0 & 0 & $- 1$  & 0 & 1 \\
\hline
$l_{82}$ & 1  & 0 & 1 & 0 & 2 & 2  \\
\hline
$l_{13}$ & $- 1$ & 0 & 0 & 1 & 0 & 1  \\
\hline
$l_{23}$ & 0 & $- 1$ & 0 & 1 & 0 & 1 \\
\hline
$l_{33}$ & 0 & 0 & 1 & $- 1$  & 0 & $- 1$  \\
\hline
$l_{43}$ & 1 & 0 & 0 & 0 & 1 & 1 \\
\hline
$l_{53}$ & 0 & 1 & 0 & 0 & 1 & 1   \\
\hline
$l_{63}$ & 0 & 0  & 0 & 1 & 1 & 2 \\
\hline
$l_{73}$ & 0 & 0 & 0  & $- 1$ & 0 & 1 \\
\hline
$l_{83}$ & 1  & 1 & 0 & 0 & 2 & 2  \\
\hline
\end{tabular}
\end{center}

\begin{center}
\begin{tabular}{|c|c|c|c|c|c|c|}
\hline
& $E_1$ & $E_2$ & $E_3$ & $E_7$ & $\pi^*\mathcal{O}_Q(1)$ & $- K_T$ \\
\hline
$l_{14}$ & 1 & 0 & 0 & 1 & 2 & 3 \\
\hline
$l_{24}$ & 0 & $- 1$ & 0 & 1 & 0 & 1 \\
\hline
$l_{34}$ & 0 & 0 & $- 1$ & 1 & 0 & 1  \\
\hline
$l_{44}$ & $- 1$ & 0 & 0 & 0 & $- 1$ & $- 1$ \\
\hline
$l_{54}$ & 1 & 1 & 0 & 0 & 2 & 2   \\
\hline
$l_{64}$ & 1 & 0  & 1 & 0 & 2 & 2 \\
\hline
$l_{74}$ & 0 & 0 & 0 & $- 1$ & 0 & 1 \\
\hline
$l_{84}$ & 1 & 1 & 1 & $- 1$ & 2 & 1 \\

\hline
$l_{15}$ & $- 1$ & 0 & 0 & 1 & 0 & 1 \\
\hline
$l_{25}$ & 0 & 1 & 0 & 1 & 2 & 3 \\
\hline
$l_{35}$ & 0  & 0 & $- 1$ & 1 & 0 & 1  \\
\hline
$l_{45}$ &1 & 1 & 0 & 0 & 2 & 2 \\
\hline
$l_{55}$ & 0 & $- 1$ & 0 & 0 & $- 1$ & $- 1$   \\
\hline
$l_{65}$ & 0 & 1  &  0 & 0 & 2 & 2 \\
\hline
$l_{75}$ & 0 & 0  & 0 & $- 1$ & 0 & 1 \\
\hline
$l_{85}$ & 1 & 1  & 1 & $- 1$ & 2 & 1  \\

\hline
$l_{16}$ & $- 1$ & 0 & 0 & 1 & 0 & 1 \\
\hline
$l_{26}$ & 0 & $- 1$ & 0 & 1  & 0 & 1 \\
\hline
$l_{36}$ & 0  & 0 & 0 & 1  &  2 & 3  \\
\hline
$l_{46}$ & 1 & 0 & 1 & 0 & 2 & 2 \\
\hline
$l_{56}$ & 0 & 1 & 1 & 0 & 2 & 2   \\
\hline
$l_{66}$ & 0 & 0 & $- 1$  & 0 & $- 1$ & $- 1$  \\
\hline
$l_{76}$ & 0 & 0  & 0 & $- 1$  & 0 & 1 \\
\hline
$l_{86}$ & 1 & 1  & 1 & $- 1$ & 2 & 1  \\

\hline
$l_{17}$ & $- 1$ & 0 & 0 & 0 & 0 & 2 \\
\hline
$l_{27}$ & 0  & $- 1$  & 0 & 0  & 0  & 2 \\
\hline
$l_{37}$ & 0  & 0 & $- 1$  & 0 & 0  & 2  \\
\hline
$l_{47}$ & 1 & 0 & 0 & 0 & 1 & 1 \\
\hline
$l_{57}$ & 0 & 1 & 0 & 0 & 1 & 1   \\
\hline
$l_{67}$ & 0 & 0  & 1 & 0 & 1 & 1 \\
\hline
$l_{77}$ & 0 & 0 & 0 & 1 & 0 & $- 1$ \\
\hline
$l_{87}$ & 1 & 1  & 1 & 0 & 2 & 0  \\
\hline
$l_{18}$ & 0 & 1 & 1 & 0 & 2 & 2 \\
\hline
$l_{28}$ & 1 & 0 & 1 & 0  & 2 & 2  \\
\hline
$l_{38}$ & 1  & 1 & 0 & 0 & 2  &  2 \\
\hline
$l_{48}$ & 1 & 0 & 0 & 0 & 1  & 1 \\
\hline
$l_{58}$ & 0 & 1 & 0 & 0 & 1 & 1   \\
\hline
$l_{68}$ & 0 & 1  & 0 & 0 & 1 & 1  \\
\hline
$l_{78}$ & $-1$ & $-1$  & $-1$ & 0  & $- 2$ & 0 \\
\hline
$l_{88}$ & $-1$ & $-1$  & $- 1$ & 1 & $-2$ & $- 1$  \\

\hline
\end{tabular}
\end{center}

\newpage

Note: There are mistakes in the table for the intersection pairing in \cite{WGBTMM95} for $B_2 = 5$ with $n^o\ 1$:

\begin{center}
\begin{tabular}{|c|c|c|c|c|c|c|}
\hline
& $E_1$ & $E_2$ & $E_3$ & $E_7$ & $\pi^*\mathcal{O}_Q(1)$ & $- K_T$ \\
\hline
$l_{25}$ & 0 & 1 & 0 & 1 & 2 & $- 3$ \\
\hline
\end{tabular}
\end{center}

The first part of the mistakes, writing $- 3$ instead of $3$, is a simple typographical error.

\begin{center}
\begin{tabular}{|c|c|c|c|c|c|c|}
\hline
& $E_1$ & $E_2$ & $E_3$ & $E_7$ & $\pi^*\mathcal{O}_Q(1)$ & $- K_T$ \\
\hline
$l_{17}$ & $- 1$ & 0 & 0 & 1 & 0 & 1 \\
\hline
$l_{27}$ & 0  & $- 1$  & 0 & 1  & 0  & 1 \\
\hline
$l_{37}$ & 0  & 0 & $- 1$  & 3 & 0  & $- 1$  \\
\hline
$l_{47}$ & 1 & 0 & 0 & 1 & 1 & 0 \\
\hline
$l_{57}$ & 0 & 1 & 0 & 1 & 1 & 0   \\
\hline
$l_{67}$ & 0 & 0  & 1 & $- 1$ & 1 & 2 \\
\hline
\end{tabular}
\end{center}

The second part of the mistakes, wrong numbers in the column indicating the intersection pairing with $E_7$ (and accordingly wrong numbers in the column indicating the intersection pairing with $- K_T$, is apparently caused by somehow not realizing that the strict transform, after the flop of the 4-fold, of the divisor corresponding to $E_7$ is away from the strict transform of the zero section and hence that those numbers should all have been $0$.

\vskip.1in

These mistakes are corrected in the new table.

\vskip.1in

\item[(3)] When we contract $l_1$ (resp. $l_2, l_3, l_4, l_5, l_6$, we get $n^o\ 4$ (resp. $4, 4, 12, 12, 12$) with $B_2 = 4$, while the contraction of $l_7$ or $l_8$ leads to a 3-fold which is not a Fano 3-fold.  The analysis of the KKMR decomposition is still reduced to the previous ones with lower Picard numbers.  (See the comments in (1) above.)  Thus we omit the details here.

\item[(4)] ($E_1$), ($E_2$), ($E_5$) and their inverses, and ($F$) as described in \S 5 as ``the flop of type ($F$) in the item (4)''.  For example, we have two minimal models, obtained by flopping the extremal rays $l_7$ and $l_8$.  The Weyl chambers associated with these two minimal models share the codimension-one face, and connected by a flop of type ($F$) as described in \S 5 as ``the flop of type ($F$) in the item (4)''.  The label ($F$) was missing in \cite{WGBTMM95}, as the extemal ray $l_8$ was also missing.  Moreover, we also have flops of Types (Others) associated with $\boxed{B_2 = 4}$ $n^o\ 12$.

\vskip.1in

\item[(5)] $\mathrm{WG}_T = A_2$

\end{enumerate}

\newpage

\section{Analysis of the Fano 3-fold denoted as $n^o\ 13$ with $B_2 = 4$ in Table 4}

In this section, we present an analysis of the Fano 3-fold $T$ denoted as $n^o\ 13$  with $B_2 = 4$ in \cite{MMManu03}, which was missing from the original classification table in \cite{MMManu81}.

\vskip.1in

$\boxed{B_2 = 4}$ : $n^o\ 13$

\vskip.1in

\begin{enumerate}

\item[(0)] (Description of $T$) The Fano 3-fold $T$ is obtained as the blow up $\pi: T \rightarrow \mathbb{P}_1^1 \times \mathbb{P}_2^1 \times \mathbb{P}_3^1$ with center a curve $C$ of tridegree $(1, 1, 3)$.

\item[(1)] (Description of the extremal rays on $T$)

\begin{itemize}

\item[ $l_1$] : the strict transform of the fiber $p_{23}^{-1}(p_{23}(P))$ with $P \in C$ where $p_{23}: \mathbb{P}_1^1 \times \mathbb{P}_2^1 \times \mathbb{P}_3^1 \rightarrow \mathbb{P}_2^1 \times {P}_3^1$ $\cdots$ ($E_1$)

\item[$l_2$] : the strict transform of the fiber $p_{13}^{-1}(p_{13}(P))$ with $P \in C$ where $p_{13}: \mathbb{P}_1^1 \times \mathbb{P}_1^1 \times \mathbb{P}_3^1 \rightarrow \mathbb{P}_1^1 \times {P}_3^1$ $\cdots$ ($E_1$)

\item[$l_3$] : the strict transform of the fiber $p_{12}^{-1}(p_{12}(P))$ with $P \in C$ where $p_{12}: \mathbb{P}_1^1 \times \mathbb{P}_1^1 \times \mathbb{P}_2^1 \rightarrow \mathbb{P}_1^1 \times {P}_2^1$ $\cdots$ ($E_1$)

\item[$l_4$] : the ruling of the exceptional divisor for $\pi$ over $C$ $\cdots$ ($E_1$)

\item[$l_5$] : Observe that $\Delta \times \mathbb{P}_3^1 \cong \mathbb{P}^1 \times \mathbb{P}^1$, since the diagonal $\Delta \subset \mathbb{P}_1^1 \times  \mathbb{P}_2^1$ is isomorphic to $\mathbb{P}^1$.  The strict transform of the fiber of the first projection $\{p\} \times \mathbb{P}^1\ (p \in \Delta)$ is nothing but $l_3$.  The strict transform of the fiber of the second projection $\mathbb{P}^1 \times \{q\}\ (q \in {P}_3^1)$ is $l_5$.

In fact, let $\mathrm{cont}_{l_3}: T \rightarrow Z$ be the contraction of the extemal ray $l_3$.  Observe that the image of $l_5$, which we denote by $C_Z$, generates an extremal ray of flopping type on $Z$.  Let $Z \rightarrow W$ be the contraction of this extremal ray of flopping type.  By composing $T \rightarrow Z$ with $Z \rightarrow W$, we obtain the morphism $T \rightarrow W$.  By construction, we see that the exceptional locus of the morphism  is $D$ (the strict transform of $\Delta \times {P}_3^1$, which is still isomorphic to $\mathbb{P}^1 \times \mathbb{P}^1$), $\dim N^1(T/W) = 2$, and that $- K_T$ is relatively ample for the morphism $T \rightarrow W$.  Now $\overline{\mathrm{NE}}(T/W)$ has two extremal rays, one is generated by the fiber $\{p\} \times \mathbb{P}^1$ of the first projection, the other by the fiber $\mathbb{P}^1 \times \{q\}$ of the second projection.  Since $\overline{\mathrm{NE}}(T/W)$ corresponds to the face of $\overline{\mathrm{NE}}(T)$ (cf. Lemma 3-2-4 in \cite{KMM}), the fiber $\mathbb{P}^1 \times \{q\}$ of the second projection gives rise to an extremal ray $l_5$ of $\overline{\mathrm{NE}}(T)$.   (Note that the fiber $f = \mathbb{P}^1 \times \{q\}$ of the second projection is cut out by the strict transform of $\pi_3^*(q)$ and by $D$.  From this and via the identification $f \cong \mathbb{P}^1$, it follows immediately that $\mathit{N}_{f/T} = \mathcal{O}_{\mathbb{P}^1} \oplus \mathcal{O}_{\mathbb{P}^1}(-1)$ and that $f \cdot K_T = -1$.)

\vskip.1in

Note: We can see why $C_Z$ generates an extremal ray of flopping type on $Z$ as follows.  Set $G_1 = \pi_3^*(p), G_2 = \pi_3^*(q)$ for two distinct points $p, q \in \mathbb{P}_3^1$.  Then $C_Z$ is cut out by the images of $G_1$ and $G_2$ under the morphism $T \rightarrow Z$.  Let $G$ be the general member of the pencil spanned by the images of $G_1$ and $G_2$.  Then we see that $C_Z \cdot G = l_5 \cdot \mathrm{cont}_{l_3}^*(G) = l_5 \cdot (G_1 + D) = l_5 \cdot D = \Delta \times \{q\} \cdot \{p_{12}^*(\Delta) - E_4\} = 2 - 3 = -1$ and that $C_Z$ is the only curve which has negative intersection with $G$, as $C_Z$ is the fixed locus of the pencil.  Since $\overline{\mathrm{NE}}(Z)$ is a polyhedral cone each of whose edges is generated by a curve, this implies that $C_Z$ generates one edge.  It is straightforward to see that $C_Z \cdot K_Z = 0$.  Now it follows that the contraction of this extremal ray, negative with respect to $K_Z + \epsilon G$ for some $0 < \epsilon << 1$, is a contraction of flopping type, with respect to $K_Z$, whose exceptional locus only consists of $C_Z$.  (Note that the pair $(Z,\epsilon \Delta)$ for $0 < \epsilon << 1$ has Kawamata log terminal singularities (cf. Remark 4-4-3 in \cite{Matsuki02}).)

\end{itemize}

\includegraphics[width=5in]{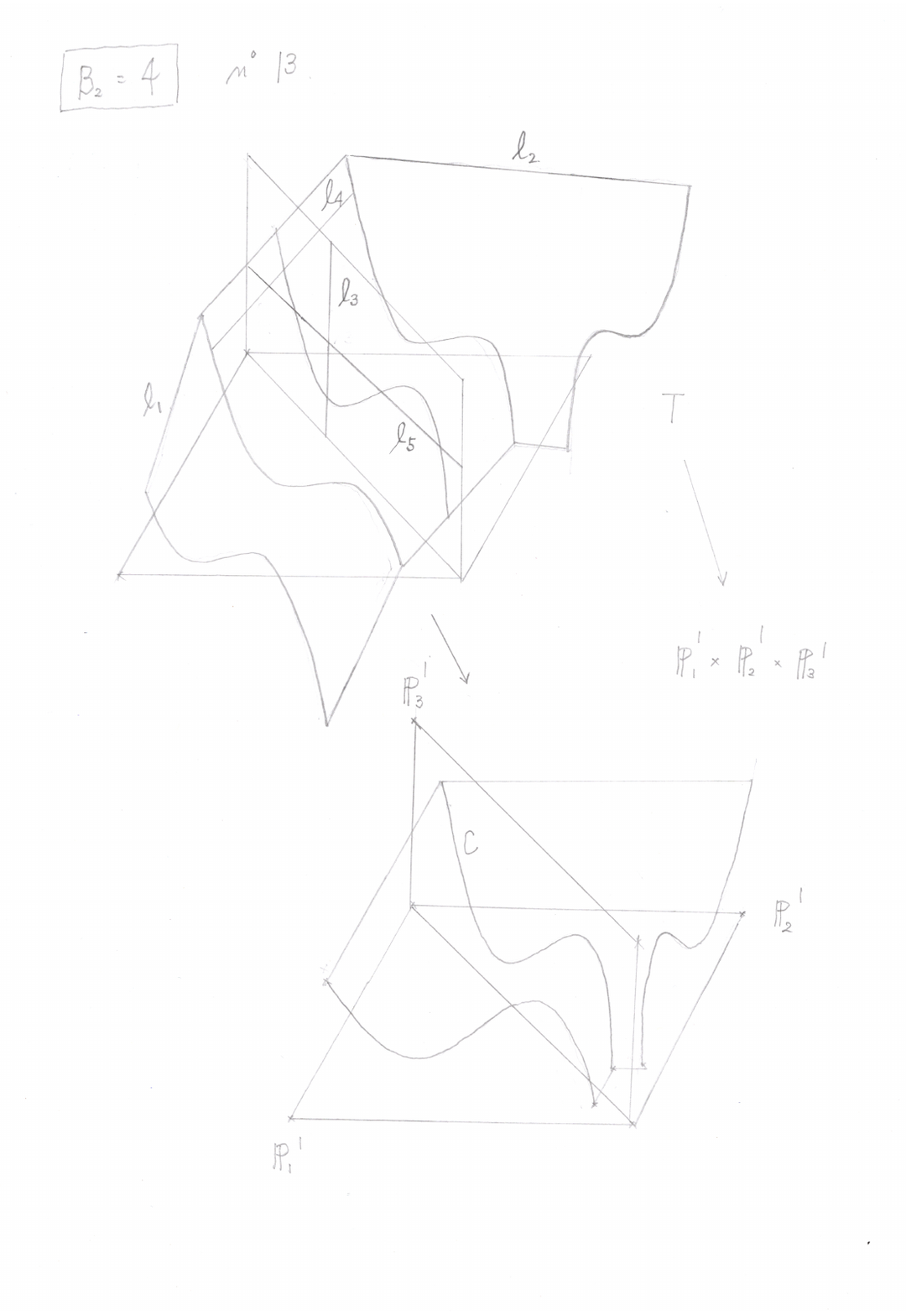}

\item[(2)] Table for the intersection pairing

\vskip.1in

\begin{center}
\begin{tabular}{|c|c|c|c|c|c|}
\hline
& $\pi_1^*\mathcal{O}_{\mathbb{P}^1}(1)$ & $\pi_2^*\mathcal{O}_{\mathbb{P}^1}(1)$ & $\pi_3^*\mathcal{O}_{\mathbb{P}^1}(1)$ & $E_4$ & $- K_T$ \\
\hline
$l_1$ & 1 & 0 & 1 & 1 & 1 \\
\hline
$l_2$ & 0 & 1 & 0 & 1 & 1 \\
\hline
$l_3$ & 0 & 0 & 1 & 1 & 1 \\
\hline
$l_4$ & 0 & 0 & 0 & $- 1$ & 1 \\
\hline
$l_5$ & 1 & 1 & 0 & 3 & 1 \\
\hline
$l_{11}$ & $- 1$ & 0 & 0 & $- 1$ & $- 1$ \\
\hline
$l_{21}$ & 0 & 1 & 0 & 1  & 1  \\
\hline
$l_{31}$ & 0 & 0 & 1 & 1 & 1  \\
\hline
$l_{41}$ & 1 & 0 & 0 & 0 & 2 \\
\hline
$l_{51}$ & 1 & 1 & 0 & 3 & 1 \\
\hline
$l_{12}$ & 1 & 0  & 0 & 1 & 1 \\
\hline
$l_{22}$ & 0 & $- 1$ & 0  & $- 1$ & $- 1$ \\
\hline
$l_{32}$ & 0 & 0 & 1 & 1 & 1  \\
\hline
$l_{42}$ & 0 & 1 & 0 & 0 & 2 \\
\hline
$l_{52}$ & 1 & 1 & 0 & 3 & 1 \\
\hline
$l_{13}$ & 1 & 0 & 0 & 1 & 1  \\
\hline
$l_{23}$ & 0 & 1 & 0 & 1  & 1  \\
\hline
$l_{33}$ & 0 & 0 & $- 1$  & $- 1$  & $- 1$  \\
\hline
$l_{43}$ & 0 & 0 & 1 & 0 & 2 \\
\hline
$l_{53}$ & 1 & 1 & $- 1$ & 2 & 0 \\
\hline
$l_{14}$ & 1  & 0 & 0 & 0 & 2 \\
\hline
$l_{24}$ & 0 & 1 & 0 & 0  & 2  \\
\hline
$l_{34}$ & 0 & 0 & 1 & 0 & 2 \\
\hline
$l_{44}$ & 0 & 0 & 0 & 1 &  $- 1$ \\
\hline
$l_{54}$ & 1 & 1 & 0 & 0 & 4 \\
\hline
$l_{15}$ & 1  & 0 & 0 & 1 & 1 \\
\hline
$l_{25}$ & 0 & 1 & 0 & 1  & 1  \\
\hline
$l_{35}$ & $- 1$ & $- 1$ & 1 & 2 & 0 \\
\hline
$l_{45}$ & 1 & 1 & 0 & 2 & 2  \\
\hline
$l_{55}$ & $- 1$  & $- 1$ & 0 & $- 3$ & $- 1$ \\
\hline
\end{tabular}
\end{center}

\newpage

\item[(3)] We only describe the nef cone $\overline{\mathrm{Amp}}(T)$.

\includegraphics[width=5in]{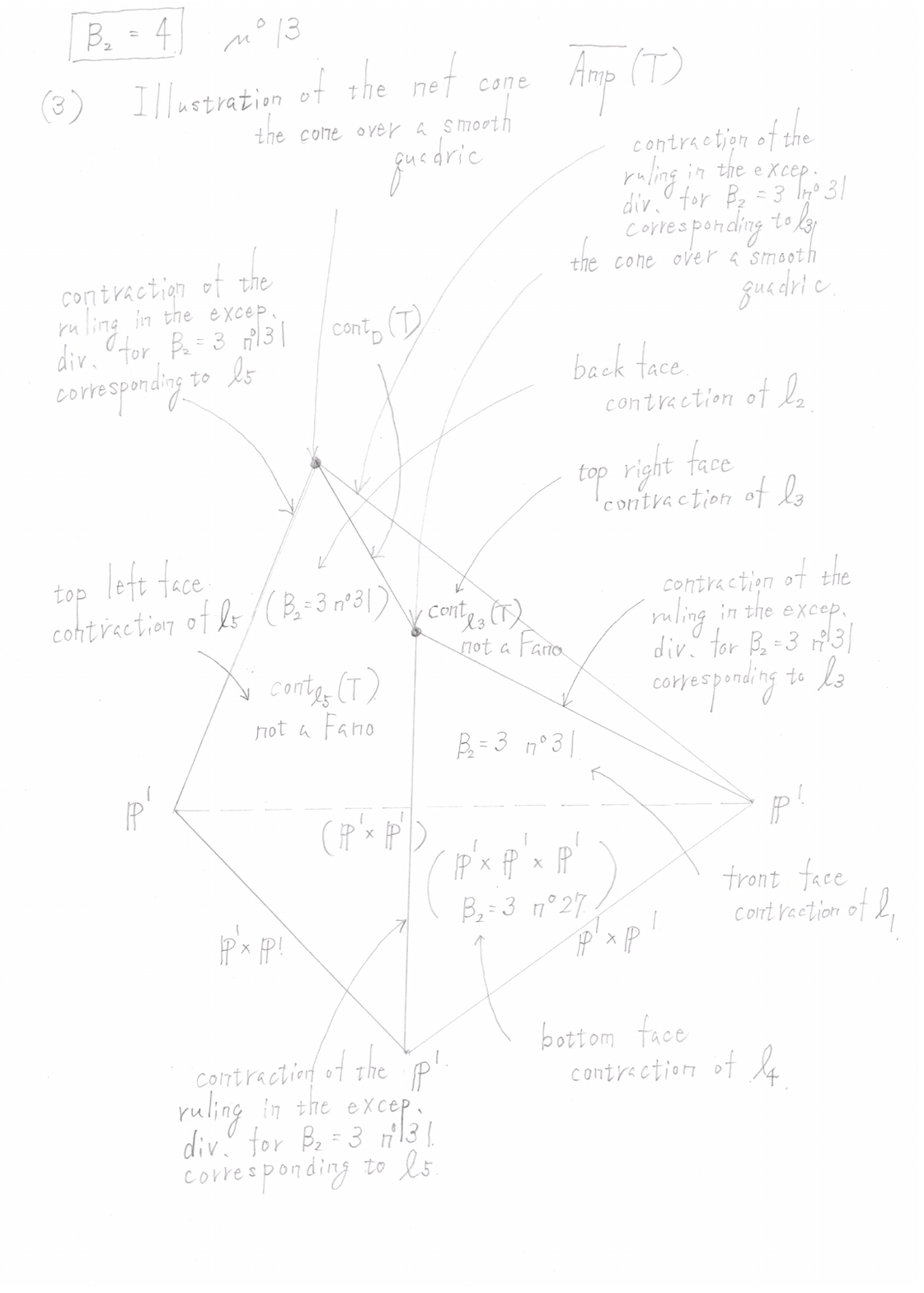}

\item[(4)] ($E_1$) and its inverse, and ($F$) (See \S 5 for the explanation of the type ($F$) appearing in the item (4)).

\item[(5)] $\mathrm{WG}_T = A_1$

\end{enumerate}

\vskip.1in

The description of the Fano 3-fold denoted as $n^o\ 13$  with $B_2 = 4$ in \cite{MMManu03} is quite similar to the one of the Fano 3-fold denoted as $n^o\ 3$  with $B_2 = 4$ in \cite{MMManu81}, which we present here again for the comparison and for the reference of the reader, with a slightly different description from the one given in \cite{WGBTMM95}.  The only difference is that, starting from $\mathbb{P}_1^1 \times \mathbb{P}_2^1 \times \mathbb{P}_3^1$, the former blows up a curve of tridegree $(1,1,3)$, while the latter blows up a curve of tridegree $(1,1,2)$.  However, it is this small difference that brings quite a change in the analysis: the former has 5 extremal rays, while the latter has only 4.  Accordingly, the shape of the nef cone of the former, which is a pyramid, is remarkably different from the one for the latter, which is a tetrahedron.

\vskip.1in

$\boxed{B_2 = 4}$ : $n^o\ 3$

\vskip.1in

\begin{enumerate}

\item[(0)] (Description of $T$) The Fano 3-fold $T$ is obtained as the blow up $\pi: T \rightarrow \mathbb{P}_1^1 \times \mathbb{P}_2^1 \times \mathbb{P}_3^1$ with center a curve $C$ of tridegree $(1, 1, 2)$.

\item[(1)] (Description of the extremal rays on $T$)

\begin{itemize}

\item[ $l_1$] : the strict transform of the fiber $p_{23}^{-1}(p_{23}(P))$ with $P \in C$ where $p_{23}: \mathbb{P}_1^1 \times \mathbb{P}_2^1 \times \mathbb{P}_3^1 \rightarrow \mathbb{P}_2^1 \times {P}_3^1$ $\cdots$ ($E_1$)

\item[$l_2$] : the strict transform of the fiber $p_{13}^{-1}(p_{13}(P))$ with $P \in C$ where $p_{13}: \mathbb{P}_1^1 \times \mathbb{P}_1^1 \times \mathbb{P}_3^1 \rightarrow \mathbb{P}_1^1 \times {P}_3^1$ $\cdots$ ($E_1$)

\item[$l_3$] : the strict transform of the fiber $p_{12}^{-1}(p_{12}(P))$ with $P \in C$ where $p_{12}: \mathbb{P}_1^1 \times \mathbb{P}_1^1 \times \mathbb{P}_3^1 \rightarrow \mathbb{P}_1^1 \times {P}_2^1$ $\cdots$ ($E_1$)

\item[$l_4$] : the ruling of the exceptional divisor for $\pi$ over $C$ $\cdots$ ($E_1$)

\end{itemize}

\includegraphics[width=3.5in]{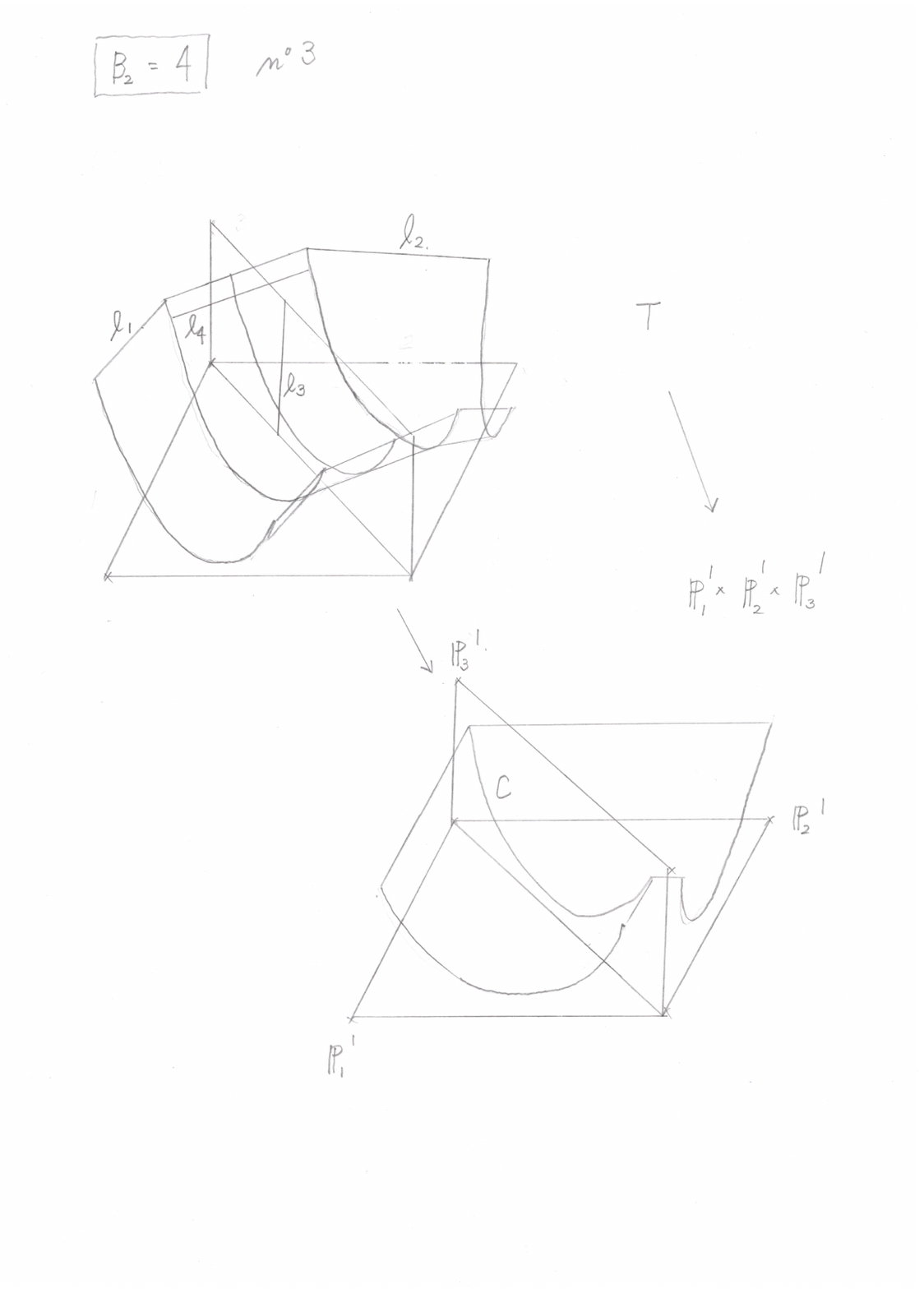}

\item[(2)]

\begin{center}
\begin{tabular}{|c|c|c|c|c|c|}
\hline
& $\pi_1^*\mathcal{O}_{\mathbb{P}^1}(1)$ & $\pi_2^*\mathcal{O}_{\mathbb{P}^1}(1)$ & $\pi_3^*\mathcal{O}_{\mathbb{P}^1}(1)$ & $E_4$ & $- K_T$ \\
\hline
$l_1$ & 1 & 0 & 0 & 1 & 1 \\
\hline
$l_2$ & 0 & 1 & 0 & 1 & 1 \\
\hline
$l_3$ & 0 & 0 & 1 & 1 & 1 \\
\hline
$l_4$ & 0 & 0 & 0 & $- 1$ & 1 \\
\hline
$l_{11}$ & $- 1$ & 0 & 0 & $- 1$ & $- 1$ \\
\hline
$l_{21}$ & 0 & 1 & 0 &  1 & 1  \\
\hline
$l_{31}$ & 0 & 0 & 1 & 1 & 1  \\
\hline
$l_{41}$ & 1 & 0 & 0 & 0 & 2 \\
\hline
$l_{12}$ & 1 & 0 & 0 & 1 & 1 \\
\hline
$l_{22}$ & 0 & $- 1$ & 0 & $- 1$ & $- 1$ \\
\hline
$l_{32}$ & 0 & 0 & 1 & 1 & 1  \\
\hline
$l_{42}$ & 0 & 1 & 0 & 0 & 2 \\
\hline
$l_{13}$ & 1 & 0 & 0 & 1 & 1  \\
\hline
$l_{23}$ & 0 & 1 & 0 &  1 & 1  \\
\hline
$l_{33}$ & 0 & 0 & $- 1$ & $- 1$ & $- 1$ \\
\hline
$l_{43}$ & 0 & 0 & 1& 0 & $2$ \\
\hline
$l_{14}$ & 1 & 0 & 0 & 0 & 2 \\
\hline
$l_{24}$ & 0 & 1 & 0 &  0 & 2  \\
\hline
$l_{34}$ & 0 & 0 & 1 & 0 & 2 \\
\hline
$l_{44}$ & 0 & 0 & 0 & 1 & $- 1$ \\
\hline
\end{tabular}
\end{center}

\vskip.1in

Remark: In the table of the intersection parings in \cite{WGBTMM95} for $B_2 = 4$ with $n^o\ 3$, there are some obvious mistakes in the rows for $l_4, l_{21}, l_{12}, l_{43}$ as below.  Here in the above table we correct those mistakes.

\vskip.11in

Mistakes in the intersection pairing table in \cite{WGBTMM95} for $B_2 = 4$ with $n^o\ 3$:

\begin{center}
\begin{tabular}{|c|c|c|c|c|c|}
\hline
& $\pi_1^*\mathcal{O}_{\mathbb{P}^1}(1)$ & $\pi_2^*\mathcal{O}_{\mathbb{P}^1}(1)$ & $\pi_3^*\mathcal{O}_{\mathbb{P}^1}(1)$ & $E_4$ & $- K_T$ \\
\hline
$l_4$ & 0 & 0 & 0 & $- 1$ & $- 1$ \\
\hline
$l_{21}$ & 1 & 1 & 0 & 2 & 2 \\
\hline
$l_{12}$ & 1 & 1 & 0 & 2 & 2 \\
\hline
$l_{43}$ & 0 & 0 & $- 1$ & 0 & $- 2$ \\
\hline
\end{tabular}
\end{center}

\newpage

\item[(3)] We only describe the nef cone $\overline{\mathrm{Amp}}(T)$.

\includegraphics[width=5in]{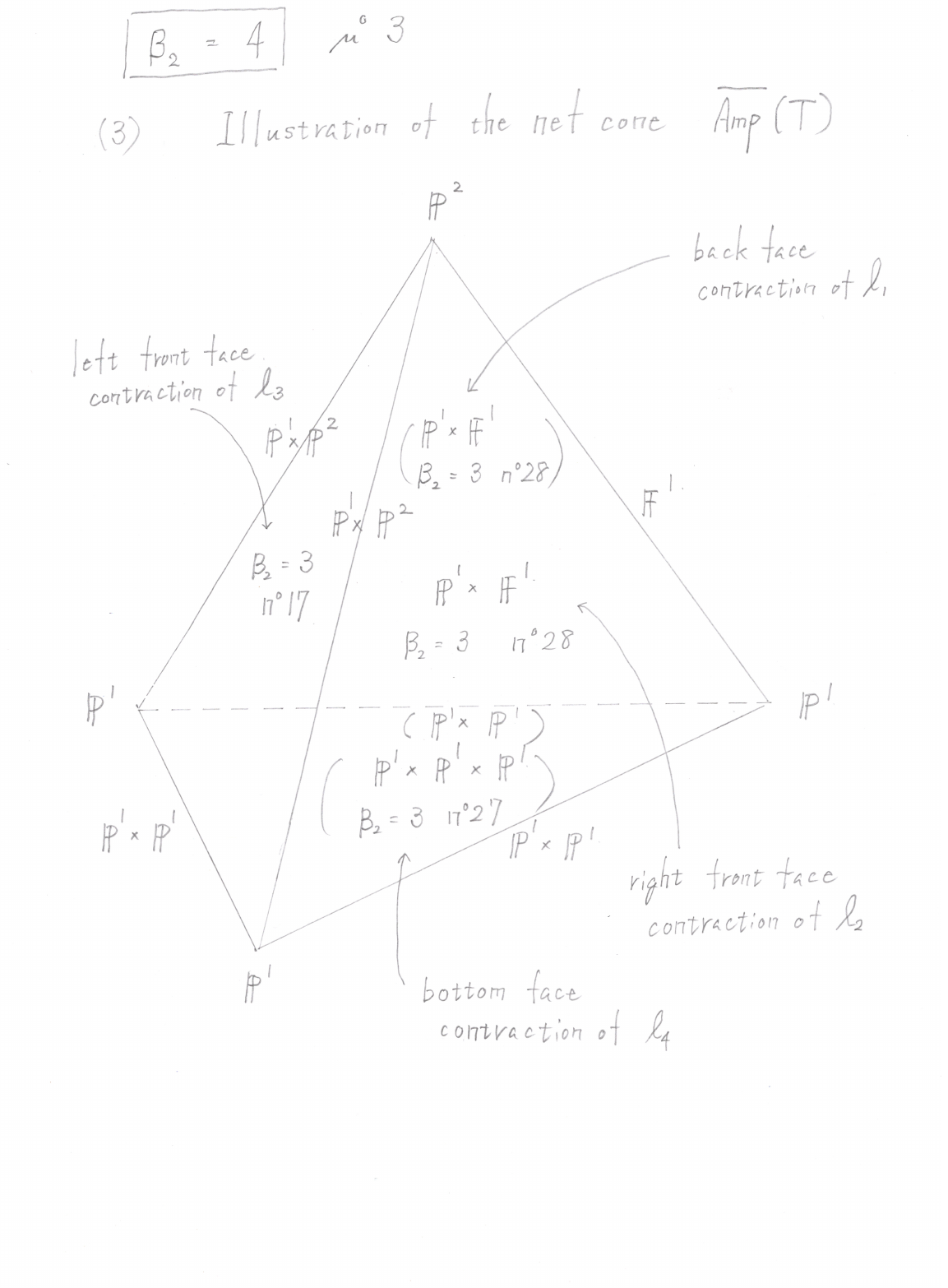}

\item[(4)] ($E_1$) and its inverse

\item[(5)] $\mathrm{WG}_T = A_1$

\vskip.1in

Note: The reason why the Weyl Group for $n^o\ 3$ with $B_2 = 4$, reflecting the symmetry of the KKMR decomposition, is $A_1$ and NOT $A_2$ does NOT come from the fact that the tridegree of the curve is $(1,1,2)$ having a heavier weight on the 3rd factor.  The true reason is that the contraction morphism $\mathbb{P}^1 \times \mathbb{F}_1 \rightarrow \mathbb{F}_1$, corresponding to the 1-dimensional edge labeled $\mathbb{F}_1$ coming from the top vertex, is not birational and hence it does NOT give rise to the flopping going from the chamber associated to $\mathbb{P}^1 \times \mathbb{F}_1$ to another chamber, while in contrast the contraction morphism $n^o\  17 \rightarrow \mathbb{P}^1 \times \mathbb{P}_2$  is birational and hence it does give rise to the flopping going from the chamber associated with $n^o\ 17$ to another chamber associated with $\mathbb{P}^1 \times \mathbb{P}_2$.  This breaks the symmetry among the three extremal rays $l_1,l_2,l_3$, only preserving the symmetry between $l_1$ and $l_2$.

\end{enumerate}

\section{Clarification of the statements about the types of our 4-fold flops}

In the list for the analysis of the Fano 3-folds in \cite{WGBTMM95}, we claim to discuss

\begin{enumerate}

\item the list of all the extremal rays,

\item the intersection pairings of the extremal rational curves (on the original side and flopped side) and their strict transforms (or their images) with the chosen $\mathbb{Z}$-basis of $\mathrm{Pic}(T)$,

\item the KKMR decomposition,

\item the types of flops we need to ``get'' all the minimal models.

\end{enumerate}

However, in the item (4) above, the meaning of the word ``get'' is obscure, and the classification of the types \cite{WGBTMM95} is misleading in the way we describe below.  We would like to present the clarification of the types of flops listed in the item (4).

\subsection{\bf Origin of the naming of the types of flops}: The naming ($E_1$), ($E_2$), ($E_3$), ($E_4$), ($E_5$) of the types of flops originates from the analysis of the first flop

$$\begin{array}{rcccl}
 \mathit{Spec} \oplus_{m \geq 0} \mathit{Sym}^m(K_T^{\vee}) = Y && \dasharrow && Y^{+} \\
&\varphi \searrow && \swarrow & \\
&& X && \\
\end{array}$$

which corresponds to the divisorial contraction of the Fano 3-fold $T$ classified by Mori \cite{Mori82} labeled accordingly as ($E_1$), ($E_2$), ($E_3$), ($E_4$), ($E_5$).  Of course the Fano 3-fold $T$ is identified with the zero section embedded in the 4-fold $Y$, which is the total space of the canonical bundle of $T$.

\subsection{\bf Theorem III-2-1}  Theorem III-2-1 generalizes the above situation to the flopping contraction $\varphi':Y' \rightarrow X'$, which, when restricted to the reference 3-fold $T' \subset Y'$, gives rise to the divisorial contraction classified by Mori.  In our situation, the reference 3-fold $T'$ is not necessarily the zero section as above (embedded in the total space of the canonical bundle of $T$), but possibly its strict transform or some other divisor in $Y'$ (after some flops).  We also introduce the types ($F$) and ($G$), which do not directly arise from the divisorial contraction of some reference 3-fold.

\subsection{\bf Clarification of the meaning of the word ``get''} In \cite{WGBTMM95}, we say ``We give the classification of the types of flops we need to ``get'' all the minimal models.''  In order to clarify the meaning, we should have said ``We give the classification of the types of flops we need to ``connect'' all the minimal models, where the word ``connect'' means going from one minimal model to another whose nef cones (Weyl chambers) share codimension-one face.''

In this sense, we need the flops of the types ($E_1$), ($E_2$), ($E_3$), ($E_4$), ($E_5$), and {\bf their inverses}, and the flops of type ($F$) and ($G$), and moreover, we need the flops of Type (Others), which we did not mention in THEOREM III-2-4 (See 5.5. below.).

In the item (4) of the list in \cite{WGBTMM95}, we forgot to mention the inverses of the flops of types ($E_1$), ($E_2$), ($E_3$), ($E_4$), ($E_5$), while the flops of types ($F$) and ($G$) go both directions by definition and hence need no mention of their inverses.  (We need no mention of the inverses for the flops of Type (Others).)

\subsection{\bf Labeling ($F$) in the item (4)} The description of the flops of type ($F$) on Page 34 of \cite{WGBTMM95} is quite general, and in fact includes the flops of type ($E_1$) and their inverses.  Therefore, it is baffling that we find the label ($F$) in the item (4) in \cite{WGBTMM95} only for some Fano 3-folds, especially together with the label ($E_1$), and that we do not always find the label ($F$) even when the label ($E_1$) is listed in the item (4).  We need the clarification of the label ($F$) appearing in the item (4) in \cite{WGBTMM95}.

Actually the flops which are labeled of type ($F$) in the item (4) in \cite{WGBTMM95} are of a very specific type as described below, and indeed form a special case of the more general type ($F$) described on Page 34 of \cite{WGBTMM95}.

\vskip.1in

({\bf Description of the flops labeled as type ($F$) in the item (4)})

\vskip.1in

Here we give an explanation of the flops labeled as type ($F$) in the item (4), first dealing with the case of the Fano 3-fold $T$ denoted as $n^o\ 31$ with $B_2 = 3$.  The Fano 3-fold $T$ is the blow up of the cone over a smooth quadric surface $S \cong \mathbb{P}^1 \times \mathbb{P}^1 \subset \mathbb{P}^3$ with center the vertex, i.e., $T = \mathbb{P}_{\mathbb{P}^1 \times \mathbb{P}^1}(\mathit{O} \oplus \mathit{O}(1,1))$.  We have the following three extremal rays:

\begin{itemize}

\item[$l_1$:] the one generated by one ruling in the exceptional divisor ($\cong \mathbb{P}^1 \times \mathbb{P}^1$) of the blow up over the vertex

\item[$l_2$:] the one generated by the other ruling in the exceptional divisor ($\cong \mathbb{P}^1 \times \mathbb{P}^1$) of the blow up over the vertex

\item[$l_3$:] the one generated by the fiber of the $\mathbb{P}^1$-bundle $T = \mathbb{P}_{\mathbb{P}^1 \times \mathbb{P}^1}(\mathit{O} \oplus \mathit{O}(1,1)) \rightarrow \mathbb{P}^1 \times \mathbb{P}^1$.

\end{itemize}

Then the KKMR decomposition is given as follows:

\newpage

\includegraphics[width=5in]{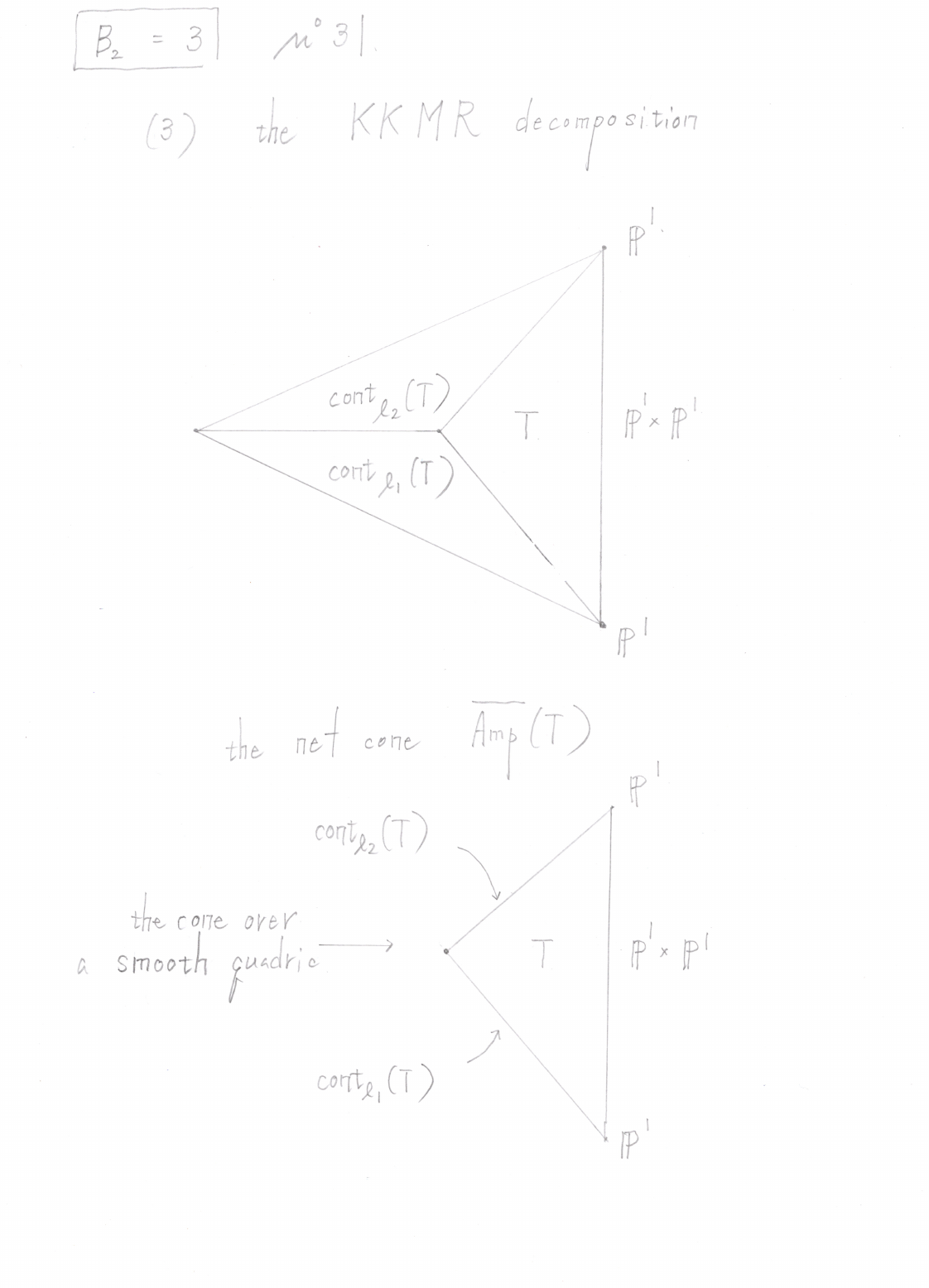}

Going from the chamber labeled as $T$ to $\mathrm{cont}_{l_1}(T)$ is a flop of type $(E_1)$, and going from the chamber labeled as $T$ to $\mathrm{cont}_{l_2}(T)$ is another flop of type $(E_1)$.

Going from the chamber labeled as $\mathrm{cont}_{l_1}(T)$ to the one labeled as $\mathrm{cont}_{l_2}(T)$ (and also going the other way around) is the flop labeled as ($F$) in the item (4), which is the composition of the inverse of a flop of type ($E_1$) and another flop of type ($E_1$).

Secondly we give a more precise description of the flops labeled as type ($F$) \text{in the item (4).}

We only describe the flop labeled as type ($F$) in the item (4) in the case $f : Y = \mathit{Spec}\ \oplus_{m \geq 0}\mathit{Sym}^m(K_T^{\vee}) \rightarrow X = \mathrm{Spec}\ \oplus_{m \geq 0}H^0(T,\mathit{Sym}^m(K_T^{\vee}))$ is the contraction of the zero section of the total space of the canonical bundle of a Fano 3-fold $T$.  The zero section is identified with $T$.  (The description in the other cases is identical to the one in this case.  Or rather to say, the description is not using any condition specific to this particular case.)

Inside of a Fano 3-fold $T$ we have a surface $S \cong \mathbb{P}^1 \times \mathbb{P}^1 \subset T$.  Assume that the two distinct rulings of $S \cong \mathbb{P}^1 \times \mathbb{P}^1$ give rise to two distinct extremal rays $l_1$ and $l_2$, each of which is of type ($E_1$).  The contraction of the extremal ray $l_i \ (i = 1,2)$ for the Fano 3-fold $T$ gives rise to the corresponding flopping contraction of the 4-fold $\varphi_i:Y \rightarrow Z_i\ (i = 1,2)$.  The flopped morthism $\varphi_i^+:Y_i^+ \rightarrow Z_i\ (i = 1,2)$ is obtained as follows.  First blow up $S$ inside of the 4-fold $Y$ to have the morphism $\mathrm{Blp}_S: Y' \rightarrow Y$, whose exceptional divisor is denoted by $D$.  Then contract $D$ in the other direction (than $\mathrm{Blp}_S$) to obtain $(\mathrm{Bldn}_D)_i: Y' \rightarrow Y_i^+$.  The induced morphism $\varphi_i^+:Y_i^+ \rightarrow Z_i$ is the flopped morphism.

$$\begin{array}{rcccccccr}
&&&& Y' &&&& \\
&&&&&&&& \\
&&&& \mathrm{Blp}_S &&&& \\
&& \swarrow (\mathrm{Bldn}_D)_1  && \downarrow  \hskip.04in && (\mathrm{Bldn}_D)_2 \searrow && \\
&&&& \\
Y_1^+ && \overset{\psi_1}\dashleftarrow &&  Y && \overset{\psi_2}\dashrightarrow && \hskip.2in Y_2^+ \\
&&&&&&&& \\
&\varphi_1^+ \searrow && \swarrow \varphi_1 && \varphi_2 \searrow && \swarrow \varphi_2^+ & \\
&&&&&&&& \\
&& Z_1 &&&& Z_2 && \\
\end{array}$$

Each flop $\psi_i: Y \dashrightarrow Y_i^+\ (i = 1,2)$ is of type ($E_1$).

Now one can further blow down $(\mathrm{Bldn}_D)_i$ in the other direction of the two rulings to obtain $\sigma_i: Y_i^+ \rightarrow W$, where $W$ is common for $i = 1, 2$.  The resulting 4-fold $W$ is equi-singular along the curve (See  Remark (Atiyah's flop) below.) corresponding to the fiber of the blow up $\mathrm{Blp}_S$.

$$\begin{array}{ccccc}
&& Y' && \\
&&&& \\
& \swarrow (\mathrm{Bldn}_D)_1  &  & (\mathrm{Bldn}_D)_2 \searrow & \\
&&&& \\
Y_1^+ && \dashrightarrow && Y_2^+ \\
&&&& \\
&\sigma_1 \searrow && \swarrow \sigma_2 & \\
&&&& \\
&& W && \\
\end{array}$$

The obtained flop $Y_1^+ \dashrightarrow Y_2^+$ (or symmetrically $Y_2^+ \dashrightarrow Y_1^+$) is the flop labeled as type ($F$) in the item (4), which is the composition of the inverse of the flop $\psi_1$ of type ($E_1$) and the flop $\psi_2$ also of type ($E_1$) )or symmetrically the composition of the inverse of the flop $\psi_2$ of type ($E_1$) and the flop $\psi_1$ also of type ($E_1$).  Actually $(\mathrm{Bldn}_D)_2 \circ \{(\mathrm{Bldn}_D)_1\}^{-1}$ is the blow up and down  appearing in the description the flop of type ($F$) on Page 34 of \cite{WGBTMM95}.  Note that the 4-fold $W$ does not give rise to the full dimensional chamber, but corresponds to the common codimension-one face shared by $\overline{\mathrm{Amp}}(Y_1^+/X)$ and $\overline{\mathrm{Amp}}(Y_2^+/X)$ in the KKMR decomposition.

\vskip.1in

Remark (Atiyah's flop (cf. Page 157 in \cite{Matsuki02})): We should emphasize again that the flops of type ($E_1$) and their inverses, and the flops described above labeled as type ($F$) in the item (4) are all special cases of the flops of type ($F$) described on Page 34 of \cite{WGBTMM95}, which are all analytically isomorphic to $(\mathrm{Atiyah's\ flop}) \times \mathbb{A}^1$.

\subsection{\bf Correction to THEOREM III-2-4} The statement of the theorem is correct as far as it claims the following: Let $T$ be a Fano 3-fold, and $f: Y := \mathit{Spec}\oplus_{m \geq 0}\mathit{Sym}(K_T^{\vee}) \rightarrow X := \mathrm{Spec}\oplus_{m \geq 0}H^0(T,\mathit{Sym}^m(K_T^{\vee})$ the natural morphism.  Then all the minimal models of $Y$ over $X$ are obtained from $Y$ by a sequence of flops of Types \text{($E_1$), ($E_2$), ($E_3$), ($E_4$), ($E_5$).}

However, if we use the word ``connect two minimal models'' in a more restrictive way, meaning that going from one minimal model to another whose nef cones (Weyl chambers) share codimension-pne face (see 5.3) , then we need not only the flops of Types ($E_1$), ($E_2$), ($E_3$), ($E_4$), ($E_5$) and their inverses, flops of Types ($F$) and ($G$), but also {\bf flops of Types (Others)}, which we did not mention in the statement of \text{THEOREM III-2-4 in \cite{WGBTMM95}.}

The flops of Types ($F$), ($G$) and (Others) modify the loci away from the strict transform of the zero section.

The flops of Types (Others) appear with Fano 3-folds with $\boxed{B_2 = 3}$ $n^o\ 11$, $n^o\ 16$, $n^o\ 18$, $n^o\ 23$, $n^o\ 30$.  Each one of flops of Types (Others) is obtained as a sequence of flops of Types  ($E_1$), ($E_2$), ($E_3$), ($E_4$), ($E_5$), and can be constructed as such.  Therefore, we omit the description of the construction, as we did with the flops of Types ($F$) and ($G$).

\newpage

\section{Mistakes found in the tables of the intersection pairings and the corrections}

\vskip.1in

$\boxed{B_2 = 2}$

\vskip.03in

\underline{$n^o\ 1$.}

\vskip.03in

About (2).

\vskip.03in

\hskip.1in Mistake:

\begin{center}
\begin{tabular}{|c|c|c|c|}
\hline
{\ \ \ \ \ }   &  $E$ & $Blp^*(- \frac{1}{2}K_V)$ & $- K_T$ \\
\hline
\end{tabular}
\end{center}

\hskip.1in Correction:

\begin{center}
\begin{tabular}{|c|c|c|c|}
\hline
{\ \ \ \ \ }   &  $E$ & $Blp^*(- \frac{1}{2}K_{V_1})$ & $- K_T$ \\
\hline
\end{tabular}
\end{center}

\vskip.1in

Note: How to compute the numbers in the last two rows in the table, associated with the 4-fold flop of type ($E_1$).

\vskip.1in

In the following, we explain how to compute the numbers in the last two rows in the table, associated with the 4-fold flop of type ($E_1$), using its explicit construction given in Page 31 of \cite{WGBTMM95}.  (See the attached picture.)

\newpage

\includegraphics[width=5in]{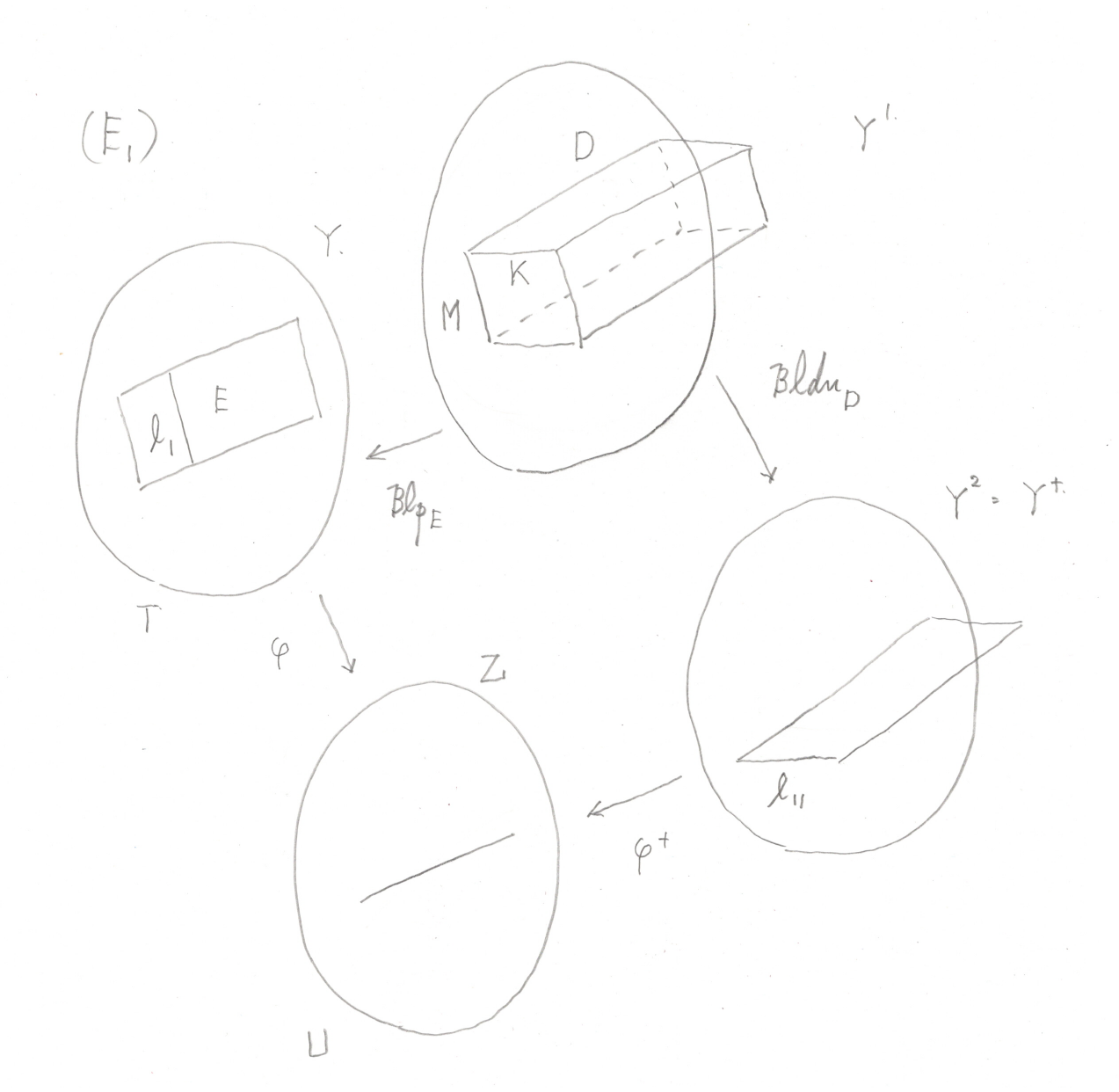}

Let $\varphi:Y \rightarrow Z$ be the contraction of flopping type on the 4-fold $Y$.  Then $\varphi|_T$ is the divisorial contraction (of the extremal ray $l_1$) $\phi:T \rightarrow U$, which contracts the exceptional divisor $E = E_{\varphi}$ onto a smooth curve $\phi|_E: E \rightarrow \phi(E)$ giving the rulng for $E$.  To obtain the flop of the small contraction $\varphi:Y \rightarrow Z$, we first blow up $E$ in the 4-fold $Blp_E: Y^1 \rightarrow Y$.  It is easy to see that the normal bundle $\mathcal{N}_{E/Y}$ restricted to any ruling ($\cong \mathbb{P}^1$) is isomorphic to $\mathcal{O}_{\mathbb{P}^1}(-1) \oplus \mathcal{O}_{\mathbb{P}^1}(-1)$.  Thus the exceptional divisor $D$ for $Blp_E$ has normal bundle $\mathcal{N}_{D/Y^1}$ restricting to $\mathcal{O}(-1,-1)$ on each inverse image under $Blp_E$ of any ruling, which is isomorphic to $\mathbb{P}^1 \times \mathbb{P}^1$.  Thus $D$ can be blown down in the different direction to get $Bldn_D: Y^1 \rightarrow Y^2$.  Now the image $Bldn_D(D)$ is a ruled surface on $Y^2$, which can be contracted to get $Z$.  The resulting morphism $\varphi^+:Y^2 = Y^+ \rightarrow Z$ is the desired flop.

We need to give some explanation of the last bottom two rows of the table.  The symbol $l_{11}$ denotes the image curve of $\mathbb{P}^1 \times \mathbb{P}^1$ under the second projection, where the fisrt projection corresponds to $Blp_E$ and where the second projection corresponds to $Bldn_D$.  The symbol $l_{21}$ denotes the strict transform of the extremal rational curve denoted by $l_2$.  The symbol $E$ and $Blp^*(- \frac{1}{2}K_{V_1})$ in the top row actually refer to $p^{-1}(E)$ and $p^{-1}(Blp^*(- \frac{1}{2}K_{V_1}))$, where $p:Y \rightarrow T$ is the natural map from $Y$ regarded as the total space of the canonical bundle $K_T$ to $T$ as above.  Therefore, in the last bottom two rows of the table, we are supposed to compute the intersection numbers of $l_{11}$ and $l_{21}$ with $\{p^{-1}(E)\}^+$ and $\{p^{-1}(Blp^*(- \frac{1}{2}K_{V_1}))\}^+$, the strict tarnsforms of the divisors $p^{-1}(E)$ and $p^{-1}(Blp^*(- \frac{1}{2}K_{V_1}))$ on $Y^+$, which we denote simply by $E^+$ and $Blp^*(- \frac{1}{2}K_{V_1})^+$ by abuse of notation.

In the above sense, we would like to compute the intersection numbers

$$\begin{array}{ll}
(l_{11},E^+) & (l_{11},Blp^*(- \frac{1}{2}K_{V_1})^+) \\
(l_{21},E^+) & (l_{21},Blp^*(- \frac{1}{2}K_{V_1})^+). \\
\end{array}$$

The intersection number with $(- K_T)^+$ can be easily computed from the above using the formula
$$K_T = Blp^*K_{V_1} + E,$$
which implies
$$(- K_T)^+ = 2 Blp^*(- \frac{1}{2}K_{V_1})^+ - E^+.$$
Let $S$ be a member of $|- \frac{1}{2}K_{V_1}|$.  Then the adjunction formula gives
$$K_S = (S + K_{V_1})|_S = \frac{1}{2}K_{V_1}|_S,$$
which implies
$$(Blp(l_2),\frac{1}{2}K_{V_1})_{V_1} = (Blp(l_2),K_S)_S = -1.$$

Our strategy is to compute these intersection numbers, using the intersection pairings between some appropriate curves and divisors on $Y^2$.

Let $\psi: Y^2 \rightarrow Y$ and $\psi^+:Y^2 \rightarrow Y^+$ denote the corresponding morphisms.

We set
$$\left\{\begin{array}{rcl}
M &:& \text{a fiber of the second projection }\mathbb{P}^1 \times \mathbb{P}^1 \rightarrow \mathbb{P}^1 \\
K &:& \text{a fiber of the first projection }\mathbb{P}^1 \times \mathbb{P}^1 \rightarrow \mathbb{P}^1 \\
l_{21} &:& \text{the strict transform of }l_2 \text{ on }Y^2 \text{ by abuse of notation} \\
\end{array}\right.$$

We have
$$\left\{\begin{array}{lcl}
(\psi^+)^*E^+ &=& \psi^*E + aD \\
(\psi^+)^*Blp^*(- \frac{1}{2}K_{V_1})^+ &=& \psi^*Blp^*(- \frac{1}{2}K_{V_1}) + cD \\
&& \text{ for some numbers }a,c. \\
\end{array}\right.$$

In order to determine the numbers $a,c$, we use the conditions

$$\left\{\begin{array}{lclcl}
(M,(\psi^+)^*E^+) &=& (M,\psi^*E + aD) &=& 0,\\
(M,(\psi^+)^*Blp^*(- \frac{1}{2}K_{V_1})^+) &=& (M,\psi^*Blp^*(- \frac{1}{2}K_{V_1}) + cD) &=& 0.\\
\end{array}\right.$$

Since
$$\left\{\begin{array}{ll}
(M,\psi^*E) = -1, & (M,D) = -1, \\
(M,\psi^*Blp^*(- \frac{1}{2}K_{V_1})) = 0, & (M,D) = -1,\\
\end{array}\right.$$

we conclude
$$a = -1, c = 0.$$

Finally we compute

$$\left\{\begin{array}{rcl}
(l_{11},E^+) &=& (K,\psi^*E - 1D) \\
&=& 0 - 1\cdot (-1) = 1 \\
(l_{11},Blp^*(- \frac{1}{2}K_{V_1})^+) &=& (K,\psi^*Blp^*(- \frac{1}{2}K_{V_1}) + 0D) \\
&=& 0 + 0 \cdot (-1) = 0 \\
\end{array}\right.$$ 

since
$$(K,\psi^*E) = -1, (K,\psi^*Blp^*(- \frac{1}{2}K_{V_1})) = 0, (K,D) = -1.$$

Similarly we compute

$$\left\{\begin{array}{rcl}
(l_{21},E^+) &=& (l_{21},\psi^*E - 1D) \\
&=& 1 -1  \cdot 1 = 0 \\
(l_{21},Blp^*(- \frac{1}{2}K_{V_1})^+) &=& (l_{21},\psi^*Blp^*(- \frac{1}{2}K_{V_1}) + 0D) \\
&=& 1 + 0 \cdot 1 = 1 \\
\end{array}\right.$$

since
$$(l_{21},\psi^*E) = 0, (l_{21},\psi^*Blp^*(- \frac{1}{2}K_{V_1})) = 1, (l_{21},D) = 1.$$

Accordingly, we have
$$\left\{\begin{array}{l}
(l_{11},(- K_T)^+) = (l_{11},2 Blp^*(- \frac{1}{2}K_{V_1})^+ - E^+\}) = 2 \cdot 0 - 1 = - 1 \\
(l_{21},(- K_T)^+) = (l_{21},2 Blp^*(- \frac{1}{2}K_{V_1})^+ - E^+) = 2 \cdot 1 - 0 = 2.\\
\end{array}\right.$$

This completes the explanation of how to compute the numbers in the last two rows in the table.

\vskip.03in

\underline{$n^o\ 2$.}

\vskip.03in

\hskip.1in NO Mistake.

\vskip.03in

\underline{$n^o\ 3$.}

\vskip.03in

About (2).

\vskip.03in

\hskip.1in Mistake:

\begin{center}
\begin{tabular}{|c|c|c|c|}
\hline
{\ \ \ \ \ }   &  $E$ & $Blp^*(- \frac{1}{2}K_V)$ & $- K_T$ \\
\hline
\end{tabular}
\end{center}

\hskip.1in Correction:

\begin{center}
\begin{tabular}{|c|c|c|c|}
\hline
{\ \ \ \ \ }   &  $E$ & $Blp^*(- \frac{1}{2}K_{V_2})$ & $- K_T$ \\
\hline
\end{tabular}
\end{center}

\vskip.03in

\underline{$n^o\ 4$.}

\vskip.03in

\hskip.1in NO Mistake.

\vskip.03in

\underline{$n^o\ 5$.}

\vskip.03in

\hskip.1in NO Mistake.

\vskip.03in

\underline{$n^o\ 6$ (6.A).}

\vskip.03in

\hskip.1in NO Mistake.

\vskip.03in

\underline{$n^o\ 6$ (6.B).}

\vskip.03in

About the description.

\vskip.03in

\hskip.1in Mistake: $W_6$ is a smooth divisor Of bidegree ...

\hskip.1in Correction: $W_6$ is a smooth divisor of bidegree ...

\hskip.1in Mistake: $l_1$: an irreducible component of a reducible fiber ora reduced part of ...

\hskip.1in Correction: $l_1$: an irreducible component of a reducible fiber or a reduced part of ...

\vskip.03in

\hskip.1in NO Mistake in the table of the intersection pairings.

\vskip.03in

\underline{$n^o\ 7$.}

\vskip.03in

\hskip.1in NO Mistake.

\vskip.03in

\underline{$n^o\ 8$.}

\vskip.03in

About (2).

\vskip.03in

\hskip.1in Mistake:

\begin{center}
\begin{tabular}{|c|c|c|c|}
\hline
  &  $d^*D$ & $d^*Blp^*\mathcal{O}_{\mathbb{P}^3}(1)$ & $- K_T$ \\
\hline
$l_1$ & $- 1$ & 1 & 1 \\
\hline
$l_2$ & $- 1$ & 0 & 1 \\
\hline
$l_{11}$ & 0 & 1 & 2 \\
\hline
$l_{21}$ & 1 & 0 & $- 1$ \\
\hline
\end{tabular}
\end{center}

\vskip.03in

Note: The subscripts for the last two lines $l_{11}$ and $l_{21}$ should have been $l_{12}$ and $l_{22}$, respectively.

\vskip.03in

\hskip.1in Correction:

\begin{center}
\begin{tabular}{|c|c|c|c|}
\hline
  &  $d^*D$ & $d^*Blp^*\mathcal{O}_{\mathbb{P}^3}(1)$ & $- K_T$ \\
\hline
$l_1$ & 1 & 1 & 1 \\
\hline
$l_2$ & $- 1$ & 0 & 1 \\
\hline
$l_{12}$ & 0 & 1 & 2 \\
\hline
$l_{22}$ & 1 & 0 & $- 1$ \\
\hline
\end{tabular}
\end{center}

\vskip.03in

Note: How to compute the numbers in the last two rows in the table, associated with the 4-fold flop of type ($E_3$) or ($E_4$) .

\vskip.03in

In the following, we explain how to compute the numbers in the last two rows in the table, associated with the 4-fold flop of type ($E_3$) or ($E_4$), using its explicit construction given in Page 32 of \cite{WGBTMM95}.  (See the attached picture.)

\includegraphics[width=3.5in]{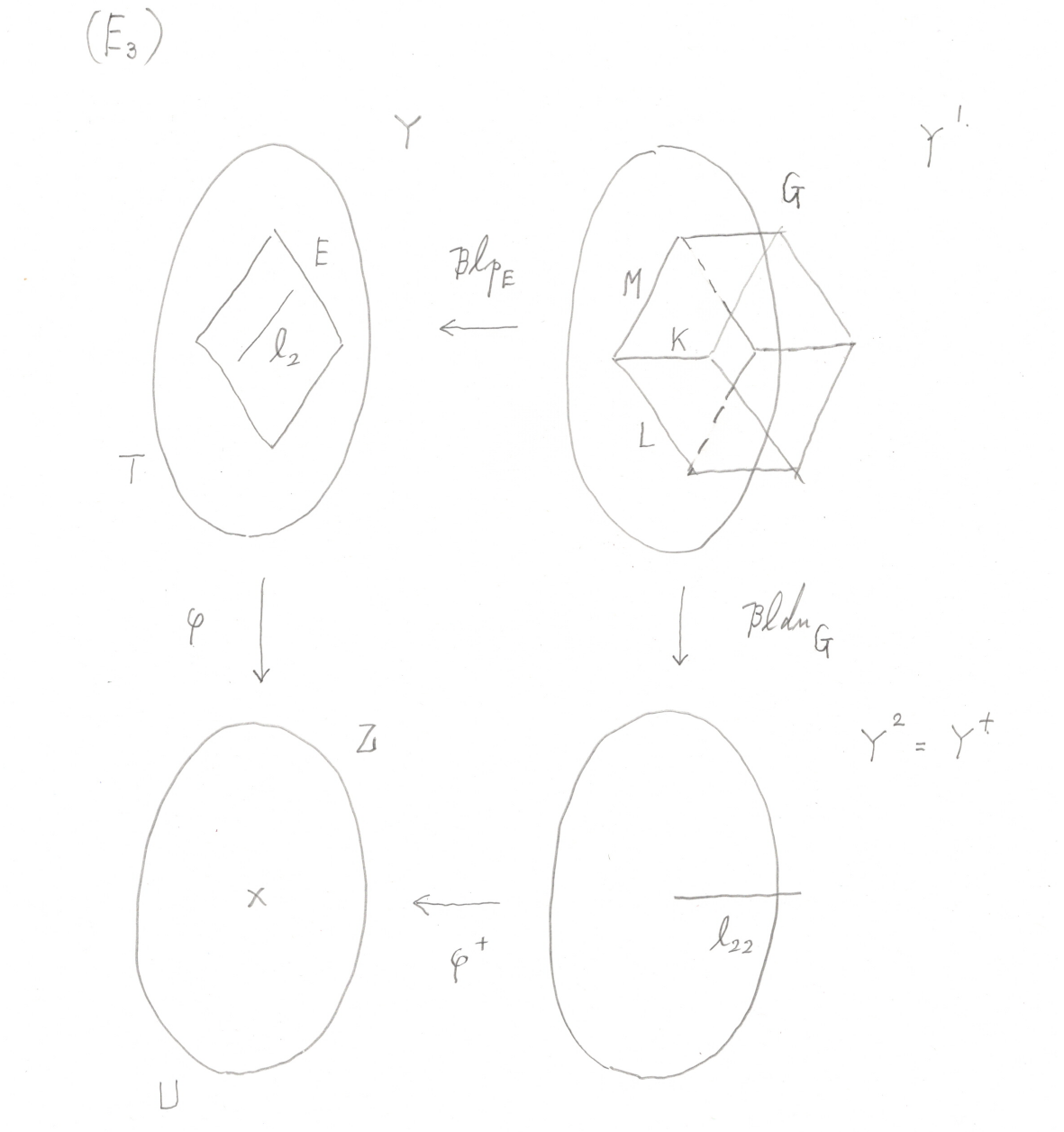}\includegraphics[width=3.5in]{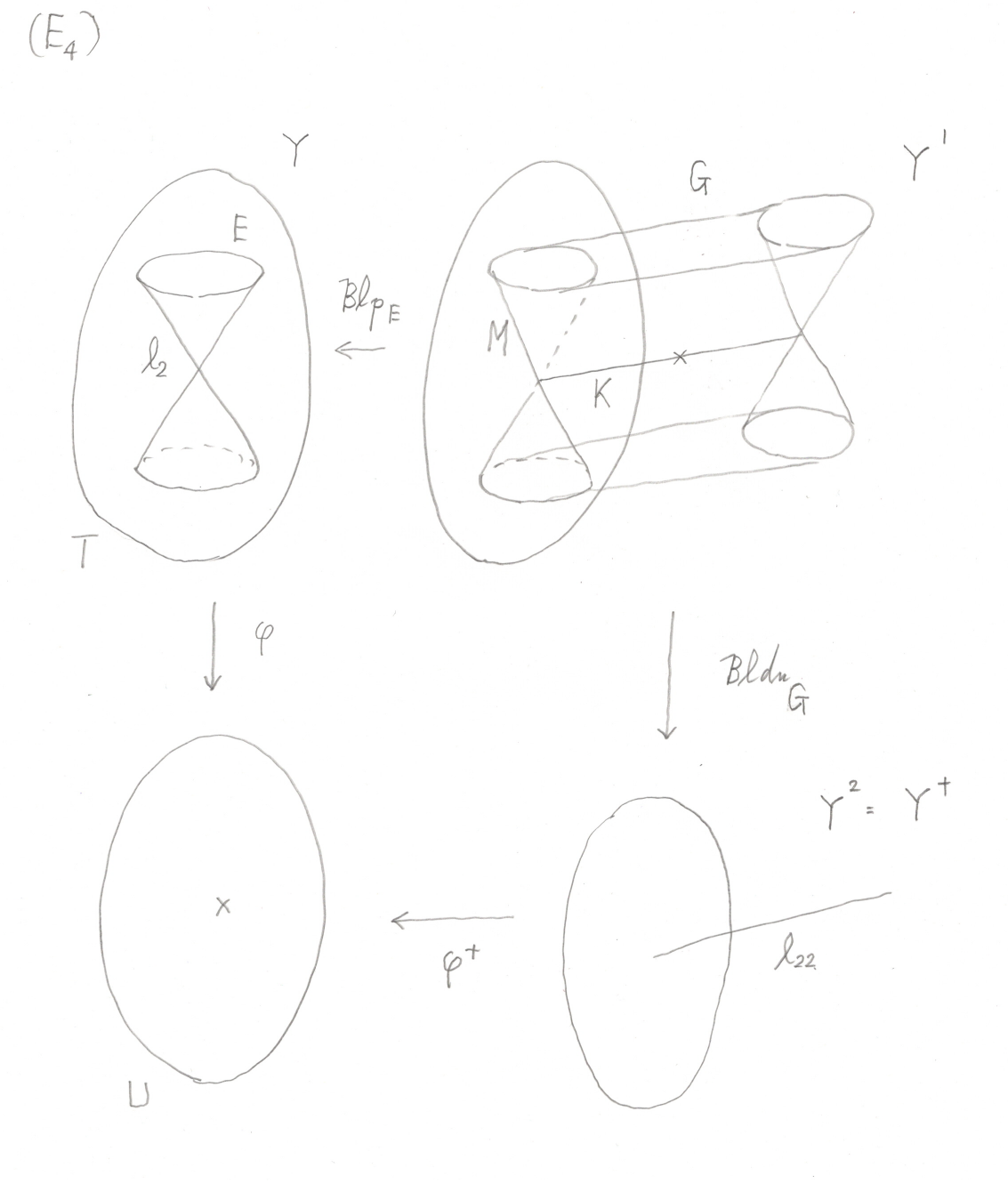}

Since the computation for type ($E_4$) (the degenerate case of type ($E_3$)) is identical ot the comotation for type ($E_3$), we only present the computation for type ($E_3$) below.

Let $\varphi:Y \rightarrow Z$ be the contraction of flopping type on the 4-fold $Y$.  Then $\varphi|_T$ is the contraction (of extremal ray $l_2$) $\phi: T \rightarrow U$, which contracts the exceptional divisor $E = E_{\varphi} \cong \mathbb{P}^1 \times \mathbb{P}^1$ with normal bundle $\mathcal{N}_{E/T} \cong \mathcal{O}(-1,-1)$ onto an ordinary double point of $U$.  (Note that $E = d^*D$.)  To obtain the flop of the small contraction $\varphi: Y \rightarrow Z$, we first blow up $E$ in $Y$.  Since $E$ is the complete intersection of $T$ (the zero section) with $V$, the normal bundle is the direct sum $\mathcal{N}_{E/Y} \cong \mathcal{O}(-1,-1) \oplus \mathcal{O}(-1,-1)$  (Note that $V = p^{-1}(E)$ is the divisor in $Y$, where $p:Y \rightarrow T$ is the natural map from $Y$ regarded as the total space of the canonical bundle $K_T$ to $T$.  This divisor $V$ is denoted by $T_1$ in Page 32 of of \cite{WGBTMM95}.)  Thus the blow up $Blp_E: Y^1 \rightarrow Y$ has the exceptional divisor $G \cong \mathbb{P}^1 \times \mathbb{P}^1 \times \mathbb{P}^1$ with normal bundle $\mathcal{N}_{G/Y_1} \cong \mathcal{O}_{E \times \mathbb{P}^1}(-1,-1,-1)$.  Now $G$ can be contracted $Bldn_G:Y^1 \rightarrow Y^2$ corresponding to the projection $E \times\mathbb{P}^1 \rightarrow \mathbb{P}^1$.  The target space $Y^2$ is singular along the image of $G$.  Then we can contract this singular locus ($\cong \mathbb{P}^1$) to get the desired flop $\varphi^+:Y^2 = Y^+ \rightarrow Z$.  (Note that on Page 32 of \cite{WGBTMM95} the exceptional divisor $G$ here is denoted by the symbol $D$.  However, on Page 54 of \cite{WGBTMM95} the symbol $D$ is already used for the exceptional divisor for $V_7 \rightarrow \mathbb{P}^3$.  We use the different symbol $G$ here to avoid any confusion.)

We need to give some explanation of the last bottom two rows of the table.  The symbol $l_{12}$ denotes the strict transform of the extremal rational curve denoted by $l_1$, while $l_{22}$ denotes the image curve of $D$, which is the exceptional locus of $\varphi^+$.  The symbol $d^*D$ and $d^*Blp^*\mathcal{O}_{\mathbb{P}^3}(1)$ in the top row actually refer to $p^{-1}(d^*D)$ and $p^{-1}(d^*Blp^*\mathcal{O}_{\mathbb{P}^3}(1))$, wheer $p:Y \rightarrow T$ is the natural map from $Y$ regarded as the total sopace of the canonical bundle $K_T$ to $T$ as above.  Therefore, in the last bottom two rows of the table, we are supposed to comoute the intersection numbers of $l_{12}$ and $l_{22}$ with $\{p^{-1}(d^*D)\}^+$ and $\{p^{-1}(d^*Blp^*\mathcal{O}_{\mathbb{P}^3}(1))\}^+$, the strict transforms of the divisors $p^{-1}(d^*D)$ and $p^{-1}(d^*Blp^*\mathcal{O}_{\mathbb{P}^3}(1))$ on $Y^+$, which we simply denote by $d^*D^+$ and $d^*Blp^*\mathcal{O}_{\mathbb{P}^3}(1)^+$ by abuse of notation.

In the above sense, we would like to compute the intersection numbers
$$\begin{array}{ll}
(l_{12},d^*D^+) & (l_{12},d^*Blp^*\mathcal{O}_{\mathbb{P}^3}(1)^+) \\
(l_{22},d^*D^+) & (l_{22},d^*Blp^*\mathcal{O}_{\mathbb{P}^3}(1)^+). \\
\end{array}$$

The intersection number with $(- K_T)^+$ can be easily computed from the above using the formula
$$\begin{array}{rcl}
K_{V_7} &=& -4 Blp^*\mathcal{O}_{\mathbb{P}^3}(1) + 2D \\
K_T &=& d^*K_{V_7} + d^*(\frac{1}{2}|- K_{V_7}|) \\
&=& d^*(-4 Blp^*\mathcal{O}_{\mathbb{P}^3}(1) + 2D) + d^*(2 Blp^*\mathcal{O}_{\mathbb{P}^3}(1) - D) \\
&=& d^*(-2 Blp^*\mathcal{O}_{\mathbb{P}^3}(1) + D),\\
\end{array}$$
which implies
$$(- K_T)^+ = 2 d^*Blp^*\mathcal{O}_{\mathbb{P}^3}(1)^+ - d^*D^+.$$

Our strategy is to compute these intersection numbers, using the intersection pairings between some appropriate curves and divisors on $Y^1$.

Let $\psi = Blp_E: Y^1 \rightarrow Y$ and $\psi^+ = Bldn_G:Y^1 \rightarrow Y^+$ denote the corresponding morphisms.  Let $T^1$ be the strict transform of $T$ (as the divisor on $Y$) on $Y^1$.  We set

$$\left\{\begin{array}{rcl}
M &:& \text{a fiber of the first projection }G \cap T^1 \cong E \cong \mathbb{P}^1 \times \mathbb{P}^1 \rightarrow \mathbb{P}^1 \\
L &:& \text{a fiber of the second projection }G \cap T^1 \cong E \cong \mathbb{P}^1 \times \mathbb{P}^1 \rightarrow \mathbb{P}^1 \\
K &:& \text{a fiber of }G \cong E \times \mathbb{P}^1 \rightarrow E \\
l_{12} &:& \text{the strict transform of }l_1 \text{ on }Y^1 \text{ by abuse of notation} \\
\end{array}\right.$$

We have
$$\left\{\begin{array}{lcl}
(\psi^+)^*d^*D^+ &=& \psi^*d^*D + aD \\
(\psi^+)^*d^*Blp^*\mathcal{O}_{\mathbb{P}^3}(1)^+ &=& \psi^*d^*Blp^*\mathcal{O}_{\mathbb{P}^3}(1) + cD \\
&& \text{ for some numbers }a,c. \\
\end{array}\right.$$

In order to determine the numbers $a,c$, we use the conditions

$$\left\{\begin{array}{lclcl}
(M,(\psi^+)^*d^*D^+) &=& (M,\psi^*d^*D + aG) &=& 0,\\
(M,(\psi^+)^*d^*Blp^*\mathcal{O}_{\mathbb{P}^3}(1)^+) &=& (M,\psi^*d^*Blp^*\mathcal{O}_{\mathbb{P}^3}(1) + cG) &=& 0.\\
\end{array}\right.$$

(Remark that the conditions
$$\left\{\begin{array}{lclcl}
(L,(\psi^+)^*d^*D^+) &=& (L,\psi^*d^*D + aG) &=& 0,\\
(L,(\psi^+)^*d^*Blp^*\mathcal{O}_{\mathbb{P}^3}(1)^+) &=& (L,\psi^*d^*Blp^*\mathcal{O}_{\mathbb{P}^3}(1) + cG) &=& 0.\\
\end{array}\right.$$
follow from the above in the case of type ($E_3$).)

Since
$$\left\{\begin{array}{ll}
(M,\psi^*d^*D) = -1, & (M,G) = -1, \\
(M,\psi^*d^*Blp^*\mathcal{O}_{\mathbb{P}^3}(1)) = 0, & (M,G) = -1,\\
\end{array}\right.$$

we conclude
$$a = -1, c = 0.$$

Finally we compute

$$\left\{\begin{array}{rcl}
(l_{12},d^*D^+) &=& (l_{12},\psi^*d^*D - 1G) \\
&=& 1 - 1\cdot 1 = 0 \\
(l_{12},d^*Blp^*\mathcal{O}_{\mathbb{P}^3}(1)^+) &=& (l_{12},\psi^*d^*Blp^*\mathcal{O}_{\mathbb{P}^3}(1) + 0G) \\
&=& 1 + 0 \cdot 1 = 0 \\
\end{array}\right.$$ 

since
$$(l_{12},\psi^*d^*D) = 1, (l_{12},\psi^*d^*Blp^*\mathcal{O}_{\mathbb{P}^3}(1)) = 1, (l_{12},G) = 1.$$

Similarly we compute

$$\left\{\begin{array}{rcl}
(l_{22},d^*D^+) &=& (l_{22},\psi^*d^D - 1G) \\
&=& 0 - 1  \cdot (-1) = 1 \\
(l_{22},d^*Blp^*\mathcal{O}_{\mathbb{P}^3}(1)^+) &=& (K,\psi^*d^*Blp^*\mathcal{O}_{\mathbb{P}^3}(1) + 0G) \\
&=& 0 + 0 \cdot (-1) = 0 \\
\end{array}\right.$$

since
$$(K,\psi^*d^*D) = 0, (K,\psi^*Blp^*\mathcal{O}_{\mathbb{P}^3}(1)) = 0, (K,G) = -1.$$

Accordingly, we have
$$\left\{\begin{array}{l}
(l_{12},(- K_T)^+) = (l_{12},2 d^*Blp^*\mathcal{O}_{\mathbb{P}^3}(1)^+ - d^*D^+) = 2 \cdot 1 - 0 = 2\\
(l_{22},(- K_T)^+) = (l_{22},2 d^*Blp^*\mathcal{O}_{\mathbb{P}^3}(1)^+ - d^*D^+) = 2 \cdot 0 - 1 = - 1.\\
\end{array}\right.$$

This completes the explanation of how to compute the numbers in the last two rows in the table.

\underline{$n^o\ 9$.}

\vskip.03in

\hskip.1in NO Mistake.

\vskip.03in

\underline{$n^o\ 10$.}

\vskip.03in

\hskip.1in NO Mistake.

\vskip.03in

\underline{$n^o\ 11$.}

\vskip.03in

\hskip.1in NO Mistake.

\vskip.03in

\underline{$n^o\ 12$.}

\vskip.03in

\hskip.1in NO Mistake.

\vskip.03in

\underline{$n^o\ 13$.}

\vskip.03in

\hskip.1in NO Mistake.

\newpage

\underline{$n^o\ 14$.}

\vskip.03in

\hskip.1in NO Mistake.

\vskip.03in

\underline{$n^o\ 15$.}

\vskip.03in

About (2).

\vskip.03in

\hskip.1in Mistake:

\begin{center}
\begin{tabular}{|c|c|c|c|}
\hline
{\ \ \ \ \ }   &  $D$ & $Blp^*\mathcal{O}_{\mathbb{P}^3}(1)$ & $- K_T$ \\
\hline
\end{tabular}
\end{center}

\hskip.1in Correction:

\begin{center}
\begin{tabular}{|c|c|c|c|}
\hline
{\ \ \ \ \ }   &  $E$ & $Blp^*\mathcal{O}_{\mathbb{P}^3}(1)$ & $- K_T$ \\
\hline
\end{tabular}
\end{center}

\hskip1.5in $E$: the exceptional divisor of the blow up $Blp: T \rightarrow \mathbb{P}^3$

\vskip.03in

About (4).

\hskip.1in Mistake: It is missing ($E_1$).

(15.a) ($E_3$)

(15.b) ($E_4$)

\hskip.1in Correction:

($E_1$)

(15.a) ($E_3$)

(15.b) ($E_4$)

\vskip.03in

\underline{$n^o\ 16$.}

\vskip.03in

\hskip.1in NO Mistake.

\vskip.03in

\underline{$n^o\ 17$.}

\vskip.03in

\hskip.1in NO Mistake.

\vskip.03in

\underline{$n^o\ 18$.}

\vskip.03in

\hskip.1in NO Mistake.

\vskip.03in

\underline{$n^o\ 19$.}

\vskip.03in

About (2).

\vskip.03in

\hskip.1in Mistake:

\begin{center}
\begin{tabular}{|c|c|c|c|}
\hline
  &  $E$ & $Blp^*\mathcal{O}_{V_4}(1)$ & $- K_T$ \\
\hline
$l_1$ & $- 1$ & 0 & 1 \\
\hline
$l_2$ & 1 & 1 & 1 \\
\hline
$l_{11}$ & 1 & 0 & $- 1$ \\
\hline
$l_{21}$ & 0 & 1 & 2  \\
\hline
$l_{12}$ & 2 & 3 & 4 \\
\hline
$l_{21}$ & $- 1$ & $- 1$ & $- 1$  \\
\hline
\end{tabular}
\end{center}

\vskip.03in

Note: The subscript for the last line $l_{21}$ should have been $l_{22}$.

\newpage

\hskip.1in Correction:

\begin{center}
\begin{tabular}{|c|c|c|c|}
\hline
  &  $E$ & $Blp^*\mathcal{O}_{V_4}(1)$ & $- K_T$ \\
\hline
$l_1$ & $- 1$ & 0 & 1 \\
\hline
$l_2$ & 1 & 1 & 1 \\
\hline
$l_{11}$ & 1 & 0 & $- 1$ \\
\hline
$l_{21}$ & 0 & 1 & 2  \\
\hline
$l_{12}$ & 2 & 3 & 4 \\
\hline
$l_{22}$ & $- 1$ & $- 1$ & $- 1$  \\
\hline
\end{tabular}
\end{center}

\vskip.03in

\underline{$n^o\ 20$.}

\vskip.03in

About (1).

Mistake: $l_2$: ... $l_2$ ia an irreducible component ...

Correction: $l_2$: ... $l_2$ is an irreducible component ...

\vskip.03in

\hskip.1in NO Mistake in the table of the intersection pairings.

\vskip.03in

\underline{$n^o\ 21$.}

\vskip.03in

About (1).

Mistake: $l_2$: ... (which areise as the exceptional curves ...)

Correction: $l_2$: ... (which arise as the exceptional curves ...)

\vskip.03in

About (2).

\vskip.03in

\hskip.1in Mistake:

\begin{center}
\begin{tabular}{|c|c|c|c|}
\hline
  &  $E$ & $Blp^*\mathcal{O}_Q(1)$ & $- K_T$ \\
\hline
$l_1$ & $- 1$ & 0 & 1 \\
\hline
$l_2$ & 2 & 1 & 1 \\
\hline
$l_{11}$ & 1 & 0 & $- 1$ \\
\hline
$l_{21}$ & 0 & 1 & 3  \\
\hline
$l_{12}$ & 3 & 2 & 3 \\
\hline
$l_{21}$ & $- 2$ & $- 1$ & $- 1$  \\
\hline
\end{tabular}
\end{center}

\vskip.03in

Note: The subscript for the last line $l_{21}$ should have been $l_{22}$.

\vskip.03in

Correction:

\begin{center}
\begin{tabular}{|c|c|c|c|}
\hline
  &  $E$ & $Blp^*\mathcal{O}_Q(1)$ & $- K_T$ \\
\hline
$l_1$ & $- 1$ & 0 & 1 \\
\hline
$l_2$ & 2 & 1 & 1 \\
\hline
$l_{11}$ & 1 & 0 & $- 1$ \\
\hline
$l_{21}$ & 0 & 1 & 3  \\
\hline
$l_{12}$ & 3 & 2 & 3 \\
\hline
$l_{22}$ & $- 2$ & $- 1$ & $- 1$  \\
\hline
\end{tabular}
\end{center}

\vskip.03in

\underline{$n^o\ 22$.}

\vskip.03in

\hskip.1in NO Mistake.

\newpage

\underline{$n^o\ 23$.}

\vskip.03in

\hskip.1in NO Mistake.

\vskip.03in

\underline{$n^o\ 24$.}

\vskip.03in

\hskip.1in NO Mistake.

\vskip.03in

\underline{$n^o\ 25$.}

\vskip.03in

\hskip.1in NO Mistake.

\vskip.03in

\underline{$n^o\ 26$.}

\vskip.03in

About (1).

Mistake: $l_2$: ... $H \cap V_5$ is a Del Pezzo surface of degree 5 containing $C$.

Correction: $l_2$: ... $H \cap V_5$ is a Del Pezzo surface of degree 5 containing $L$.

\vskip.03in

About (2).

\vskip.03in

\hskip.1in Mistake:

\begin{center}
\begin{tabular}{|c|c|c|c|}
\hline
  &  $E$ & $Blp^*\mathcal{O}_{V_5}(1)$ & $- K_T$ \\
\hline
$l_1$ & $- 1$ & 0 & 1 \\
\hline
$l_2$ & 1 & 1 & 1 \\
\hline
$l_{11}$ & 1 & 0 & $- 1$ \\
\hline
$l_{21}$ & 0 & 1 & 2  \\
\hline
$l_{12}$ & 1 & 2 & 3 \\
\hline
$l_{21}$ & $- 1$ & $- 1$ & $- 1$  \\
\hline
\end{tabular}
\end{center}

\vskip.03in

Note: The subscript for the last line $l_{21}$ should have been $l_{22}$.

\vskip.03in

Correction:

\vskip.03in

\begin{center}
\begin{tabular}{|c|c|c|c|}
\hline
  &  $E$ & $Blp^*\mathcal{O}_{V_5}(1)$ & $- K_T$ \\
\hline
$l_1$ & $- 1$ & 0 & 1 \\
\hline
$l_2$ & 1 & 1 & 1 \\
\hline
$l_{11}$ & 1 & 0 & $- 1$ \\
\hline
$l_{21}$ & 0 & 1 & 2  \\
\hline
$l_{12}$ & 1 & 2 & 3 \\
\hline
$l_{22}$ & $- 1$ & $- 1$ & $- 1$  \\
\hline
\end{tabular}
\end{center}

\vskip.03in

\underline{$n^o\ 27$.}

\vskip.03in

\hskip.1in NO Mistake.

\vskip.03in

\underline{$n^o\ 28$.}

\vskip.03in

About (2).

\vskip.03in

\hskip.1in NO Mistake in the table of the intersection pairings.

\vskip.03in

Note: How to compute the numbers in the last two rows in the table, associated with the 4-fold flop of type ($E_5$).

\vskip.1in

In the following, we explain how to compute the numbers in the last two rows in the table, associated with the 4-fold flop of type ($E_5$), using its explicit construction given in Page 33 of \cite{WGBTMM95}.  (See the attached picture.)

\includegraphics[width=4in]{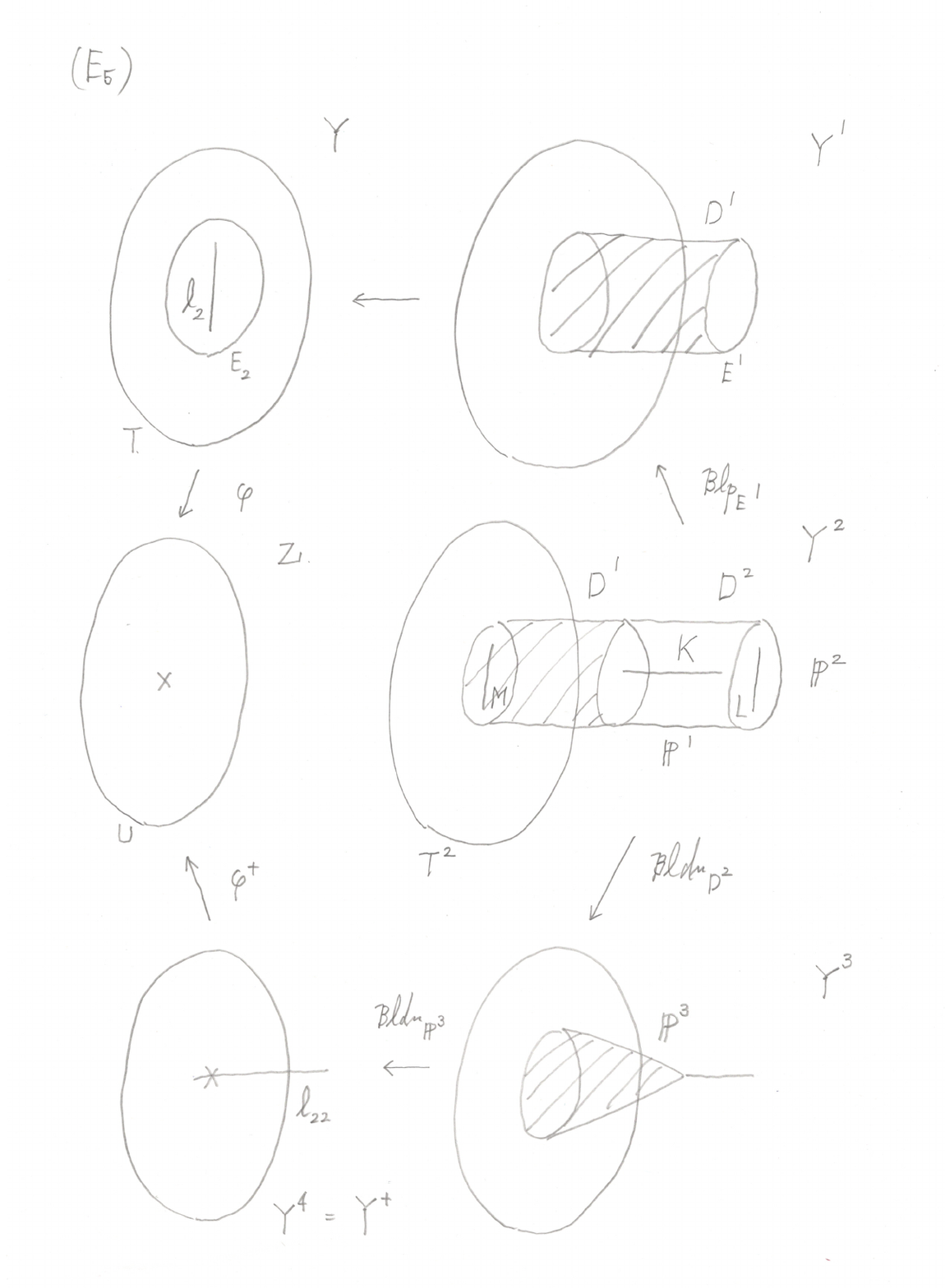}

Let $\varphi:Y \rightarrow Z$ be the contraction of flopping type on the 4-fold $Y$.  Then $\varphi|_T$ is the divisorial contraction (of the extremal ray $l_2$) $\phi:T \rightarrow U$, which contracts the exceptional divisor $E_2 \cong E_{\varphi} \cong \mathbb{P}^2$ (We note that the subscript in $E_2$ was inserted to distinguish it from the exceptional divisor $E$ of the contraction of the extremal ray $l_1$.) wih normal bundle $\mathcal{N}_{E_2/T} \cong \mathcal{O}_{\mathbb{P}^2}(-2)$.  To obtain the flop of the small contraction $\varphi:Y \rightarrow Z$, we first blow up $E_2$ in $Y$.  Since $E_2$ is the complete intersection of $T$ with $V$, the normal bundle is the direct sum $\mathcal{N}_{E_2/Y} \cong \mathcal{O}_{\mathbb{P}^2}(-2) \oplus \mathcal{O}_{\mathbb{P}^2}(-1)$  (Note that $V = p^{-1}(E_2)$ is the divisor in $Y$, where $p:Y \rightarrow T$ is the natural map from $Y$ regarded as the total space of the canonical bundle $K_T$ to $T$.  This divisor $V$ is denoted by $T_1$ in Page 33 of of \cite{WGBTMM95}.)  Thus the blow up $Blp_{E_2}:Y^1 \rightarrow Y$ has the exceptional divisor $D^1 \cong \mathbb{P}(\mathcal{O}_{\mathbb{P}^2}(-2) \oplus \mathcal{O}_{\mathbb{P}^2}(-1))$.  The natural map $D^1 \rightarrow \mathbb{P}^2$ has a section $\mathbb{P}^2 \cong E^1 \subset D^1$, which is the intersection of $D^1$ with the strict transform of $V$.  Thus $E^1$ has the normal bundle  $\mathcal{N}_{E^1/Y^1} \cong \mathcal{O}_{\mathbb{P}^2}(-1) \oplus \mathcal{O}_{\mathbb{P}^2}(-1)$.  Secondly, we blow up this section $E^1$ to have $Blp_{E^1}:Y^2 \rightarrow Y^1$ with the exceptional divisor $D^2 \cong \mathbb{P}^2 \times \mathbb{P}^1$ with normal bundle $\mathcal{N}_{D^2/Y^2} \cong \mathcal{O}(-1,-1)$.  Now contract $D^2$ into the other direction, namely into the direction of the second projection to get $Bldn_{D^2}: Y^2 \rightarrow Y^3$.  The image of the strict transform of $D^1$ under $Bldn_{D^2}$ becomes isomorphic to $\mathbb{P}^3$ with normal bundle $\mathcal{O}_{\mathbb{P}^3}(-2)$, only to be contracted onto a singular point $Bldn_{\mathbb{P}^3}: Y^3 \rightarrow Y^4$.  Finally we can contract the strict transform of the image curve $l_{22}$ of $D^2$ to get the desired flop $\varphi^+:Y^4 = Y^+ \rightarrow Z$.

We need to give some explanation of the last bottom two rows of the table.  The symbol $l_{12}$ denotes the strict transform of the extremal rational curve denoted by $l_1$, while $l_{22}$ denotes the image curve of $D^2$, which is the exceptional locus of $\varphi^+$.  The symbol $E$ and $Blp^*\mathcal{O}_{\mathbb{P}^3}(1)$ in the top row actually refer to $p^{-1}(E)$ and $p^{-1}(Blp^*\mathcal{O}_{\mathbb{P}^3}(1))$, where $p:Y \rightarrow T$ is the natural map from $Y$ regarded as the total space of the canonical bundle $K_T$ to $T$ as above.  Therefore, in the last bottom two rows of the table, we are supposed to compute the intersection numbers of $l_{12}$ and $l_{22}$ with $\{p^{-1}(E)\}^+$ and $\{p^{-1}(Blp^*\mathcal{O}_{\mathbb{P}^3}(1))\}^+$, the strict tarnsforms of the divisors $p^{-1}(E)$ and $p^{-1}(Blp^*\mathcal{O}_{\mathbb{P}^3}(1))$ on $Y^+$, which we denote simply by $E^+$ and $Blp^*\mathcal{O}_{\mathbb{P}^3}(1)^+$ by abuse of notation.

In the above sense, we would like to compute the intersection numbers

$$\begin{array}{ll}
(l_{12},E^+) & (l_{12},Blp^*\mathcal{O}_{\mathbb{P}^3}(1)^+) \\
(l_{22},E^+) & (l_{22},Blp^*\mathcal{O}_{\mathbb{P}^3}(1)^+). \\
\end{array}$$

The intersection number with $(- K_T)^+$ can be easily computed from the above using the formula
$$K_T = Blp^*K_{\mathbb{P}^3} + E = Blp^*\mathcal{O}_{\mathbb{P}^3}(-4) + E,$$
which implies
$$(- K_T)^+ = 4 Blp^*\mathcal{O}_{\mathbb{P}^3}(1)^+ - E^+.$$

Our strategy is to compute these intersection numbers, using the intersection pairings between some appropriate curves and divisors on $Y^2$.

Let $\psi: Y^2 \rightarrow Y$ and $\psi^+:Y^2 \rightarrow Y^+$ denote the corresponding morphisms.  We denote the strict transform of $D^1$ on $Y^2$ also by $D^1$ by abuse of notation.  Let $T^2$ be the strict transform of $T$ (as the divisor on $Y$) on $Y^2$.  We set

$$\left\{\begin{array}{rcl}
M &:& \text{a line in }D^1 \cap T^2 \cong \mathbb{P}^2 \\
L &:& \text{a line in a section of }D^2 \cong \mathbb{P}^2 \times \mathbb{P}^1 \rightarrow \mathbb{P}^2 \\
K &:& \text{a fiber of }D^2 \cong \mathbb{P}^2 \times \mathbb{P}^1 \rightarrow \mathbb{P}^2 \\
l_{12} &:& \text{the strict transform of }l_1 \text{ on }Y^2 \text{ by abuse of notation} \\
\end{array}\right.$$

We have
$$\left\{\begin{array}{lcl}
(\psi^+)^*E^+ &=& \psi^*E + aD^1 + bD^2 \\
(\psi^+)^*Blp^*\mathcal{O}_{\mathbb{P}^3}(1)^+ &=& \psi^*Blp^*\mathcal{O}_{\mathbb{P}^3}(1) + cD^1 + dD^2 \\
&& \text{ for some numbers }a,b,c,d. \\
\end{array}\right.$$

In order to determine the numbers $a,b$, we use the conditions

$$\left\{\begin{array}{lclcl}
(M,(\psi^+)^*E^+) &=& (M,\psi^*E + aD^1 + bD^2) &=& 0,\\
(L,(\psi^+)^*E^+ &=& (L,\psi^*E + aD^1 + bD^2) &=& 0.\\
\end{array}\right.$$

Since
$$\left\{\begin{array}{ccc}
(M,\psi^*E) = 3, & (M,D^1) = -2, & (M,D^2) = 0,\\
(L,\psi^*E) = 3, & (L,D^1) = 0, & (L,D^2) = -1,\\
\end{array}\right.$$

we conclude
$$a = 3/2, b = 3.$$

In order to determine the numbers $c,d$, we use the conditions

$$\left\{\begin{array}{rclcl}
(M,(\psi^+)^*Blp^*\mathcal{O}_{\mathbb{P}^3}(1)^+) &=& (M,\psi^*Blp^*\mathcal{O}_{\mathbb{P}^3}(1) + cD^1 + dD^2) &=& 0,\\
(L,(\psi^+)^*Blp^*\mathcal{O}_{\mathbb{P}^3}(1)^+) &=& (L,\psi^*Blp^*\mathcal{O}_{\mathbb{P}^3}(1) + cD^1 + dD^2) &=& 0.\\
\end{array}\right.$$

Since
$$\left\{\begin{array}{ccc}
(M,\psi^*Blp^*\mathcal{O}_{\mathbb{P}^3}(1)) = 1, & (M,D^1) = -2, & (M,D^2) = 0,\\
(L,\psi^*Blp^*\mathcal{O}_{\mathbb{P}^3}(1)) = 1, & (L,D^1) = 0, & (L,D^2) = -1,\\
\end{array}\right.$$

we conclude

$$c = 1/2, d = 1.$$

Finally we compute

$$\left\{\begin{array}{rcl}
(l_{12},E^+) &=& (l_{12},\psi^*E + 3/2D^1 + 3D^2) \\
&=& (-1) + 3/2 \cdot 1 + 3 \cdot 0 = 1/2 \\
(l_{12},Blp^*\mathcal{O}_{\mathbb{P}^3}(1)^+) &=& (l_{12},\psi^*Blp^*\mathcal{O}_{\mathbb{P}^3}(1) + 1/2D^1 + 1D^2) \\
&=& 0 + 1/2 \cdot 1 + 1 \cdot 0 = 1/2 \\
\end{array}\right.$$ 

since
$$(l_{12},\psi^*E) = -1, (l_{12},\psi^*Blp^*\mathcal{O}_{\mathbb{P}^3}(1)) = 0, (l_{12},D^1) = 1, (l_{12},D^2) = 0.$$

Similarly we compute

$$\left\{\begin{array}{rcl}
(l_{22},E^+) &=& (K,\psi^*E + 3/2D^1 + 3D^2) \\
&=& 0 + 3/2 \cdot 1 + 3 \cdot (-1) = - 3/2 \\
(l_{22},Blp^*\mathcal{O}_{\mathbb{P}^3}(1)^+) &=& (K,\psi^*Blp^*\mathcal{O}_{\mathbb{P}^3}(1) + 1/2D^1 + 1D^2) \\
&=& 0 + 1/2 \cdot 1 + 1 \cdot (-1) = - 1/2 \\
\end{array}\right.$$

since
$$(K,\psi^*E) = 0, (K,\psi^*Blp^*\mathcal{O}_{\mathbb{P}^3}(1)) = 0, (K,D^1) = 1, (K,D^2) = -1.$$

Accordingly, we have
$$\left\{\begin{array}{l}
(l_{12},(- K_T)^+) = (l_{12},-4 Blp^*\mathcal{O}_{\mathbb{P}^3}(1)^+ + E^+) = - \{(-4)(1/2) + 1/2\} = 3/2 \\
(l_{22},(- K_T)^+) = (l_{22},-4 Blp^*\mathcal{O}_{\mathbb{P}^3}(1)^+ + E^+) = - \{(-4)(-1/2) + (- 3/2)\} = - 1/2 \\
\end{array}\right.$$

This completes the explanation of how to compute the numbers in the last two rows in the table.

\vskip.03in

About (4).

\hskip.1in Mistake: ($E_1$)

Note:  It is missing ($E_5$).

\hskip.1in Correction: ($E_1$) ($E_5$)

\vskip.03in

\underline{$n^o\ 29$.}

\vskip.03in

\hskip.1in NO Mistake.

\vskip.03in

\underline{$n^o\ 30$.}

\vskip.03in

About (2).

\vskip.03in

\hskip.1in NO Mistake in the table of the intersection pairings.

Note: In the following, we explain how to compute the numbers in the last two rows in the table, associated with the 4-fold flop of type ($E_2$), using its explicit construction given in Page 33 of \cite{WGBTMM95}.  (See the attached picture.)

\includegraphics[width=4in]{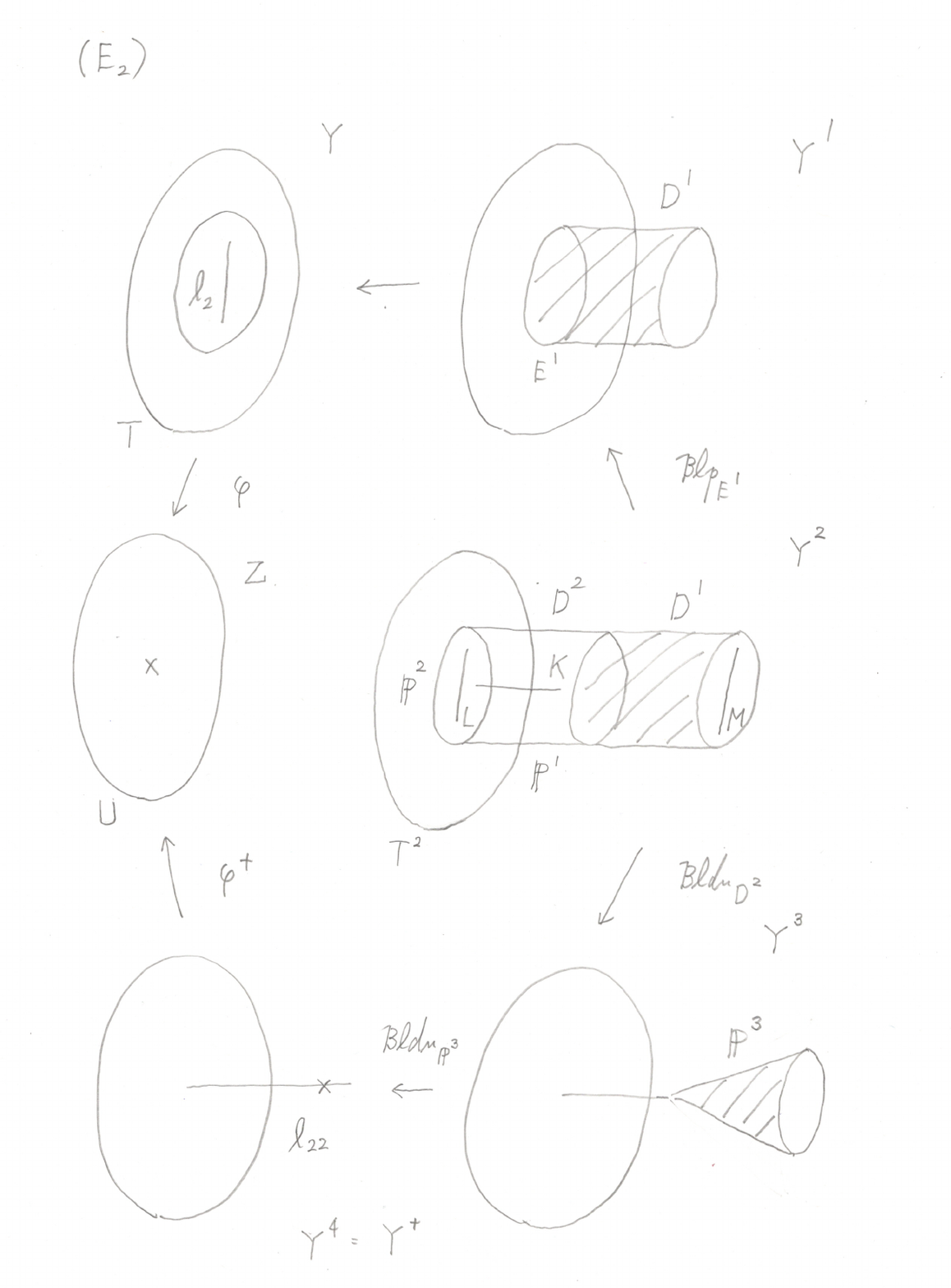}

Let $\varphi:Y \rightarrow Z$ be the contraction of flopping type on the 4-fold $Y$.  Then $\varphi|_T$ is the divisorial contraction (of the extremal ray $l_2$) $\phi:T \rightarrow U$, which contracts the exceptional divisor $E_2 \cong E_{\varphi} \cong \mathbb{P}^2$ (We note that the subscript in $E_2$ was inserted to distinguish it from the exceptional divisor $E$ of the contraction of the extremal ray $l_1$.) wih normal bundle $\mathcal{N}_{E_2/T} \cong \mathcal{O}_{\mathbb{P}^2}(-1)$.  To obtain the flop of the small contraction $\varphi:Y \rightarrow Z$, we first blow up $E_2$ in $Y$.  Since $E_2$ is the complete intersection of $T$ with $V$, the normal bundle is the direct sum $\mathcal{N}_{E_2/Y} \cong \mathcal{O}_{\mathbb{P}^2}(-1) \oplus \mathcal{O}_{\mathbb{P}^2}(-2)$  (Note that $V = p^{-1}(E_2)$ is the divisor in $Y$, where $p:Y \rightarrow T$ is the natural map from $Y$ regarded as the total space of the canonical bundle $K_T$ to $T$.  This divisor $V$ is denoted by $T_1$ in Page 31 of of \cite{WGBTMM95}.)  Thus the blow up $Blp_{E_2}:Y^1 \rightarrow Y$ has the exceptional divisor $D^1 \cong \mathbb{P}(\mathcal{O}_{\mathbb{P}^2}(-1) \oplus \mathcal{O}_{\mathbb{P}^2}(-2))$.  Secondly, we blow up the intersection $E^1 \cong \mathbb{P}^2$ of $D^1$ and the strict transform of the zero section.  Since $E^1$ has the normal bundle $\mathcal{N}_{E^1/Y^1} \cong \mathcal{O}_{\mathbb{P}^2}(-1) \oplus \mathcal{O}_{\mathbb{P}^2}(-1)$, the blow up $Blp_{E^1}:Y^2 \rightarrow Y^1$ has the exceptional divisor the exceptional divisor $D^2 \cong \mathbb{P}^2 \times \mathbb{P}^1$ with normal bundle $\mathcal{N}_{D^2/Y^2} \cong \mathcal{O}(-1,-1)$.  Now contract $D^2$ into the other direction, namely into the direction of the second projection to get $Bldn_{D^2}: Y^2 \rightarrow Y^3$.  The image of the strict transform of $D^1$ under $Bldn_{D^2}$ becomes isomorphic to $\mathbb{P}^3$ with normal bundle $\mathcal{O}_{\mathbb{P}^3}(-2)$, only to be contracted onto a singular point $Bldn_{\mathbb{P}^3}: Y^3 \rightarrow Y^4$.  Finally we can contract the strict transform of the image curve $l_{22}$ of $D^2$ to get the desired flop $\varphi^+:Y^4 = Y^+ \rightarrow Z$.

We need to give some explanation of the last bottom two rows of the table.  The symbol $l_{12}$ denotes the strict transform of the extremal rational curve denoted by $l_1$, while $l_{22}$ denotes the image curve of $D^2$, which is the exceptional locus of $\varphi^+$.  The symbol $E$ and $Blp^*\mathcal{O}_{\mathbb{P}^3}(1)$ in the top row actually refer to $p^{-1}(E)$ and $p^{-1}(Blp^*\mathcal{O}_{\mathbb{P}^3}(1))$, where $p:Y \rightarrow T$ is the natural map from $Y$ regarded as the total space of the canonical bundle $K_T$ to $T$ as above.  Therefore, in the last bottom two rows of the table, we are supposed to compute the intersection numbers of $l_{12}$ and $l_{22}$ with $\{p^{-1}(E)\}^+$ and $\{p^{-1}(Blp^*\mathcal{O}_{\mathbb{P}^3}(1))\}^+$, the strict tarnsforms of the divisors $p^{-1}(E)$ and $p^{-1}(Blp^*\mathcal{O}_{\mathbb{P}^3}(1))$ on $Y^+$, which we denote simply by $E^+$ and $Blp^*\mathcal{O}_{\mathbb{P}^3}(1)^+$ by abuse of notation.

In the above sense, we would like to compute the intersection numbers

$$\begin{array}{ll}
(l_{12},E^+) & (l_{12},Blp^*\mathcal{O}_{\mathbb{P}^3}(1)^+) \\
(l_{22},E^+) & (l_{22},Blp^*\mathcal{O}_{\mathbb{P}^3}(1)^+). \\
\end{array}$$

The intersection number with $(- K_T)^+$ can be easily computed from the above using the formula
$$K_T = Blp^*K_{\mathbb{P}^3} + E = Blp^*\mathcal{O}_{\mathbb{P}^3}(-4) + E,$$
which implies
$$(- K_T)^+ = 4 Blp^*\mathcal{O}_{\mathbb{P}^3}(1)^+ - E^+.$$

Our strategy is to compute these intersection numbers, using the intersection pairings between some appropriate curves and divisors on $Y^2$.

Let $\psi: Y^2 \rightarrow Y$ and $\psi^+:Y^2 \rightarrow Y^+$ denote the corresponding morphisms.  We denote the strict transform of $D^1$ on $Y^2$ also by $D^1$ by abuse of notation.  Let $T^2$ be the strict transform of $T$ (as the divisor on $Y$) on $Y^2$.  We set

$$\left\{\begin{array}{rcl}
M &:& \text{a line in the section of }D^1  \cong \mathbb{P}(\mathcal{O}_{\mathbb{P}^2}(-1) \oplus \mathcal{O}_{\mathbb{P}^2}(-2)) \rightarrow \mathbb{P}^2,\\
&& \text{the section which is away from }D^2 \cap D^1 \\
L &:& \text{a line in }D^2 \cap T^2\\
K &:& \text{a fiber of }D^2 \cong \mathbb{P}^2 \times \mathbb{P}^1 \rightarrow \mathbb{P}^2 \\
l_{12} &:& \text{the strict transform of }l_1 \text{ on }Y^2 \text{ by abuse of notation} \\
\end{array}\right.$$

We have
$$\left\{\begin{array}{rcl}
(\psi^+)^*E^+ &=& \psi^*E + aD^1 + bD^2 \\
(\psi^+)^*Blp^*\mathcal{O}_{\mathbb{P}^3}(1)^+ &=& \psi^*Blp^*\mathcal{O}_{\mathbb{P}^3}(1) + cD^1 + dD^2 \\
&& \text{ for some numbers }a,b,c,d. \\
\end{array}\right.$$

In order to determine the numbers $a,b$, we use the conditions

$$\left\{\begin{array}{rclcl}
(M,(\psi^+)^*E^+) &=& (M,\psi^*E + aD^1 + bD^2) &=& 0,\\
(L,(\psi^+)^*E^+) &=& (L,\psi^*E + aD^1 + bD^2) &=& 0.\\
\end{array}\right.$$

Since
$$\left\{\begin{array}{ccc}
(M,\psi^*E) = 2, & (M,D^1) = -2, & (M,D^2) = 0,\\
(L,\psi^*E) = 2, & (L,D^1) = 0, & (L,D^2) = -1,\\
\end{array}\right.$$

we conclude
$$a = 1, b = 2.$$

In order to determine the numbers $c,d$, we use the conditions

$$\left\{\begin{array}{rclcl}
(M,(\psi^+)^*Blp^*\mathcal{O}_{\mathbb{P}^3}(1)^+) &=& (M,\psi^*Blp^*\mathcal{O}_{\mathbb{P}^3}(1) + cD^1 + dD^2) &=& 0,\\
(L,(\psi^+)^*Blp^*\mathcal{O}_{\mathbb{P}^3}(1)^+) &=& (L,\psi^*Blp^*\mathcal{O}_{\mathbb{P}^3}(1) + cD^1 + dD^2) &=& 0.\\
\end{array}\right.$$

Since
$$\left\{\begin{array}{ccc}
(M,\psi^*Blp^*\mathcal{O}_{\mathbb{P}^3}(1)) = 1, & (M,D^1) = -2, & (M,D^2) = 0,\\
(L,\psi^*Blp^*\mathcal{O}_{\mathbb{P}^3}(1)) = 1, & (L,D^1) = 0, & (L,D^2) = -1,\\
\end{array}\right.$$

we conclude

$$c = 1/2, d = 1.$$

Finally we compute

$$\left\{\begin{array}{rcl}
(l_{12},E^+) &=& (l_{12},\psi^*E + 1D^1 + 2D^2) \\
&=& (-1) + 1 \cdot 0 + 2 \cdot 1 = 1 \\
(l_{12},Blp^*\mathcal{O}_{\mathbb{P}^3}(1)^+) &=& (l_{12},\psi^*Blp^*\mathcal{O}_{\mathbb{P}^3}(1) + 1/2D^1 + 1D^2) \\
&=& 0 + 1/2 \cdot 0 + 1 \cdot 1 = 1 \\
\end{array}\right.$$ 

since
$$(l_{12},\psi^*E) = -1, (l_{12},\psi^*Blp^*\mathcal{O}_{\mathbb{P}^3}(1)) = 0, (l_{12},D^1) = 0, (l_{12},D^2) = 1.$$

Similarly we compute

$$\left\{\begin{array}{rcl}
(l_{22},E^+) &=& (l_{22},\psi^*E + 1D^1 + 2D^2) \\
&=& 0 + 1 \cdot 1 + 2 \cdot (-1) = - 1 \\
(l_{22},Blp^*\mathcal{O}_{\mathbb{P}^3}(1)^+) &=& (K,\psi^*Blp^*\mathcal{O}_{\mathbb{P}^3}(1) + 1/2D^1 + 1D^2) \\
&=& 0 + 1/2 \cdot 1 + 1 \cdot (-1) = - 1/2 \\
\end{array}\right.$$

since
$$(K,\psi^*E) = 0, (K,\psi^*Blp^*\mathcal{O}_{\mathbb{P}^3}(1)) = 0, (K,D^1)  = 1, (K,D^2) = -1.$$

Accordingly, we have
$$\left\{\begin{array}{l}
(l_{12},(- K_T)^+) = (l_{12},-4 Blp^*\mathcal{O}_{\mathbb{P}^3}(1)^+ + E^+) = - \{(-4)(1/2) + 1/2\} = 3 \\
(l_{22},(- K_T)^+) = (l_{22},-4 Blp^*\mathcal{O}_{\mathbb{P}^3}(1)^+ + E^+) = - \{(-4)(-1/2) + (- 3/2)\} = - 1 \\
\end{array}\right.$$

This completes the explanation of how to compute the numbers in the last two rows in the table.

\vskip.03in

\underline{$n^o\ 31$.}

\vskip.03in

\hskip.1in NO Mistake.

\vskip.03in

\underline{$n^o\ 32$.}

\vskip.03in

\hskip.1in NO Mistake.

\vskip.03in

\underline{$n^o\ 33$.}

\vskip.03in

\hskip.1in NO Mistake.

\vskip.03in

\underline{$n^o\ 34$.}

\vskip.03in

\hskip.1in NO Mistake.

\vskip.03in

\underline{$n^o\ 35$.}

\vskip.03in

About (2).

\vskip.03in

\hskip.1in Mistake:

\begin{center}
\begin{tabular}{|c|c|c|c|}
\hline
  &  $E$ & $Blp^*\mathcal{O}_{\mathbb{P}^3}(1)$ & $- K_T$ \\
\hline
$l_1$ & $- 1$ & 0 & 1 \\
\hline
$l_2$ & 1 & 1 & 2 \\
\hline
$l_{11}$ & $\frac{1}{2}$ & 0 & $- 1$ \\
\hline
$l_{21}$ & 0 & 1 & 4  \\
\hline
\end{tabular}
\end{center}

\vskip.03in

\hskip.1in Correction:

\begin{center}
\begin{tabular}{|c|c|c|c|}
\hline
  &  $E$ & $Blp^*\mathcal{O}_{\mathbb{P}^3}(1)$ & $- K_T$ \\
\hline
$l_1$ & $- 1$ & 0 & 2 \\
\hline
$l_2$ & 1 & 1 & 2 \\
\hline
$l_{11}$ & $\frac{1}{2}$ & 0 & $- 1$ \\
\hline
$l_{21}$ & 0 & 1 & 4  \\
\hline
\end{tabular}
\end{center}

\vskip.03in

About (4).

\hskip.1in Mistake: ($E_1$)

\hskip.1in Correction: ($E_2$)

\newpage

\underline{$n^o\ 36$.}

\vskip.03in

About (2).

\vskip.03in

\hskip.1in Mistake:

\begin{center}
\begin{tabular}{|c|c|c|c|}
\hline
  &  $E$ & $Blp^*\mathcal{O}_{W_4}(1)$ & $- K_T$ \\
\hline
$l_1$ & $- 2$ & 0 & 2 \\
\hline
$l_2$ & 1 & 1 & 2 \\
\hline
$l_{11}$ & 1 & 0 & $- 1$ \\
\hline
$l_{21}$ & 0 & 1 & 3  \\
\hline
\end{tabular}
\end{center}

\vskip.03in

\hskip.1in Correction:

\begin{center}
\begin{tabular}{|c|c|c|c|}
\hline
  &  $E$ & $Blp^*\mathcal{O}_{W_4}(1)$ & $- K_T$ \\
\hline
$l_1$ & $- 2$ & 0 & 1 \\
\hline
$l_2$ & 1 & 1 & 2 \\
\hline
$l_{11}$ & 1 & 0 & $- \frac{1}{2}$ \\
\hline
$l_{21}$ & 0 & 1 & $\frac{5}{2}$  \\
\hline
\end{tabular}
\end{center}

\newpage

$\boxed{B_2 = 3}$

\vskip.03in

\underline{$n^o\ 1$.}

\vskip.03in

\hskip.1in NO Mistake.

\vskip.03in

\underline{$n^o\ 2$.}

\vskip.03in

About (1).

Mistake:

$T \cap Y = D$ is a divisor such that $c|_T: D \rightarrow \mathbb{P}^1 \times \mathbb{P}^1$ gives a double cover ...

Correction:

$T \cap Y = D$ is a divisor such that $c|_D: D \rightarrow \mathbb{P}^1 \times \mathbb{P}^1$ gives a double cover ...

\vskip.03in

About (2).

\vskip.03in

\hskip.1in Mistake:

\begin{center}
\begin{tabular}{|c|c|c|c|c|}
\hline
  &  $c^*p_1^*\mathcal{O}_{\mathbb{P}^1}(1)$ & $c^*p_2^*\mathcal{O}_{\mathbb{P}^1}(1)$ & $D$ & $- K_T$ \\
\hline
$l_{33}$ & $- 1$ & 0 & 3 & $- 1$ \\
\hline
\end{tabular}
\end{center}

\vskip.03in

\hskip.1in Correction:

\begin{center}
\begin{tabular}{|c|c|c|c|c|}
\hline
  &  $c^*p_1^*\mathcal{O}_{\mathbb{P}^1}(1)$ & $c^*p_2^*\mathcal{O}_{\mathbb{P}^1}(1)$ & $D$ & $- K_T$ \\
\hline
$l_{33}$ & $- 1$ & 0 & 1 & $- 1$ \\
\hline
\end{tabular}
\end{center}

About (3).

The picture of the KKMR decomposition should be asymmetric.

\includegraphics[width=3.3in]{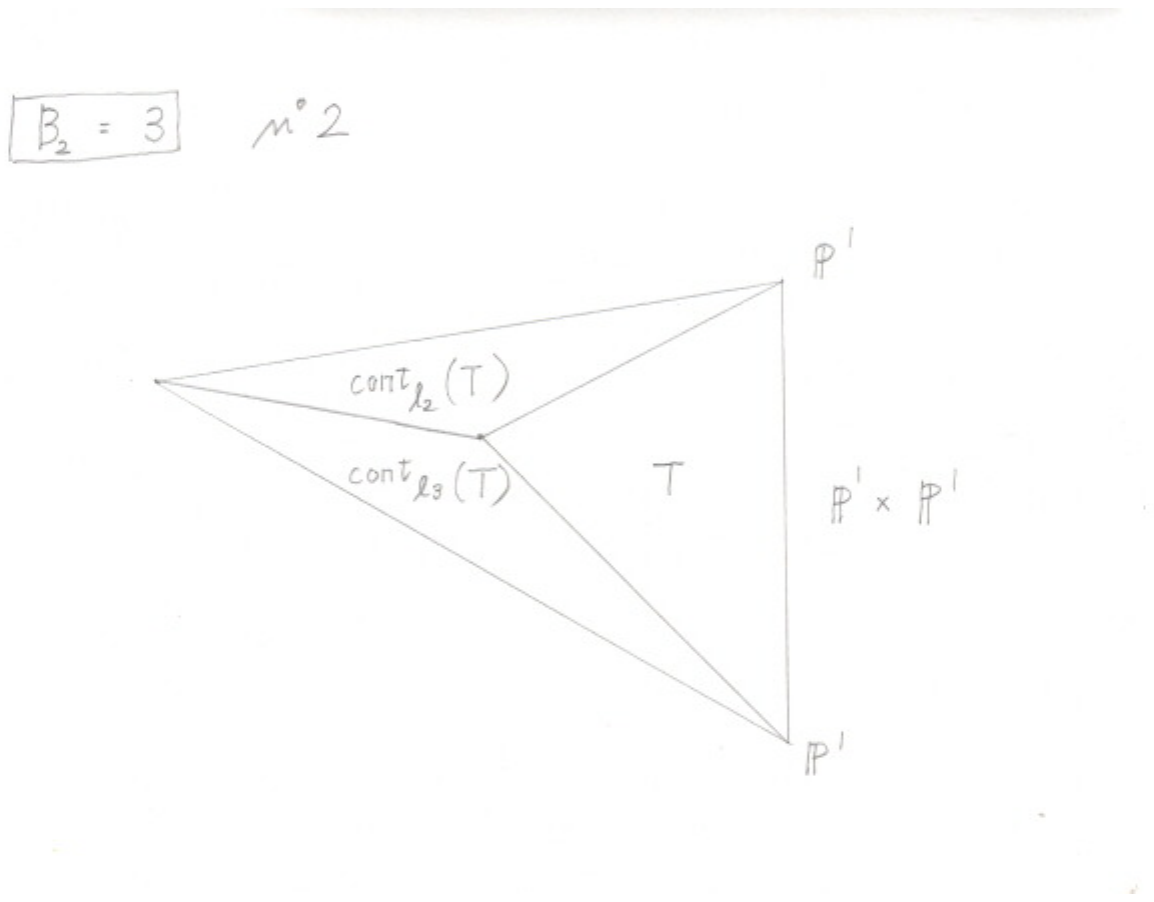}

\vskip.1in

\underline{$n^o\ 3$.}

\vskip.03in

About (1).

\hskip.1in Mistake:

$l_1$: an irreducible component of a reducible fiber ora reduced part of of a multiple fiber of the conic bundle $p_{11}|_T$

\vskip.03in

\hskip.1in Correction:

$l_1$: an irreducible component of a reducible fiber or a reduced part of of a multiple fiber of the conic bundle $p_{12}|_T$

\vskip.1in

\vskip.03in

About (2).

\vskip.03in

Mistake:

\begin{center}
\begin{tabular}{|c|c|c|c|c|}
\hline
  &  $p_1^*\mathcal{O}_{\mathbb{P}^1}(1)$ & $p_2^*\mathcal{O}_{\mathbb{P}^1}(1)$ & $p_3^*\mathcal{O}_{\mathbb{P}^2}(1)$  & $- K_T$ \\
\hline
$l_1$ & $- 0$ & 0 & 1 & 1 \\
\hline
\end{tabular}
\end{center}

Correction:

\begin{center}
\begin{tabular}{|c|c|c|c|c|}
\hline
  &  $p_1^*\mathcal{O}_{\mathbb{P}^1}(1)$ & $p_2^*\mathcal{O}_{\mathbb{P}^1}(1)$ & $p_3^*\mathcal{O}_{\mathbb{P}^2}(1)$  & $- K_T$ \\
\hline
$l_1$ & $0$ & 0 & 1 & 1 \\
\hline
\end{tabular}
\end{center}

\vskip.03in

\underline{$n^o\ 4$.}

\vskip.03in

About (2).

\vskip.03in

\hskip.1in NO Mistake in the table of the intersection pairings.

\vskip.1in

About the picture.

\hskip.1in  Mistake and Correction: In the right lower corner, $\mathbb{P}^2 \times \mathbb{P}^1$ should be replaced by $\mathbb{P}^1 \times \mathbb{P}^2$.  The arrow $\mathbb{F}^1 \leftarrow \mathbb{P}^1$ should be reversed.

\vskip.1in

\underline{$n^o\ 5$.}

\vskip.1in

\hskip.1in NO Mistake.

\vskip.1in

\underline{$n^o\ 6$.}

\vskip.1in

\hskip.1in NO Mistake.

\vskip.1in

\underline{$n^o\ 7$.}

\vskip.1in

About the description.

\hskip.1in Mistake:

$T$: the blow up of $W_6$ with center an elliptic curve which is an intersection of two members of $|- \frac{1}{2}K_W|$

\hskip.1in Correction:

$T$: the blow up of $W_6$ with center an elliptic curve which is an intersection of two members of $|- \frac{1}{2}K_{W_6}|$

\vskip.1in

About the picture.

\hskip.1in  Mistake and Correction: The labelings ``$\mathrm{cont}_{l_2}$'' and ``$\mathrm{cont}_{l_3}$'' should be switched.

\newpage

About (2).

\vskip.03in

\hskip.1in Mistake:

\begin{center}
\begin{tabular}{|c|c|c|c|c|}
\hline
  &  $E$ & $\pi_1^*\mathcal{O}_{\mathbb{P}^2}(1)$ & $\pi_2^*\mathcal{O}_{\mathbb{P}^2}(1)$ & $- K_T$ \\
\hline
$l_{32}$ & 3 & 1 & 2 & 2 \\
\hline
$l_{23}$ & 3 & 2 & 1 & 3 \\
\hline
\end{tabular}
\end{center}

\hskip.1in Correction:

\begin{center}
\begin{tabular}{|c|c|c|c|c|}
\hline
  &  $E$ & $\pi_1^*\mathcal{O}_{\mathbb{P}^2}(1)$ & $\pi_2^*\mathcal{O}_{\mathbb{P}^2}(1)$ & $- K_T$ \\
\hline
$l_{32}$ & 2 & 1 & 1 & 2 \\
\hline
$l_{23}$ & 2 & 1 & 1 & 2 \\
\hline
\end{tabular}
\end{center}

\vskip.1in

\underline{$n^o\ 8$.}

\vskip.03in

\hskip.1in NO Mistake.

\vskip.03in

\underline{$n^o\ 9$.}

\vskip.03in

About the picture.

Mistake: There is no definition of the morphism $\pi$.

Correction: The morphism $\mathbb{P}^2 \overset{\pi}\longleftarrow \mathbb{P}_{\mathbb{P}^2}(\mathcal{O} \oplus \mathcal{O}(2))$ is the natural map from the $\mathbb{P}^1$-bundle on the right-hand bottom corner of the picture to the base space in the center.

\vskip.03in

About (3).

Mistake and Correction: In the picture of the KKMR decomposition, the top chamber is labeled as ``$\text{cont}_{l_2l_3}(T)$''.  The lable should be replaced with ``$\text{cont}_{l_2l_4}(T)$''

\vskip.03in

\underline{$n^o\ 10$.}

\vskip.03in

About the picture.

Mistake and Correction: At the destination of the arrow $\overset{\pi}\longrightarrow$, the symbol $Q$ representing the quadric is missing.

\vskip.03in

\underline{$n^o\ 11$.}

\vskip.03in

About the picture.

\vskip.03in

\hskip.1in Mistake and Correction: The extremal rational curve $l_2$ should pass through the intersection of $E_2$ and $E_1$, which is the section of $\mathbb{F}^1 \rightarrow \mathbb{P}^1$.

\vskip.03in

About (2).

\hskip.1in Mistake:

\begin{center}
\begin{tabular}{|c|c|c|c|c|}
\hline
  &  $E_1$ & $E_2$ & $\pi^*\mathcal{O}_{\mathbb{P}^3}(1)$ & $- K_T$ \\
\hline
$l_{23}$ & 1 & $-1$ & 0 & 1 \\
\hline
\end{tabular}
\end{center}

\hskip.1in Correction:

\begin{center}
\begin{tabular}{|c|c|c|c|c|}
\hline
  &  $E_1$ & $E_2$ & $\pi^*\mathcal{O}_{\mathbb{P}^3}(1)$ & $- K_T$ \\
\hline
$l_{23}$ & 3 & 1 & 2 & 3 \\
\hline
\end{tabular}
\end{center}

\vskip.03in

About (3) the KKMR decomposition.

Note: The flopping of the extremal ray $l_1$, followed by the flopping of the extremal ray $l_{21}$ on the 4-fold $Y$, gives $\mathbb{P}^3$ when restricted to the strict transfiorm of the zero section.  Similarly, the flopping of the extremal ray $l_2$, followed by the flopping of the extremal ray $l_{12}$ on the 4-fold $Y$, gives the same $\mathbb{P}^3$ when restricted to the strict transform of the zero section.  However, The flopping of the extremal ray $l_1$, followed by the flopping of the extremal ray  $l_{21}$, reaches the chamber in the upper right hand corner, while the flopping of the extremal ray $l_2$, followed by the flopping of the extremal ray  $l_{12}$, reaches the chamber in the lower right hand corner.  These two chambers, though labeled with the same $\mathbb{P}^3$, correspond to two different minimal models in the KKMR decomposition, and connected  by a flop of one of the Types labeled as (Others).

\vskip.03in

About (4).

\hskip.1in Mistake: ($E_1$)($E_2$)($E_5$)

\hskip.1in Correction: ($E_1$)($E_2$)(Others)

Note: There is no flop of type ($E_5$).

\vskip.03in

\underline{$n^o\ 12$.}

\vskip.03in

About (2).

\hskip.1in Mistake:

\begin{center}
\begin{tabular}{|c|c|c|c|c|}
\hline
  &  $E_1$ & $E_2$ & $\pi^*\mathcal{O}_{\mathbb{P}^3}(1)$ & $- K_T$ \\
\hline
$l_{22}$ & 0 & 1 & 0 & 1  \\
\hline
$l_{23}$ & 2 & 3 & 2 & 3 \\
\hline
\end{tabular}
\end{center}

\vskip.03in

\hskip.03in Correction:

\begin{center}
\begin{tabular}{|c|c|c|c|c|}
\hline
  &  $E_1$ & $E_2$ & $\pi^*\mathcal{O}_{\mathbb{P}^3}(1)$ & $- K_T$ \\
\hline
$l_{22}$ & 0 & 1 & 0 & $- 1$  \\
\hline
$l_{33}$ & 1 & 1 & 1 & 2 \\
\hline
\end{tabular}
\end{center}

\vskip.03in

\underline{$n^o\ 13$.}

\vskip.03in

About (2).

\vskip.03in

\hskip.1in Mistake:

\begin{center}
\begin{tabular}{|c|c|c|c|c|}
\hline
  &  $E_1$ & $\pi_1^*\mathcal{O}_{\mathbb{P}^2}(1)$ & $\pi_2^*\mathcal{O}_{\mathbb{P}^2}(1)$ & $- K_T$ \\
\hline
$l_1$ & $- 1$ & 0 & 0 & 1 \\
\hline
$l_2$ & 1 & 1 & 0 & 1 \\
\hline
$l_3$ & 1 & 0 & 1 & 1 \\
\hline
$l_{11}$ & 1 & 0 & 0 & $- 1$ \\
\hline
$l_{21}$ & 0 & 1 & 0 & 2  \\
\hline
$l_{31}$ & 0 & 0 & 1 & 2 \\
\hline
$l_{12}$ & 0 & 1 & 0 & 2 \\
\hline
$l_{22}$ & $- 1$ & $- 1$ & 0 & $- 1$ \\
\hline
$l_{32}$ & 2 & 1 & 1 & 2 \\
\hline
$l_{13}$ & 0 & 0 & 1 & 2 \\
\hline
$l_{23}$ & 2 & 1 & 1 & 2 \\
\hline
$l_{33}$ & $- 1$ & 0 & $- 1$ & $- 1$ \\
\hline
\end{tabular}
\end{center}

\newpage

\hskip.1in Correction:

\begin{center}
\begin{tabular}{|c|c|c|c|c|}
\hline
  &  $E_1$ & $\pi_1^*\mathcal{O}_{\mathbb{P}^2}(1)$ & $\pi_2^*\mathcal{O}_{\mathbb{P}^2}(1)$ & $- K_T$ \\
\hline
$l_1$ & $- 1$ & 0 & 0 & 1 \\
\hline
$l_2$ & 1 & 0 & 1 & 1 \\
\hline
$l_3$ & 1 & 1 & 0 & 1 \\
\hline
$l_{11}$ & 1 & 0 & 0 & $- 1$ \\
\hline
$l_{21}$ & 0 & 0 & 1 & 2  \\
\hline
$l_{31}$ & 0 & 1 & 0 & 2 \\
\hline
$l_{12}$ & 0 & 0 & 1 & 2 \\
\hline
$l_{22}$ & $- 1$ & 0 & $- 1$ & $- 1$ \\
\hline
$l_{32}$ & 1 & 1 & 0 & 1 \\
\hline
$l_{13}$ & 0 & 1 & 0 & 2 \\
\hline
$l_{23}$ & 1 & 0 & 1 & 1 \\
\hline
$l_{33}$ & $- 1$ & $- 1$ & 0 & $- 1$ \\
\hline
\end{tabular}
\end{center}

\vskip.03in

\underline{$n^o\ 14$.}

\vskip.03in

About (2).

\hskip.1in Mistake:

\begin{center}
\begin{tabular}{|c|c|c|c|c|}
\hline
  &  $E_1$ & $E_2$ & $\pi^*\mathcal{O}_{\mathbb{P}^3}(1)$ & $- K_T$ \\
\hline
$l_{14}$ & 1 & 0 & 1 & 3 \\
\hline
$l_{24}$ & 0 & $- 1$ & 0 & 4 \\
\hline
$l_{34}$ & 2 & 1 & 2 & 2 \\
\hline

\end{tabular}
\end{center}

\vskip.03in

\hskip.1in Correction:

\begin{center}
\begin{tabular}{|c|c|c|c|c|}
\hline
  &  $E_1$ & $E_2$ & $\pi^*\mathcal{O}_{\mathbb{P}^3}(1)$ & $- K_T$ \\
\hline
$l_{14}$ & $\frac{1}{2}$ & 0 & $\frac{1}{2}$ & $\frac{3}{2}$ \\
\hline
$l_{24}$ & 0 & $- 1$ & 0 & 2 \\
\hline
$l_{34}$ & 1 & 1 & 1 & 1 \\
\hline
\end{tabular}
\end{center}

\vskip.03in

NO Mistake in the other parts of the table of the intersection pairings.

\vskip.03in

About (3).

Note: $Blp_{cubic}\mathbb{P}^3 = \{n^o\ 28 \text{ with }B_2 = 2\}$.

\vskip.03in

\underline{$n^o\ 15$.}

\vskip.03in

About the picture.

\vskip.03in

\hskip.1in Mistake and Correction: $n^o\ 33$ should be $n^o\ 31$.  In the upper left corner $\mathbb{P}^2$ should be replaced by $\mathbb{P}^1$, while in the lower left corner $\mathbb{P}^1$ should be replaced by $\mathbb{P}^2$.

\newpage

About (2).

\hskip.1in Mistake:

\begin{center}
\begin{tabular}{|c|c|c|c|c|}
\hline
  &  $E_1$ & $E_2$ & $\pi^*\mathcal{O}_{Q}(1)$ & $- K_T$ \\
\hline
$l_{23}$ & 0 & $- 1$ & 0 & 1 \\
\hline
\end{tabular}
\end{center}

\vskip.1in

\hskip.1in Correction:

\begin{center}
\begin{tabular}{|c|c|c|c|c|}
\hline
  &  $E_1$ & $E_2$ & $\pi^*\mathcal{O}_{Q}(1)$ & $- K_T$ \\
\hline
$l_{23}$ & 1 & 0 & 1 & 2 \\
\hline
\end{tabular}
\end{center}

\vskip.1in

NO Mistake in the other parts of the table of the intersection pairings.

\vskip.1in

About (3) the KKMR decomposition.

\hskip.1in Mistake amd Correction: $n^o\ 33$ should be $n^o\ 31$.

\vskip.03in

\underline{$n^o\ 16$.}

\vskip.03in

About the picture.

Mistake and Correction: In the picture, it looks like the extremal rational curve $l_2$ does not intersect $E_3$, which is the strict transform of the singular quadric containing the twisted cubic and whose vertex is the point of blow up $V_7 \rightarrow \mathbb{P}^3$.  But actually $l_2$ intersects $E_3$ transversally once, and hence $l_2 \cdot E_3 = 1$.

\vskip.03in

About (2).

The explanation of the symbols $E_1$ and $E_2$ on the top row of the table is missing:

\begin{itemize}

\item[] $E_1$: the exceptional divisor of the blow up of the strict transform of \text{the twisted cubic}

\item[] $E_2$: the strict transform of the exceptional divisor of the blow up $V_7 \rightarrow \mathbb{P}^3$

\end{itemize}

\vskip.03in

NO Mistake in the table of the intersection pairings.

\vskip.03in

About (3) the KKMR decomposition.

Note: The flopping of the extremal ray $l_1$, followed by the flopping of the extremal ray $l_{21}$ on the 4-fold $Y$, gives $\mathbb{P}^3$ when restricted to the strict transfiorm of the zero section.  Similarly, the flopping of the extremal ray $l_2$, followed by the flopping of the extremal ray $l_{12}$ on the 4-fold $Y$, gives the same $\mathbb{P}^3$ when restricted to the strict transform of the zero section.  However, The flopping of the extremal ray $l_1$, followed by the flopping of the extremal ray  $l_{21}$, reaches the chamber in the lower right hand corner, while the flopping of the extremal ray $l_2$, followed by the flopping of the extremal ray  $l_{12}$, reaches the chamber in the upper right hand corner.  These two chambers, though labeled with the same $\mathbb{P}^3$, correspond to two different minimal models in the KKMR decomposition, and connected by a flop of one of the Types labeled as (Others)..

\vskip.03in

About (4).

\hskip.1in Mistake: ($E_1$)($E_2$)($E_5$)

\hskip.1in Correction: ($E_1$)($E_2$)(Others)

Note: There is no flop of type ($E_5$).

\newpage

\underline{$n^o\ 17$.}

\vskip.03in

About (1) the description of the extremal rays

\vskip.03in

\hskip.1in Mistake:

\begin{itemize}

\item[] $l_1$: the ruling of the exceptional divisor of the projection $p_{13}$ (the center of the blow up is a curve on $\mathbb{P}^1 \times \mathbb{P}^2$ which is a complete intersection two divisors of type $(1,0)$ and $(0,1)$) ... ($E_1$)

\item[] $l_2$: the ruling of the exceptional divisor of the projection $p_{23}$ (the center of the blow up is a curve on $\mathbb{P}^1 \times \mathbb{P}^2$ which is a complete intersection two divisors of type $(1,0)$ and $(0,1)$) ... ($E_1$)

\end{itemize}

\vskip.03in

\hskip.1in Correction:

\begin{itemize}

\item[] $l_1$: the ruling of the exceptional divisor of the projection $p_{13}$ (the center of the blow up is a curve on $\mathbb{P}^1 \times \mathbb{P}^2$ which is a complete intersection two divisors of bidegree $(1,1)$) ... ($E_1$)

\item[] $l_2$: the ruling of the exceptional divisor of the projection $p_{23}$ (the center of the blow up is a curve on $\mathbb{P}^1 \times \mathbb{P}^2$ which is a complete intersection two divisors of bidegree $(1,1)$) ... ($E_1$)

\end{itemize}

\vskip.03in

About (2).

\hskip.1in Mistake:

\begin{center}
\begin{tabular}{|c|c|c|c|c|}
\hline
  &  $p_1^*\mathcal{O}_{\mathbb{P}^1}(1)$ & $p_2^*\mathcal{O}_{\mathbb{P}^1}(1)$ & $p_3^*\mathcal{O}_{\mathbb{P}^2}(1)$ & $- K_T$ \\
\hline
$l_{11}$ & 0 & $- 1$ & 0 & 1 \\
\hline
$l_{21}$ & 1 & 1 & 0 & $- 2$ \\
\hline
$l_{31}$ & 0 & 0 & 1 & $- 2$ \\
\hline
$l_{12}$ & 1 & 1 & 0 & $- 2$ \\
\hline
$l_{22}$ & $- 1$ & 0 & 0 & 1 \\
\hline
$l_{32}$ & 0 & 0 & 1 & $- 2$ \\
\hline
\end{tabular}
\end{center}

\vskip.03in

\hskip.1in Correction:

\begin{center}
\begin{tabular}{|c|c|c|c|c|}
\hline
  &  $p_1^*\mathcal{O}_{\mathbb{P}^1}(1)$ & $p_2^*\mathcal{O}_{\mathbb{P}^1}(1)$ & $p_3^*\mathcal{O}_{\mathbb{P}^2}(1)$ & $- K_T$ \\
\hline
$l_{11}$ & 0 & $- 1$ & 0 & $- 1$ \\
\hline
$l_{21}$ & 1 & 1 & 0 & 2 \\
\hline
$l_{31}$ & 0 & 0 & 1 & 3 \\
\hline
$l_{12}$ & 1 & 1 & 0 & 2 \\
\hline
$l_{22}$ & $- 1$ & 0 & 0 & $- 1$ \\
\hline
$l_{32}$ & 0 & 0 & 1 & 3 \\
\hline
\end{tabular}
\end{center}

\vskip.03in

\underline{$n^o\ 18$.}

\vskip.03in

About (2).

The explanation of the symbols $E_1$ and $E_2$ on the top row of the table is missing:

\begin{itemize}

\item[] $E_1$: the exceptional divisor of the blow up of the line

\item[] $E_2$: the exceptional divisor of the blow up of the conic

\end{itemize}

\newpage

About (3) the KKMR decomposition.

Note: The flopping of the extremal ray $l_1$, followed by the flopping of the extremal ray $l_{31}$ on the 4-fold $Y$, gives the 3-fold which is obtained from $n^o\ 30$ by contracting the strict transform of the hyperplane $\supset\ \text{conic}$, when restricted to the strict transform of the zero section on the 4-fold.  Similarly, the flopping of the extremal ray $l_3$, followed by the flopping of the extremal ray $l_{13}$ on the 4-fold $Y$, gives the same 3-fold as above, when restricted to the strict transform of the zero section on the 4-fold.  However, The flopping of the extremal ray $l_1$, followed by the flopping of the extremal ray  $l_{31}$, reaches the chamber in the upper right hand corner, while the flopping of the extremal ray $l_3$, followed by the flopping of the extremal ray  $l_{13}$, reaches the chamber in the upper left hand corner.  These two chambers correspond to two different minimal models in the KKMR decomposition, and connected  by a flop of one of the Types labeled as (Others).

\vskip.03in

About (4).

\hskip.1in Mistake: ($E_1$)($E_2$)($E_5$)

\hskip.1in Correction: ($E_1$)($E_2$)(Others)

Note: There is no flop of type ($E_5$).

\vskip.03in

\underline{$n^o\ 19$.}

\vskip.03in

About (2).

The explanation of the symbols $E_1$ and $E_2$ on the top row of the table is missing:

\begin{itemize}

\item[] $E_1$: the exceptional divisor of the blow up of the point $p$

\item[] $E_2$: the exceptional divisor of the blow up of the point $q$

\end{itemize}

\hskip.1in NO Mistake in the table of the intersection pairings.

\vskip.03in

\underline{$n^o\ 20$.}

\vskip.03in

About (2).

The explanation of the symbols $E_1$ and $E_2$ on the top row of the table is missing:

\begin{itemize}

\item[] $E_1$: the exceptional divisor of the blow up of the line $L_1$

\item[] $E_2$: the exceptional divisor of the blow up of the line $L_2$

\end{itemize}

\hskip.1in NO Mistake in the table of the intersection pairings.

\vskip.03in

\underline{$n^o\ 21$.}

\vskip.03in

About (1).

The explanation of the symbol $C$ is missing, and so is the explanation of the hyperplane $H$.

\begin{itemize}

\item[] $C$: the curve of bidegree $(2,1)$ to blow up

\item[] $H$: the hyperplane in $\mathbb{P}^2$ such that $\mathbb{P}^1 \times H \subset \mathbb{P}^1 \times \mathbb{P}^2$ contains $C$.

\end{itemize}

\vskip.03in

About (2).

The explanation of the symbol $E_1$ on the top row of the table is missing:

\begin{itemize}

\item[] $E_1$: the exceptional divisor of the blow up of the curve of bidegree $(2,1)$

\end{itemize}

\hskip.1in NO Mistake in the table of the intersection pairings.

\newpage

\underline{$n^o\ 22$.}

\vskip.03in

About (2).

The explanation of the symbol $E_1$ on the top row of the table is missing:

\begin{itemize}

\item[] $E_1$: the exceptional divisor of the blow up of the conic

\end{itemize}

\hskip.1in NO Mistake in the table of the intersection pairings.

\vskip.03in

\underline{$n^o\ 23$.}

\vskip.03in

About (2).

The explanation of the symbols $E_1$ and $E_2$ on the top row of the table is missing:

\begin{itemize}

\item[] $E_1$: the exceptional divisor of the blow up of the strict transform of the conic

\item[] $E_2$: the strict transform of the exceptional divisor of the blow up  $V_7 \rightarrow \mathbb{P}^3$

\end{itemize}

\hskip.1in NO Mistake in the table of the intersection pairings.

\vskip.03in

About (3) the KKMR decomposition.

Note: The flopping of the extremal ray $l_1$, followed by the flopping of the extremal ray $l_{21}$ on the 4-fold $Y$, gives $\mathbb{P}^3$ when restricted to the strict transform of the zero section.  Similarly, the flopping of the extremal ray $l_2$, followed by the flopping of the extremal ray $l_{12}$ on the 4-fold $Y$, gives the same $\mathbb{P}^3$ when restricted to the strict transform of the zero section.  However, The flopping of the extremal ray $l_1$, followed by the flopping of the extremal ray  $l_{21}$, reaches the chamber in the upper right hand corner, while the flopping of the extremal ray $l_2$, followed by the flopping of the extremal ray  $l_{12}$, reaches the chamber in the lowe right hand corner.  These two chambers, though labeled with the same $\mathbb{P}^3$, correspond to two different minimal models in the KKMR decomposition, and connected  by a flop of a type which we do not specify in (4).

The flopping of the extremal ray $l_2$, followed by the flopping of the extremal ray $l_{32}$ on the 4-fold $Y$, gives $\mathbb{Q}$ when restricted to the strict transform of the zero section.  Similarly, the flopping of the extremal ray $l_3$, followed by the flopping of the extremal ray $l_{23}$ on the 4-fold $Y$, gives the same $\mathbb{Q}$ when restricted to the strict transform of the zero section.  However, The flopping of the extremal ray $l_3$, followed by the flopping of the extremal ray  $l_{23}$, reaches the chamber in the upper left hand corner, while the flopping of the extremal ray $l_2$, followed by the flopping of the extremal ray  $l_{32}$, reaches the chamber in the lowe left hand corner.  These two chambers, though labeled with the same $\mathbb{Q}$, correspond to two different minimal models in the KKMR decomposition, and connected  by a flop of one of the Types labeled as (Others).

\vskip.03in

About (4).

\hskip.1in Mistake: ($E_1$)($E_2$)($E_5$)

\hskip.1in Correction: ($E_1$)($E_2$)(Others)

Note: There is no flop of type ($E_5$).

\newpage

\underline{$n^o\ 24$.}

\vskip.03in

About (1).

Mistake:

\begin{itemize}

\item[] $l_2$: the strict transform of a fiber $W_6 \rightarrow \mathbb{P}^2$

\end{itemize}

Correction:

\begin{itemize}

\item[] $l_2$: the strict transform of a fiber of $p_1:W_6 \rightarrow \mathbb{P}^2$

\end{itemize}

About (2).

The explanation of the symbol $E_1$ on the top row of the table is missing:

\begin{itemize}

\item[] $E_1$: the exceptional divisor of the blow up $W_6 \times_{\mathbb{P}^2}\mathbb{F}^1 \rightarrow W_6$

\end{itemize}

\hskip.1in NO Mistake in the table of the intersection pairings.

\vskip.03in

\underline{$n^o\ 25$.}

\vskip.03in

About (2).

The explanation of the symbols $E_1$ and $E_2$ on the top row of the table is missing:

\begin{itemize}

\item[] $E_1$: the exceptional divisor of the blow up of the line $L_1$

\item[] $E_2$: the exceptional divisor of the blow up of the line $L_2$

\end{itemize}

\hskip.1in NO Mistake in the table of the intersection pairings.

\vskip.03in

\underline{$n^o\ 26$.}

\vskip.03in

About (2).

The explanation of the symbols $E_1$ and $E_2$ on the top row of the table is missing:

\begin{itemize}

\item[] $E_1$: the exceptional divisor of the blow up of the line

\item[] $E_2$: the exceptional divisor of the blow up of the point

\end{itemize}

\vskip.03in

\hskip.1in Mistake:

\begin{center}
\begin{tabular}{|c|c|c|c|c|}
\hline
  &  $E_1$ & $E_2$ & $\pi^*\mathcal{O}_{\mathbb{P}^3}(1)$ & $- K_T$ \\
\hline
$l_{12}$ & 1 & 0 & 0 & -1 \\
\hline
$l_{22}$ & 0 & $\frac{1}{2}$ & 0 & -1 \\
\hline
$l_{32}$ & 0 & 0 & 1 & 4 \\
\hline

\end{tabular}
\end{center}

\vskip.03in

\hskip.1in Correction:

\begin{center}
\begin{tabular}{|c|c|c|c|c|}
\hline
  &  $E_1$ & $E_2$ & $\pi^*\mathcal{O}_{\mathbb{P}^3}(1)$ & $- K_T$ \\
\hline
$l_{12}$ & -1 & 0 & 0 & 1 \\
\hline
$l_{22}$ & 0 & $\frac{1}{2}$ & 0 & -1 \\
\hline
$l_{32}$ & 1 & 0 & 1 & 3 \\
\hline

\end{tabular}
\end{center}

\vskip.03in

\underline{$n^o\ 27$.}

\vskip.03in

\hskip.1in NO Mistake.

\newpage

\underline{$n^o\ 28$.}

\vskip.03in

About (2).

The explanation of the symbol $E_1$ on the top row of the table is missing:

\begin{itemize}

\item[] $E_1$: the exceptional divisor of the map $\mathbb{P}^1 \times \mathbb{F}^1 \rightarrow \mathbb{P}^1 \times \mathbb{P}^2$, i.e., 

$E_1 = \mathbb{P}^1 \times \{\text{the exceptional curve of the blow up }\mathbb{F}^1 \rightarrow \mathbb{P}^2$.

\end{itemize}

\hskip.1in NO Mistake in the table of the intersection pairings.

\vskip.03in

\underline{$n^o\ 29$.}

\vskip.03in

About (2).

We gave the toric description of the Fano 3-fold $T$, and the total space of its canonical bundle $Y$ accordingly in \cite{WGBTMM95}.  We can use this troic description to compute the intersection pairings of the extremal rational curves with the divisors (even though the author made some mistakes in carrying out the toric computatuion in \cite{WGBTMM95}).  However, they can be computed in an elementary way for the flop of type ($E_1$) (resp. of type ($E_5$)) as demonstrated with $\boxed{B_2 = 2}$ $n^o\ 1$ (resp. with $\boxed{B_2 = 2}$ $n^o\ 28$).  We remark that an elementary way to compute the intersection pairings for the flop of type ($E_2$) (resp. of type ($E_3$ or type ($E_4$)) as demonstrated with $\boxed{B_2 = 2}$ $n^o\ 30$ (resp. with $\boxed{B_2 = 2}$ $n^o\ 8$).

\hskip.1in Mistake:

\begin{center}
\begin{tabular}{|c|c|c|c|c|}
\hline
  &  $E_1$ & $E_3$ & $\pi^*\mathcal{O}_{\mathbb{P}^3}(1)$ & $- K_T$ \\
\hline
$l_{13}$ & $- \frac{1}{2}$ & 0 & 0 & $\frac{3}{2}$ \\
\hline
$l_{23}$ & 2 & 0 & 1 & 2 \\
\hline
$l_{33}$ & $- \frac{1}{4}$ & $\frac{1}{2}$ & 0 & $- \frac{1}{4}$ \\
\hline

\end{tabular}
\end{center}

\vskip.03in

\hskip.1in Correction:

\begin{center}
\begin{tabular}{|c|c|c|c|c|}
\hline
  &  $E_1$ & $E_3$ & $\pi^*\mathcal{O}_{\mathbb{P}^3}(1)$ & $- K_T$ \\
\hline
$l_{13}$ & $- \frac{1}{2}$ & 0 & 0 & $\frac{3}{2}$ \\
\hline
$l_{23}$ & 1 & 0 & 1 & 1 \\
\hline
$l_{33}$ & $- \frac{1}{2}$ & 1 & 0 & $- \frac{1}{2}$ \\
\hline

\end{tabular}
\end{center}

\vskip.03in

About (3) the KKMR decomposition.

Note: The flopping of the extremal ray $l_1$, followed by the flopping of the extremal ray $l_{31}$ on the 4-fold $Y$, gives $\mathbb{P}^3$ when restricted to the strict transform of the zero section.  Similarly, the flopping of the extremal ray $l_3$, followed by the flopping of the extremal ray $l_{13}$ on the 4-fold $Y$, gives the same $\mathbb{P}^3$ when restricted to the strict transform of the zero section.  However, The flopping of the extremal ray $l_1$, followed by the flopping of the extremal ray  $l_{31}$, reaches the chamber in the lower right hand corner, while the flopping of the extremal ray $l_3$, followed by the flopping of the extremal ray  $l_{13}$, reaches the chamber in the lowe left hand corner.  These two chambers, though labeled with the same $\mathbb{P}^3$, correspond to two different minimal models in the KKMR decomposition, and connected  by a flop of type ($G$).

\vskip.03in

\underline{$n^o\ 30$.}

\vskip.03in

About (2).

The explanation of the symbols $E_1$ and $E_2$ on the top row of the table is missing:

\begin{itemize}

\item[] $E_1$: the exceptional divisor of the blow up of the strict transform of the line passing through the center of the blowing up $V_7 \rightarrow \mathbb{P}^3$.

\item[] $E_2$: the strict transform of the exceptional divisor of the blow up  $V_7 \rightarrow \mathbb{P}^3$

\end{itemize}

\hskip.1in NO Mistake in the table of the intersection pairings.

\vskip.03in

About (3) the KKMR decomposition.

Note: The flopping of the extremal ray $l_1$, followed by the flopping of the extremal ray $l_{21}$ on the 4-fold $Y$, gives $\mathbb{P}^3$ when restricted to the strict transform of the zero section.  Similarly, the flopping of the extremal ray $l_2$, followed by the flopping of the extremal ray $l_{12}$ on the 4-fold $Y$, gives the same $\mathbb{P}^3$ when restricted to the strict transform of the zero section.  However, The flopping of the extremal ray $l_1$, followed by the flopping of the extremal ray  $l_{21}$, reaches the chamber in the upper right hand corner, while the flopping of the extremal ray $l_2$, followed by the flopping of the extremal ray  $l_{12}$, reaches the chamber in the lowe right hand corner.  These two chambers, though labeled with the same $\mathbb{P}^3$, correspond to two different minimal models in the KKMR decomposition, and connected  by a flop of one of the Types labeled as (Others).

\vskip.03in

About (4).

\hskip.1in Mistake: ($E_1$)($E_2$)($E_5$)

\hskip.1in Correction: ($E_1$)($E_2$)(Others)

Note: There is no flop of type ($E_5$).

\vskip.03in

\underline{$n^o\ 31$.}

\vskip.03in

About (2).

The explanation of the symbol $E$ on the top row of the table is missing:

\begin{itemize}

\item[] $E$: the exceptional divisor of the blow up of the cone over a smooth quadric surface in $\mathbb{P}^3$ with center the vertex

\end{itemize}

\hskip.03in Mistake:

\begin{center}
\begin{tabular}{|c|c|c|c|c|}
\hline
  &  $E$ & $\pi_1^*\mathcal{O}_{\mathbb{P}^1}(1)$ & $\pi_2^*\mathcal{O}_{\mathbb{P}^1}(1)$ & $- K_T$ \\
\hline
$l_{31}$ & 0& 1 & 0 & $-1$ \\
\hline
$l_{32}$ & 0& 0 & 1 & $-1$ \\
\hline
\end{tabular}
\end{center}

\hskip.03in Correction:

\begin{center}
\begin{tabular}{|c|c|c|c|c|}
\hline
  &  $E$ & $\pi_1^*\mathcal{O}_{\mathbb{P}^1}(1)$ & $\pi_2^*\mathcal{O}_{\mathbb{P}^1}(1)$ & $- K_T$ \\
\hline
$l_{31}$ & 0& 1 & 0 & 3 \\
\hline
$l_{32}$ & 0& 0 & 1 & 3 \\
\hline
\end{tabular}
\end{center}

\newpage

$\boxed{B_2 = 4}$

\vskip.03in

\underline{$n^o\ 1$.}

\vskip.03in

About (1).

Mistake:

\begin{itemize}

\item[] $l_i\ (i = 1,2,3,4)$: the ruling of the exceptional divisor of the projection $p_i: T \rightarrow \mathbb{P}^1 \times \mathbb{P}^1 \times \mathbb{P}^1$ (omitting the $i$-th factor) with the center of blow up being a curve of curve of complete intersection of divisors of of tridegrees $(1,0,0), (0,1,0), (0,0,1)$ ... $(E_1$)

\end{itemize}

Correction:

\begin{itemize}

\item[] $l_i\ (i = 1,2,3,4)$: the ruling of the exceptional divisor of the projection $p_i: T \rightarrow \mathbb{P}^1 \times \mathbb{P}^1 \times \mathbb{P}^1$ (omitting the $i$-th factor), which is nothing but the blow up of $\mathbb{P}^1 \times \mathbb{P}^1  \times \mathbb{P}^1 $ with the center of blow up being a curve of curve of complete intersection of divisors of of tridegree $(1,1,1)$ ... $(E_1$)

\end{itemize}

\vskip.03in

\underline{$n^o\ 2$.}

\vskip.03in

About (1).

Mistake:

\begin{itemize}

\item[] $l_4$: the other ruling of $S$ ... ($E_5$)

\item[] $l_5$: the ruling of the exceptional divisor ($\cong \mathbb{P}^1 \times \mathbb{P}^1$) of the blow up of the cone, having $0$ intersection with $\pi_2^*\mathcal{O}_{\mathbb{P}^1}(1)$ ... ($E_1$)

\end{itemize}

Correction:

\begin{itemize}

\item[] $l_4$: the other ruling of $S$ ... ($E_1$)

\item[] $l_5$: the ruling of the exceptional divisor ($\cong \mathbb{P}^1 \times \mathbb{P}^1$) of the blow up of the vertex of the cone, having $0$ intersection with $\pi_2^*\mathcal{O}_{\mathbb{P}^1}(1)$ ... ($E_1$)

\end{itemize}

\vskip.03in

About (2).

The explanation of the symbol $E_5$ on the top row of the table is missing:

\begin{itemize}

\item[] $E_5$: the exceptional divisor of the blow up of the vertex of the cone over a smooth quadric surface in $\mathbb{P}^3$

\end{itemize}

Mistake and Correction:

\newpage 

\begin{center}
\begin{tabular}{|c|c|c|c|c|c|}
\hline
  &  $\pi_1^*\mathcal{O}_{\mathbb{P}^1}(1)$ &  $\pi_1^*\mathcal{O}_{\mathbb{P}^1}(1)$ & $S$ & $E_5$ & $- K_T$ \\
\hline
$l_1$ & 0& 0 & 1 & 0$-1$ & 1 \\
\hline
$l_2$ & 0& 0 & 0 & 1 & 1 \\
\hline
$l_3$ & 1& 0 & $-1$ & 2 & 1 \\
\hline
$l_4$ & 0& 1 & $-1$ & 2 & 1 \\
\hline
$l_5$ & 1& 0 & 0 & 0 & 1 \\
\hline
$l_6$ & 0& 1 & 0 & 0 & 1 \\
\hline
$l_{11}$ & 0& 0 & $-1$ & 1 & $-1$ \\
\hline
$l_{21}$ & 0& 0 & 1 & 0 & 2 \\
\hline
$l_{31}$ & 1 & 0 & 1 & 0 & 3 \\
\hline
$l_{41}$ & 0& 1 & 1 & 0 & 3 \\
\hline
$l_{51}$ & 1& 0 & 0 & 0 & 1 \\
\hline
$l_{61}$ & 0& 1 & 0 & 0 & 1 \\
\hline
$l_{12}$ & 0& 0 & 1 & 0 & 2 \\
\hline
$l_{22}$ & 0& 0 & 0 & $-1$ & $-1$ \\
\hline
$l_{32}$ & 1 & 0 & $-1$ & $-2$ & $-3$ \\
\hline
$l_{42}$ & 0& 1 & $-1$ & $-2$ & $-3$ \\
\hline
$l_{52}$ & 1& 0 & 0 & $-2$ & $-1$ \\
\hline
$l_{62}$ & 0& 1 & 0 & $-2$ & $-1$ \\
\hline
$l_{13}$ & 1& 0 & 0 & 1 & 2 \\
\hline
$l_{23}$ & 0& 0 & 0 & 1 & 1 \\
\hline
$l_{33}$ & $-1$ & 0 & 1 & $-2$ & $-1$ \\
\hline
$l_{43}$ & $-1$ & 1 & 0 & 0 & 0 \\
\hline
$l_{53}$ & 1& 0 & 0 & 0 & 1 \\
\hline
$l_{63}$ & 0& 1 & 0 & 0 & 1 \\
\hline
$l_{14}$ & 0 & 1 & 0 & 1 & 2 \\
\hline
$l_{24}$ & 0& 0 & 0 & 1 & 1 \\
\hline
$l_{34}$ & 1 & $-1$ & 0 & 0 & $-1$ \\
\hline
$l_{44}$ & 0 & $-1$ & 1 & $-2$ & $-1$ \\
\hline
$l_{54}$ & 1& 0 & 0 & 0 & 1 \\
\hline
$l_{64}$ & 0& 1 & 0 & 0 & 1 \\
\hline
$l_{15}$ & 0& 0 & 1 & $-1$ & 1 \\
\hline
$l_{25}$ & 1& 0 & 0 & 1 & 2 \\
\hline
$l_{35}$ & 1 & 0 & $-1$ & 2 & 1 \\
\hline
$l_{45}$ & 0 & 1 & $-1$ & 2 & 1 \\
\hline
$l_{55}$ & $-1$ & 0 & 0 & 0 & $-1$ \\
\hline
$l_{65}$ & $-1$ & 1 & 0 & 0 & 0 \\
\hline
$l_{16}$ & 0& 0 & 1 & $-1$ & 1 \\
\hline
$l_{26}$ & 0 & 1 & 0 & 1 & 2 \\
\hline
$l_{36}$ & 1 & 0 & $-1$ & 2 & 1 \\
\hline
$l_{46}$ & 0 & 1 & $-1$ & 2 & 1 \\
\hline
$l_{56}$ & 1& $-1$& 0 & 0 & 0 \\
\hline
$l_{66}$ & 0 & $-1$ & 0 & 0 & $-1$ \\
\hline
\end{tabular}
\end{center}

\newpage

\begin{center}
\begin{tabular}{|c|c|c|c|c|c|}
\hline
  &  $\pi_1^*\mathcal{O}_{\mathbb{P}^1}(1)$ &  $\pi_1^*\mathcal{O}_{\mathbb{P}^1}(1)$ & $S$ & $E_5$ & $- K_T$ \\
\hline
$l_1$ & 0& 0 & 1 & 0 & 1 \\
\hline
$l_2$ & 0& 0 & 0 & 1 & 1 \\
\hline
$l_3$ & 1& 0 & $-1$ & 0 & 1 \\
\hline
$l_4$ & 0& 1 & $-1$ & 0 & 1 \\
\hline
$l_5$ & 1& 0 & 0 & $-1$ & 1 \\
\hline
$l_6$ & 0& 1 & 0 & $-1$ & 1 \\
\hline
$l_{11}$ & 0& 0 & $-1$ & 0 & $-1$ \\
\hline
$l_{21}$ & 0& 0 & 1 & 1 & 2 \\
\hline
$l_{31}$ & 1 & 0 & 1 & 0 & 3 \\
\hline
$l_{41}$ & 0& 1 & 1 & 0 & 3 \\
\hline
$l_{51}$ & 1& 0 & 2 & $-1$ & 1 \\
\hline
$l_{61}$ & 0& 1 & 2 & $-1$ & 1 \\
\hline
$l_{12}$ & 0& 0 & 1 & 1 & 2 \\
\hline
$l_{22}$ & 0& 0 & 0 & $-1$ & $-1$ \\
\hline
$l_{32}$ & 1 & 0 & $-1$ & 0 & 0 \\
\hline
$l_{42}$ & 0& 1 & $-1$ & 0 & 0 \\
\hline
$l_{52}$ & 1& 0 & 0 & 1 & 3 \\
\hline
$l_{62}$ & 0& 1 & 0 & 1 & 3 \\
\hline
$l_{13}$ & 1& 0 & 0 & 0 & 2 \\
\hline
$l_{23}$ & 0& 0 & 0 & 1 & 1 \\
\hline
$l_{33}$ & $-1$ & 0 & 1 & 0 & $-1$ \\
\hline
$l_{43}$ & $-1$ & 1 & 0 & 0 & 0 \\
\hline
$l_{53}$ & 1& 0 & 0 & $-1$ & 1 \\
\hline
$l_{63}$ & 0& 1 & 0 & $-1$ & 1 \\
\hline
$l_{14}$ & 0 & 1 & 0 & 1 & 2 \\
\hline
$l_{24}$ & 0& 0 & 0 & 1 & 1 \\
\hline
$l_{34}$ & 1 & $-1$ & 0 & 0 & 0 \\
\hline
$l_{44}$ & 0 & $-1$ & 1 & $-2$ & $-1$ \\
\hline
$l_{54}$ & 1& 0 & 0 & 0 & 1 \\
\hline
$l_{64}$ & 0& 1 & 0 & 0 & 1 \\
\hline
$l_{15}$ & 0& 0 & 1 & 0 & 1 \\
\hline
$l_{25}$ & 1& 0 & 0 & 0 & 2 \\
\hline
$l_{35}$ & 1 & 0 & $-1$ & 0 & 1 \\
\hline
$l_{45}$ & 0 & 1 & $-1$ & 0 & 1 \\
\hline
$l_{55}$ & $-1$ & 0 & 0 & 1 & $-1$ \\
\hline
$l_{65}$ & $-1$ & 1 & 0 & 0 & 0 \\
\hline
$l_{16}$ & 0& 0 & 1 & 0 & 1 \\
\hline
$l_{26}$ & 0 & 1 & 0 & 0 & 2 \\
\hline
$l_{36}$ & 1 & 0 & $-1$ & 0 & 1 \\
\hline
$l_{46}$ & 0 & 1 & $-1$ & 0 & 1 \\
\hline
$l_{56}$ & 1& 2& 0 & 0 & 0 \\
\hline
$l_{66}$ & 0 & 2 & 0 & 1 & $-1$ \\
\hline

\end{tabular}
\end{center}

\vskip.03in

\underline{$n^o\ 3$.}

\vskip.03in

See the complete discussion in Section 4.

\vskip.03in

\underline{$n^o\ 4$.}

\vskip.03in

About the picture.

The morphism $T \rightarrow Q$ should be labeled as $\pi$.

\vskip.03in

About (2).

The explanation of the symbols $E_1$, $E_2$ and $E_5$ on the top row of the table is missing:

\begin{itemize}

\item[] $E_1$: the strict transform of the exceptional divisor of the blow up of the point $p$

\item[] $E_2$: the strict transform of the exceptional divisor of the blow up of the point $q$

\item[] $E_5$: the exceptional divisor of the blow of the strict transform of the conic

\end{itemize}

\hskip.1in NO Mistake in the table of the intersection pairings.

\vskip.03in

About (4).

\hskip.1in Mistake: ($E_1$)($E_2$)($E_5$)

\hskip.1in Correction: ($E_1$)($E_2$), (Others) associated with $\boxed{B_2=3}$ $n^o\ 18$, $n^o\ 30$

Note: There is no flop of type ($E_5$).

\vskip.03in

\underline{$n^o\ 5$.}

\vskip.03in

About (2).

The explanation of the symbols $E_1$, $E_2$, $\pi_1$, $\pi_2$ on the top row of the table is missing:

\begin{itemize}

\item[] $E_1$: the exceptional divisor of the blow up of the curve of bidegree $(0,1)$

\item[] $E_2$: the exceptional divisor of the blow up of the curve of bidegree $(2,1)$

\item[] $\pi_1: \mathbb{P}^1 \times \mathbb{P}^2 \rightarrow \mathbb{P}^1$

\item[] $\pi_2: \mathbb{P}^1 \times \mathbb{P}^2 \rightarrow \mathbb{P}^2$

\end{itemize}

\vskip.03in

\hskip.1in Mistake:

\begin{center}
\begin{tabular}{|c|c|c|c|c|c|}
\hline
  &  $E_1$ & $E_2$ & $\pi_1^*\mathcal{O}_{\mathbb{P}^1}(1)$ & $\pi_2^*\mathcal{O}_{\mathbb{P}^2}(1)$ &$- K_T$ \\

\hline
$l_{24}$ & 0 & 1 & 0 & 1 & $-1$ \\
\hline
$l_{15}$ & $-1$ & 0 & 0 & $-1$ & $-2$ \\
\hline
$l_{25}$ & 1 & 0 & 0 & 1 & 2 \\
\hline
$l_{35}$ & $-1$ & 0 & 1 & $-1$ & 0 \\
\hline
$l_{45}$ & $-2$ & 0 & 0 & $-2$ & $-4$ \\
\hline
$l_{55}$ & $-1$ & $-1$ & 0 & $-1$ & $-1$ \\
\hline

\end{tabular}
\end{center}

\vskip.03in

\hskip.1in Correction:

\begin{center}
\begin{tabular}{|c|c|c|c|c|c|}
\hline
  &  $E_1$ & $E_2$ & $\pi_1^*\mathcal{O}_{\mathbb{P}^1}(1)$ & $\pi_2^*\mathcal{O}_{\mathbb{P}^2}(1)$ &$- K_T$ \\

\hline
$l_{24}$ & 0 & 1 & 0 & 1 & 2 \\
\hline
$l_{15}$ & 1 & 2 & 0 & 2 & 3 \\
\hline
$l_{25}$ & 1 & 0 & 0 & 1 & 2 \\
\hline
$l_{35}$ & 0 & 1 & 1 & 0 & 1 \\
\hline
$l_{45}$ & 0 & 2 & 0 & 1 & 1 \\
\hline
$l_{55}$ & $-1$ & $-1$ & 0 & $-1$ & $-1$ \\
\hline

\end{tabular}
\end{center}

\newpage

\underline{$n^o\ 6$.}

\vskip.03in

About (1).

The explanation of the symbols $p_{23}$, $p_{13}$, $p_{12}$ is missing:

\begin{itemize}

\item[] $p_{23}:\mathbb{P}_1^1 \times \mathbb{P}_2^1 \times \mathbb{P}_3^1 \rightarrow \mathbb{P}_2^1 \times \mathbb{P}_3^1$: the projection onto the 2nd and 3rd factors.

\item[] $p_{13}:\mathbb{P}_1^1 \times \mathbb{P}_2^1 \times \mathbb{P}_3^1 \rightarrow \mathbb{P}_1^1 \times \mathbb{P}_3^1$: the projection onto the 1st and 3rd factors.

\item[] $p_{12}:\mathbb{P}_1^1 \times \mathbb{P}_2^1 \times \mathbb{P}_3^1 \rightarrow \mathbb{P}_1^1 \times \mathbb{P}_2^1$: the projection onto the 1st and 2nd factors.

\end{itemize}

\vskip.03in

About (2).

The explanation of the symbols $\pi_1$, $\pi_2$, $\pi_3$, $E_4$ is missing:

\begin{itemize}

\item[] $\pi_1: T \rightarrow \mathbb{P}_1^1$: the morphism onto 1st factor.

\item[] $\pi_2: T \rightarrow \mathbb{P}_2^1$: the morphism onto 2nd factor.

\item[] $\pi_3: T \rightarrow \mathbb{P}_3^1$: the morphism onto 3rd factor.

\item[] $E_4$: the exceptional divisor of the blow up of the tridiagonal curve $C$

\end{itemize}

\hskip.1in NO Mistake in the table of the intersection pairings.

\vskip.03in

\underline{$n^o\ 7$.}

\vskip.03in

About (2).

The explanation of the symbols $E_1$, $E_2$ is missing:

\begin{itemize}

\item[] $E_1$: the exceptional divisor of the blow up of the curve $C_1 \subset W_6 \subset \mathbb{P}^2 \times \mathbb{P}^2$ of bidegree $(0,1)$

\item[] $E_2$: the exceptional divisor of the blow up of the curve $C_2 \subset W_6 \subset \mathbb{P}^2 \times \mathbb{P}^2$ of bidegree $(1,0)$

\end{itemize}

\hskip.1in NO Mistake in the table of the intersection pairings.

\vskip.03in

\underline{$n^o\ 8$.}

\vskip.03in

About (2).

The explanation of the symbols $E_1$, $\pi_1$, $\pi_2$, $\pi_3$ is missing:

\begin{itemize}

\item[] $E_1$: the exceptional divisor of the blow up of the curve $C$ of tridegree $(0,1,1)$

\item[] $\pi_1: T (\rightarrow \mathbb{P}_1^1 \times \mathbb{P}_2^1 \times \mathbb{P}_3^1) \rightarrow \mathbb{P}_1^1$ the projection onto the 1st factor

\item[] $\pi_2: T (\rightarrow \mathbb{P}_1^1 \times \mathbb{P}_2^1 \times \mathbb{P}_3^1) \rightarrow \mathbb{P}_2^1$ the projection onto the 1st factor
\
\item[] $\pi_3: T (\rightarrow \mathbb{P}_1^1 \times \mathbb{P}_2^1 \times \mathbb{P}_3^1) \rightarrow \mathbb{P}_3^1$ the projection onto the 1st factor

\end{itemize}

\hskip.1in NO Mistake in the table of the intersection pairings.

\vskip.03in

\underline{$n^o\ 9$.}

\vskip.03in

About (2).

The explanation of the symbols $E_1$, $E_2$, $E_4$ and $\pi$ on the top row of the table is missing:

\begin{itemize}

\item[] $E_1$: the strict transform of the exceptional divisor of the blow up of the line $L_1$

\item[] $E_2$: the strict transform of the exceptional divisor of the blow up of the line $L_2$

\item[] $E_4$: the exceptional divisor of the blow of the exceptional line

\item[] $\pi: T (\rightarrow Y) \rightarrow \mathbb{P}^3$ the natural morphism

\end{itemize}

\hskip.1in NO Mistake in the table of the intersection pairings.

\newpage

About (4).

\hskip.1in Mistake: ($E_1$)($E_2$)($E_5$)

\hskip.1in Correction: ($E_1$)($E_2$) and (Others) associated with $\boxed{B_2 = 3}$ $n^o\ 30$

Note: There is no flop of type ($E_5$).

\vskip.03in

\underline{$n^o\ 10$.}

\vskip.03in

$T = \mathbb{P}^1 \times S_7$ where $S_7$ is the Del Pezzo surface of degree $(K_{S_7})^2 = 7$, i.e., the blow up of $\mathbb{P}^2$ with center 2 (general) points

\vskip.03in

About the picture.

\includegraphics[width=4in]{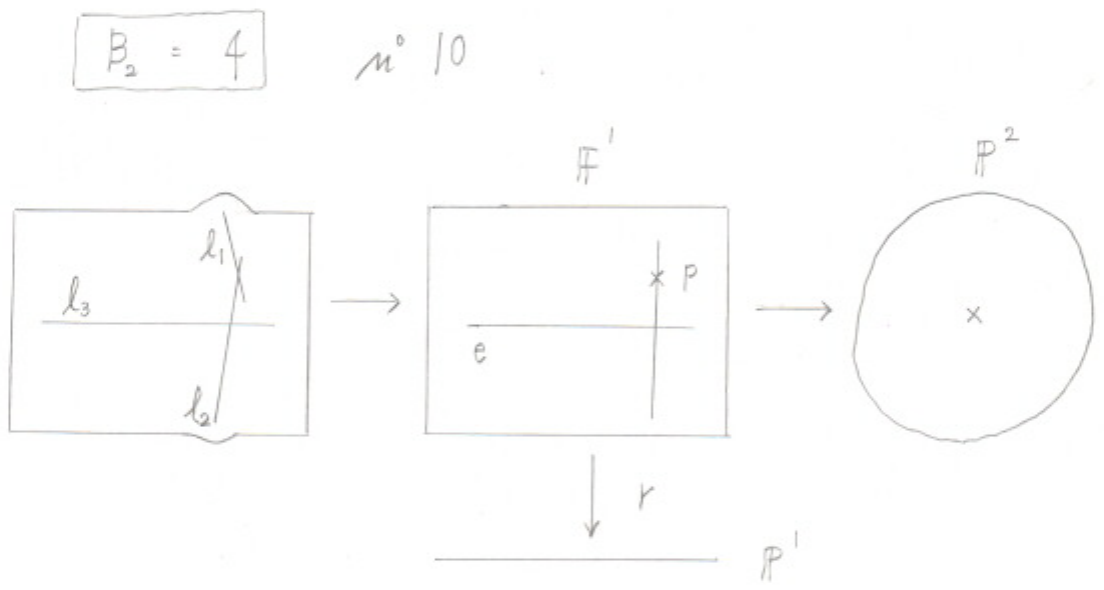}

\vskip.03in

About (1).

\vskip.03in

Extremal Rational Curves for $S_7$.

\begin{itemize}

\item[] $l_1$: the exceptional divisor for the blow up of a point $P \in \mathbb{F}^1$, which is not contained in the unique minimal section $e$ of the ruling $r:\mathbb{F}^1 \rightarrow \mathbb{P}^1$.

\item[] $l_2$: the strict transform of the ruling $r^{-1}(r(P))$

\item[] $l_3$: the strict transform of $e$

\end{itemize}

\vskip.03in

Extremal Rational Curves for $T = \mathbb{P}^1 \times S_7$.  

Let $i: S_7 \overset{\sim}\rightarrow \{0\} \times S_7 \subset \mathbb{P}^1 \times S_7 = T$ be the inclusion.

\begin{itemize}

\item[] $i_*(l_1)$: ... ($E_1$)

\item[] $i_*(l_2)$: ... ($E_1$)

\item[] $i_*(l_3)$: ... ($E_1$)

\item[] $l_4$: the fiber of the projection $\mathbb{P}^1 \times S_7 \rightarrow S_7$

\end{itemize}

\vskip.03in

About (2) and (3). 

\hskip.1in Left to the reader as an exercise.

\vskip.03in

About (4).

\hskip.1in ($E_1$)($F$)

\vskip.03in

About (5).

\hskip.1in $WG_T = A_1$.

\vskip.03in

\underline{$n^o\ 11$.}

\vskip.03in

About (2).

The explanation of the symbols $E_1$, $E_2$, $E_4$ and $p_2^*(ruling)$ on the top row of the table is missing:

\begin{itemize}

\item[] $E_1$: the exceptional divisor of the blow up of $t \times e$

\item[] $E_2$: the strict transform of $\mathbb{P}^1 \times e$

\item[] $E_4$: the strict transform of $\mathbb{P}^1 \times e$

\item[] $p_2: \mathbb{P}^1 \times \mathbb{F}_1 \rightarrow \mathbb{F}_1$, , $r: \mathbb{F}_1 \rightarrow \mathbb{P}^1$ the morphism giving the structure of a ruled surface, $ruling = r^{-1}(r(P))\ P \in \mathbb{P}^1$

\end{itemize}

\hskip.1in NO Mistake in the table of the intersection pairings.

\vskip.03in

About (4).

\hskip.1in Mistake: ($E_1$)($E_2$)($E_5$)($F$)

\hskip.1in Correction: ($E_1$)($F$)

Note: There is no flop of type ($E_2$) or ($E_5$).

\vskip.03in

\vskip.03in

\underline{$n^o\ 12$.}

\vskip.03in

About the picture.

The most upper rightarrow indicating the morphism going from $T$ to $\mathbb{P}^3$ \text{should be labeled as $\pi$.}

\vskip.03in

About (2).

The explanation of the symbols $E_1$, $E_2$, $E_3$ on the top row of the table is missing:

\begin{itemize}

\item[] $E_1$: the exceptional divisor of the blow up of one exceptional line $L_1$

\item[] $E_2$: the exceptional divisor of the blow up of another exceptional line $L_2$

\item[] $E_4$: the strict transform of the exceptional divisor of the blow up $Y \rightarrow \mathbb{P}^3$

\end{itemize}

\hskip.1in NO Mistake in the table of the intersection pairings.

About (4).

\hskip.1in Mistake: ($E_1$)($E_2$)($E_5$)($F$)

\hskip.1in Correction: ($E_1$)($E_2$)($F$) and (Others) associated with $\boxed{B_2 = 3}$ $n^o\ 30$

Note: There is no flop of type ($E_5$).

\newpage

$\boxed{B_2 = 5}$

\vskip.03in

\underline{$n^o\ 1$.}

\vskip.03in

See the detailed discussion in Section 3.

\vskip.03in

\underline{$n^o\ 2$.}

\vskip.03in

About (2).

The explanation of the symbols $E_1$, $E_2$, $E_3$, $E_5$, $\pi$ on the top row of the table is missing:

\begin{itemize}

\item[] $E_1$: the exceptional divisor of the blow up of one exceptional line $\epsilon_1$

\item[] $E_2$: the exceptional divisor of the blow up of another exceptional line $\epsilon_2$

\item[] $E_3$: the strict transform of the exceptional divisor of the blow up of the line $L_2$

\item[] $E_5$: the exceptional divisor of the blow up of the line $L_1$

\item[] $\pi: T (\rightarrow Y) \rightarrow \mathbb{P}^3$ is the natural morphism

\end{itemize}

\hskip.03in Mistake:

\begin{center}
\begin{tabular}{|c|c|c|c|c|c|c|}
\hline
  &  $E_1$ & $E_2$ & $E_3$ & $E_5$ & $\pi^*\mathcal{O}_{\mathbb{P}^3}(1)$ & $- K_T$ \\
\hline
$l_{36}$ & 0& 0 & $-1$ & 0 & 1 & 1 \\
\hline
\end{tabular}
\end{center}

\hskip.03in Correction:

\begin{center}
\begin{tabular}{|c|c|c|c|c|c|c|}
\hline
  &  $E_1$ & $E_2$ & $E_3$ & $E_5$ & $\pi^*\mathcal{O}_{\mathbb{P}^3}(1)$ & $- K_T$ \\
\hline
$l_{36}$ & 0 & 0 & $-1$ & 0 & 0 & 1 \\
\hline
\end{tabular}
\end{center}

\vskip.03in

About (4).

\hskip.1in Mistake: ($E_1$)($E_2$)($E_5$)($F$)

\hskip.1in Correction: ($E_1$)($E_2$)($F$) and (Others) associated with $\boxed{B_2 = 4}$ $n^o\ 12$

Note: There is no flop of type ($E_5$).

\vskip.03in

\underline{$n^o\ 3$.}

\vskip.03in

$T = \mathbb{P}^1 \times S_6$ where $S_6$ is the Del Pezzo surface of degree $(K_{S_7})^2 = 6$, i.e., the blow up of $\mathbb{P}^2$ with center 3 (general) points.

The analysis of the extremal rays, the KKMR deocmposition, and the types of flops, for the associated 4-fold is easily reduced to that of the Del Pezzo surface $S_6$.

See the analysis for the case $\boxed{B_2 = 4}$ $n^o\ 10$: $T = \mathbb{P}^1 \times S_7$.

\newpage

\begin{center}{\bf Table of the Weyl Groups}\end{center}

\vskip.2in

\rotatebox{-270}{

\centering
\scalebox{0.9}{

\begin{tabular}{|c|c|c|c|c|c|c|c|c|c|c|c|}
\hline
 $B_2$ & $\{e\}$ & $A_1$ & $A_1 \times A_1$ & $A_2$ & $A_1 \times A_2$ & $A_3$ & $A_4$ & $D_5$ & $E_6$ & $E_7$ & $E_8$ \\
\hline
2 & $\begin{array}{l}
 1,3,4 \\\\
5, 7,8 \\
9,10,11 \\
13,14, 15 \\
16,17,18 \\
19,20,22 \\
23,24, 25\\
26,27, 28 \\
29, 30, 31\\
33,34,35\\
36 \\
\end{array}$ & $\begin{array}{l}
2,6 \\
12,21 \\
32 \\
\end{array}$ & & & & & & & & &   \\
\hline
3 & $\begin{array}{l}
2,4,5 \\
6,8, 11 \\
12,14,15 \\
16,18,21\\
22,23,24 \\
26,28,29 \\
30 \\
\end{array}$ & $\begin{array}{l}
3,7 \\
9,10 \\
13,17\\
19, 20 \\
25, 31 \\
\end{array}$ & & $\begin{array}{l}
1 \\
27\\
\end{array}$ 

& & & & & & &   \\
\hline
4 & $\begin{array}{l}
5,9,11 \\
\end{array}$ & $\begin{array}{l}
3,4 \\
7,8 \\
10,12 \\
13 \\
\end{array}$ & $\begin{array}{l}
2\\
\end{array}$ & $\begin{array}{l}
6 \\
\end{array}$
& & $\begin{array}{l}
1 \\
\end{array}$ &
& & & &   \\
\hline

5 &  & $\begin{array}{l}
2 \\
\end{array}$ &  & $\begin{array}{l}
1 \\
\end{array}$
& $\begin{array}{l}
3 \\
\end{array}$ & &
& & & &   \\
\hline

6 &  &  &  & 
&  & & $\begin{array}{l}
\mathbb{P}^1 \times S_5 \\
\end{array}$
& & & &   \\
\hline

7 &  &  &  & 
&  & & & $\begin{array}{l}
\mathbb{P}^1 \times S_4 \\
\end{array}$
& & &   \\
\hline

8 &  &  &  & 
&  & & & & $\begin{array}{l}
\mathbb{P}^1 \times S_3 \\
\end{array}$
& &   \\
\hline

9 &  &  &  & 
&  & & & & & $\begin{array}{l}
\mathbb{P}^1 \times S_2 \\
\end{array}$
&   \\
\hline

10 &  &  &  & 
&  & & & & & & $\begin{array}{l}
\mathbb{P}^1 \times S_1 \\
\end{array}$
\\
\hline
\end{tabular}

}}

\end{document}